\UseRawInputEncoding
\documentclass[12pt]{article}
\usepackage[colorlinks=true,linkcolor=blue,urlcolor=blue,citecolor=blue]{hyperref}
\usepackage{pgfplots}
\usepackage{mathrsfs}
\usepackage{amssymb}
\usepackage{amsmath}
\usepackage{amsfonts}
\usepackage{amstext}
\usepackage{amsthm}
\usepackage{tikz}
\usepackage{graphicx}
\usepackage{subfig}
\usepackage{color}
\usepackage{indentfirst,latexsym,bm}
\usepackage{mathrsfs}
\usepackage{cite}
\usepackage[all]{xy}
\usepackage{newlfont}
\usepackage{enumerate}
\numberwithin{equation}{section}
\theoremstyle{remark}\newtheorem*{remark*}{Remark}\theoremstyle{plain}
\newtheorem*{corollary*}{Corollary}
\newtheorem*{lemma*}{Lemma}
\newtheorem*{proposition*}{Proposition}
\newtheorem*{theorem*}{Theorem}
\allowdisplaybreaks[2]
\newtheorem{theorem}{Theorem}[section]
\newtheorem{lemma}[theorem]{Lemma}
\newtheorem{proposition}[theorem]{Proposition}
\newtheorem{corollary}[theorem]{Corollary}
\theoremstyle{remark}\newtheorem{remark}[theorem]{Remark}
\newcommand{\R}{\mathbb R}
\newcommand{\dd}{\,\mathrm d}
\newcommand{\tr}{\operatorname{tr}}
\newcommand{\diver}{\operatorname{div}}
\newcommand{\supp}{\operatorname{supp}}
\newcommand{\calD}{\mathcal D}

\usepackage{needspace}
\usetikzlibrary{arrows,positioning}

\begin{document}
\title{The P\'{o}lya--Szeg\H{o} buckling conjecture\\in dimensions two through forty-seven
\thanks{Research supported by NNSF of China (No. 12371110 and No. 12301133).}}

\author{
Guowei Dai\thanks{\protect\raggedright School of Mathematical Sciences,
Dalian University of Technology, Dalian 116024, P.R. China.
E-mail: \mbox{daiguowei@dlut.edu.cn}.}
\and
Yingxin Sun\thanks{\protect\raggedright School of Mathematical Sciences,
Dalian University of Technology, Dalian 116024, P.R. China.
E-mail: \mbox{sunyingxin2023@mail.dlut.edu.cn}.}
\and
Yong Zhang\thanks{\protect\raggedright School of Mathematical Sciences,
Jiangsu University, Zhenjiang 212013, P.R. China.
E-mail: \mbox{zhangyong@ujs.edu.cn}.}
}
\date{}

\maketitle
\renewcommand{\abstractname}{Abstract}

\begin{abstract}
We prove the P\'olya--Szeg\H{o} buckling conjecture in dimensions
$2\le n\le47$: among open sets of prescribed finite volume, a ball
minimizes the first clamped buckling eigenvalue. More precisely, we
establish the sharp comparison
$\Lambda_1\left(\Omega\right)\ge
\left(j_{n/2,1}/j_{n/2-1,1}\right)^2
\lambda_1^{\mathrm D}\left(\Omega\right)$,
where $\lambda_1^{\mathrm D}$ is the first Dirichlet eigenvalue and
$j_{\nu,1}$ is the first positive zero of the Bessel function $J_\nu$ of the first kind.
The proof extends a radial calibration to all Hessian directions;
positive coefficient expansions and logarithmically convex Laplace
integrals control the resulting matrix. A Lipschitz construction
handles sign changes and arbitrary boundaries. Equality in the volume
bound holds precisely for balls modulo sets of zero $H^2$ capacity. We prove that
this calibration fails in every dimension $n\ge48$.
\end{abstract}

\section{Introduction}

In 1951, P\'olya and Szeg\H{o} conjectured that, among clamped
plates of prescribed area subjected to uniform compression in the plane of the plate,
the disk has the least critical buckling load \cite{PS}. Its
higher-dimensional form asks whether the ball minimizes the first
clamped buckling eigenvalue at fixed volume. This is a fourth-order
counterpart of the Faber--Krahn theorem for the Dirichlet Laplacian.
The question tests how much geometric information can be retained
when both the bending energy and the gradient energy enter a
variational quotient.

In the planar model, a buckling eigenfunction describes the transverse
deflection of a homogeneous thin plate, and the critical compressive
load is proportional to the first eigenvalue \cite[Section 1]{AL}.
A trial deflection gives an upper bound immediately; a lower bound
must control every admissible deflection. A sharp bound depending
only on area therefore gives a shape-independent threshold below
which this idealized plate remains stable. In higher dimensions the
same problem belongs to spectral geometry and shape optimization:
it connects a fourth-order variational problem with sharp
second-order spectral inequalities and the geometry of Sobolev spaces.

Let $\Omega\subset\mathbb R^n$ be open, with
$0<\left|\Omega\right|<\infty$, where $\left|\Omega\right|$ denotes Lebesgue measure.
For $j=1,2$, the Sobolev space $H^j\left(\Omega\right)$ consists of functions
whose weak derivatives through order $j$ are square integrable, and
$H_0^j\left(\Omega\right)$ is the closure of the smooth, compactly supported
functions $C_c^\infty\left(\Omega\right)$ in this space. Define
\begin{align}
 \Lambda_1\left(\Omega\right)&=
 \inf_{0\ne U\in H_0^2\left(\Omega\right)}
 \frac{\int_\Omega \left(\Delta U\right)^2\dd x}
      {\int_\Omega \left|\nabla U\right|^2\dd x},\label{eq:buckling}\\
 \lambda_1^{\mathrm D}\left(\Omega\right)&=
 \inf_{0\ne v\in H_0^1\left(\Omega\right)}
 \frac{\int_\Omega \left|\nabla v\right|^2\dd x}
      {\int_\Omega v^2\dd x}.\notag
\end{align}
Here $\nabla$ and $\Delta$ are the Euclidean gradient and Laplacian.
On a smooth bounded domain, the first quotient gives the eigenvalue
problem
\[
 \Delta^2U+\Lambda\Delta U=0\quad\text{in }\Omega,
 \qquad U=\partial_{\mathbf n}U=0\quad\text{on }\partial\Omega,
\]
where $\mathbf n$ is the outward unit normal. These are the clamped
boundary conditions. The closure definition in \eqref{eq:buckling}
continues to make sense for arbitrary open sets.

Write $B_R\left(x_0\right)$ for the open ball of radius $R$ centered at $x_0$,
and $B_R=B_R\left(0\right)$. Let $J_\eta$ denote the Bessel function of the
first kind of order $\eta$, and $j_{\eta,l}$ its $l$th positive zero.
The classical ball formulas are
\[
 \lambda_1^{\mathrm D}\left(B_R\right)=\frac{j_{n/2-1,1}^2}{R^2},
 \qquad \Lambda_1\left(B_R\right)=\frac{j_{n/2,1}^2}{R^2};
\]
see \cite[Sections 1 and 4]{AL} and
\cite[Theorem 3.2 and Section 6]{buckeq:DNT}.
Writing $\omega_n=\left|B_1\right|$, the volume comparison takes the explicit form
\[
 \left|\Omega\right|^{2/n}\Lambda_1\left(\Omega\right)
 \ge \omega_n^{2/n}j_{n/2,1}^2,
 \qquad \left|B_R\right|=\left|\Omega\right|.
\]
For smooth bounded domains the conjectured equality case is a ball.
For arbitrary open sets the precise equality statement involves
Sobolev capacity, as explained below. The exponent $2/n$ is forced by
the dilation law
$\Lambda_1\left(t\Omega\right)=t^{-2}\Lambda_1\left(\Omega\right)$ for $t>0$;
the Dirichlet eigenvalue obeys the same law.

The related Rayleigh conjecture concerns a different quotient,
\[
 \Gamma_1\left(\Omega\right)=\inf_{0\ne U\in H_0^2\left(\Omega\right)}
 \frac{\int_\Omega\left(\Delta U\right)^2\dd x}{\int_\Omega U^2\dd x},
\]
which describes the vibration of a clamped plate.
Talenti's comparison method \cite{Talenti} was a decisive advance:
it treated the positive and negative regions separately and reduced
the fourth-order problem to radial comparison problems on two balls.
Nadirashvili \cite{Nadirashvili} recovered the sharp planar bound by
retaining additional coupling between the rearranged problems.
Ashbaugh and Benguria \cite{AB} gave a detailed Bessel-function
analysis that also proved the three-dimensional vibration inequality.
These results established that sign changes need not prevent a sharp
fourth-order isoperimetric theorem. They do not, however, establish
the corresponding buckling inequality: replacing $\int U^2\dd x$ by
$\int\left|\nabla U\right|^2\dd x$ changes the comparison problem.

For buckling, Szeg\H{o}'s pioneering symmetrization argument
\cite{Szego50,Szego58} proves the planar result under an additional
sign assumption on a first eigenfunction. That assumption is
unavailable on general domains; the examples and the role of
admissibility are discussed in \cite[Section 8]{AL} and \cite{ABL}.
In particular, taking an absolute value does not generally preserve
$H^2$, and ordinary decreasing rearrangement does not supply the
required control of $\int\left(\Delta U\right)^2\dd x$. The obstacle is therefore
not simply a missing maximum-principle argument: both the clamped
function space and the two energies in the quotient must be retained.

A second approach compares buckling and membrane spectra.
Extending the observation of Bramble and Payne \cite{BP},
Ashbaugh and Laugesen \cite{AL} obtained, on bounded domains,
\[
 \Lambda_1\left(\Omega\right)\ge c_n\Lambda_1\left(B_R\right),\qquad
 c_n=2^{2/n}\frac{j_{n/2-1,1}^2}{j_{n/2,1}^2}<1,
 \qquad c_n\longrightarrow1\quad\left(n\longrightarrow\infty\right).
\]
Their result provides an explicit bound in every dimension and is
asymptotically sharp. It also identifies a concrete limitation of that
route: the passage through the second Dirichlet eigenvalue and its
two-ball volume bound leaves a loss at each fixed dimension.
Ashbaugh, Benguria and Laugesen \cite{ABL} explain why the analogous
two-ball reduction is sharp for vibration in dimensions two and three
but gives only a partial result for buckling. This is a limitation of
those reductions, rather than an obstruction to the conjecture itself.

Shape optimization offers a complementary approach. Ashbaugh and
Bucur \cite{AshbaughBucur} developed planar existence and relaxation
methods and examined the Willms--Weinberger characterization of smooth
optimizers. Stollenwerk \cite{Stollenwerk} proved existence of an
optimal open set in every dimension, allowing unbounded competitors.
Stollenwerk and Wagner \cite{SW} proved that a smooth minimizing domain
with connected boundary is a ball, using the second domain variation.
More recently, Lamboley and Nahon \cite{LN} obtained planar boundary
regularity near points where the limiting local area fraction is
less than one, with a description of
the possible angular singularities. These results isolate an important
difficulty: existence does not by itself provide all the regularity
and geometric hypotheses needed to characterize an optimizer.

Tang and Zhang \cite{TZ} have announced a sharp planar
buckling--Dirichlet comparison in a public working preprint. Their
construction uses radial eigenfunctions to define a local transform
of a test function and its gradient, together with two divergence
corrections. We use this calibration and weak variational framework.
Our contribution is its higher-dimensional realization and the
analytic estimates that control all Hessian components through
dimension forty-seven. We also give a coefficient proof for the
planar case.

Our aim is to preserve the ball's sharp spectral ratio while avoiding
assumptions on a minimizing eigenfunction or on an optimal domain.
The result is stronger than the volume comparison, because it retains
the first Dirichlet eigenvalue of the original set.

\Needspace{6\baselineskip}
We first compare the two eigenvalues on the same open set.

\begin{theorem}\label{thm:main}
Let $2\le n\le47$. For every open set $\Omega\subset\mathbb R^n$
of finite positive volume,
\begin{equation}\label{eq:main}
 \Lambda_1\left(\Omega\right)\ge
 \frac{j_{n/2,1}^2}{j_{n/2-1,1}^2}\lambda_1^{\mathrm D}\left(\Omega\right).
\end{equation}
The constant is optimal and is attained on balls.
\end{theorem}

The Faber--Krahn inequality compares $\lambda_1^{\mathrm D}\left(\Omega\right)$
with that of an equal-volume ball; see \cite{PS,buckeq:BDV}.
It converts Theorem~\ref{thm:main} into the desired volume bound.
For arbitrary open sets, equality in the resulting buckling bound
requires an $H^2$-capacity condition. We write $C_{2,2}$ for the Bessel capacity associated with
$H^2\left(\mathbb R^n\right)$. A set of zero such capacity is negligible for the
value trace of $H^2$ functions. A precise definition and the relevant
trace criterion are given in Section~\ref{buckeq:section}.

\Needspace{6\baselineskip}
The spectral comparison gives the following sharp volume inequality and equality class.

\begin{corollary}\label{cor:ball}
Under the hypotheses of Theorem~\ref{thm:main}, let $\left|B_R\right|=\left|\Omega\right|$.
Then
\[
 \Lambda_1\left(\Omega\right)\ge\frac{j_{n/2,1}^2}{R^2}
 =\Lambda_1\left(B_R\right).
\]
Equality holds if and only if
\[
 \Omega=B_R\left(x_0\right)\setminus E,\qquad C_{2,2}\left(E\right)=0,
\]
where $E$ is relatively closed in $B_R\left(x_0\right)$.
For $n=2,3$, this means that $\Omega$ is a translated ball.
For $n\ge4$, sets of zero capacity may be removed.
\end{corollary}

The proof starts with the first buckling mode and the first Dirichlet
mode of a ball. Their quotient determines a scalar profile $f$ and a
local transform
$F=\left|\nabla U\right|f\left(U/\left|\nabla U\right|\right)$, extended continuously where the
gradient vanishes. Two divergence terms turn the difference between
the fourth-order and second-order energies into a quadratic expression
in the Hessian $D^2U$, the matrix of second derivatives. We call this
identity a \emph{calibration}: the ball profile is a null direction,
and a nonnegative remainder yields the sharp spectral comparison
after integration.

The new difficulty in dimension $n$ lies in the $n-1$ tangential
Hessian directions. Their trace couples to the radial entry, while
their trace-free part, obtained by subtracting the scalar matrix
with the same trace, has its own coefficient. Controlling only the
planar block would leave this term untreated. We prove the required
coefficient signs by a residue expansion in which the first remaining
pole dominates a possibly signed infinite tail. Positive products for
the radial quotient then give a rational energy estimate. This already
controls the matrix through dimension ten. In the higher dimensions,
a positive Laplace representation, together with two exact endpoint
moments, gives the stronger estimate needed for the determinant.
Together, these estimates make the Hessian remainder nonnegative for
$2\le n\le47$. They follow directly from the radial equations,
independently of the spectral comparison, and control entire parameter
intervals and infinite coefficient tails.

The global step is equally necessary. The transform is made globally
Lipschitz in the independent pair $\left(U,\nabla U\right)$, with the appropriate
parities across $U=0$. Locality of weak derivatives handles
$\left\{\nabla U=0\right\}$ even when it has positive measure. Hence no flux is
lost across nodal or critical sets. Dilation and Sobolev approximation
then prove \eqref{eq:main} without assuming that a first buckling
eigenfunction is positive, or even using existence of one.
Applying Faber--Krahn then gives the volume inequality.
The equality argument combines Faber--Krahn rigidity with the clamped Sobolev trace
criterion, which explains why zero capacity, rather than zero volume,
is the appropriate qualification.

The dimension restriction is intrinsic to this particular calibration.
With $q=j_{n/2-1,1}$ and $k=j_{n/2,1}$, the leading endpoint sign of
its coupled determinant is the sign of
$3\left(k^2-q^2\right)-4\left(n-1\right)$. We prove that it is positive for $3\le n\le47$
and negative for every $n\ge48$. Thus the method has an exact
dimensional endpoint. This does not disprove the buckling conjecture
in higher dimensions, nor does the present result settle the separate
Rayleigh vibration problem.

Section~\ref{sec:calibration} contains the calibration construction,
the sharp spectral comparison, and the equality argument.
Section~\ref{sec:coefficients} proves the coefficient signs in dimensions
three through forty-seven. Section~\ref{sec:multipliers} supplies the
radial bounds and determinant estimates, and
Section~\ref{sec:dimension-threshold} identifies the exact dimensional
endpoint of the calibration. Appendix~\ref{app:planar-coefficients}
contains the planar coefficient estimates.

\section{The calibration and the sharp comparison}\label{sec:calibration}
The proof begins with the ball profiles and a correction of the Hessian
energy. We establish positivity of the resulting matrix, pass to a
global functional inequality, and then identify the equality cases.
The coefficient and multiplier estimates used here are proved in
Sections~\ref{sec:coefficients} and~\ref{sec:multipliers}, with the planar
coefficient argument in Appendix~\ref{app:planar-coefficients}.

\subsection{Radial construction}\label{sec:radial}
All functions are real. We identify $H_0^j\left(\Omega\right)$ with the closure
of the zero extensions of $C_c^\infty\left(\Omega\right)$ in $H^j\left(\mathbb R^n\right)$;
the norms agree on these test functions. Integrals without a domain
are over the ambient Euclidean space, or over the interval explicitly
specified by their limits. The norm of a matrix is its Frobenius norm,
$\left|H\right|^2=\sum_{i,j}H_{ij}^2$. A symmetric matrix $M$ is positive
semidefinite if $v^{\mathsf T}Mv\ge0$ for every real vector $v$.
A prime or the subscript $r$ denotes a radial derivative; derivatives
with respect to other variables are indicated explicitly. The notation
$O\left(r^a\right)$ has its usual local meaning, with the dimension fixed. The relation $f\asymp g$ means that
the ratio of the two positive quantities stays between positive constants.
We use the Bessel series, recurrence relations, zero interlacing, and
canonical products in \cite[Sections 10.2, 10.6 and 10.21]{DLMF}.
For Sturm comparison, which compares the spacing of zeros of
second-order equations by ordering their potentials, we use
\cite[Chapter 5]{Teschl}.

For the radial construction put $\nu=n/2-1$, $q=j_{\nu,1}$,
$k=j_{\nu+1,1}$ and $\lambda=q^2/k^2$. Put $m=n-1$. Define
$\psi_n\left(r\right)=2^\nu\Gamma\left(\nu+1\right)r^{-\nu}J_\nu\left(r\right)$, continued by
$\psi_n\left(0\right)=1$, where
$\Gamma\left(z\right)=\int_0^\infty t^{z-1}e^{-t}\dd t$ for $z>0$
is the Euler gamma function \cite[Section 5.2(i)]{DLMF}.
The same formulas define the radial profiles for real $n\ge2$,
as used in some parameter estimates below. Put
\begin{equation*}
 c=-\psi_n\left(k\right)>0,\quad\sigma=1+c,\quad\alpha=q/k,\quad
 \lambda=\alpha^2,\quad\beta=-\frac{\alpha\psi_n'\left(q\right)}c>0,
 \qquad K=\beta^2.
\end{equation*}

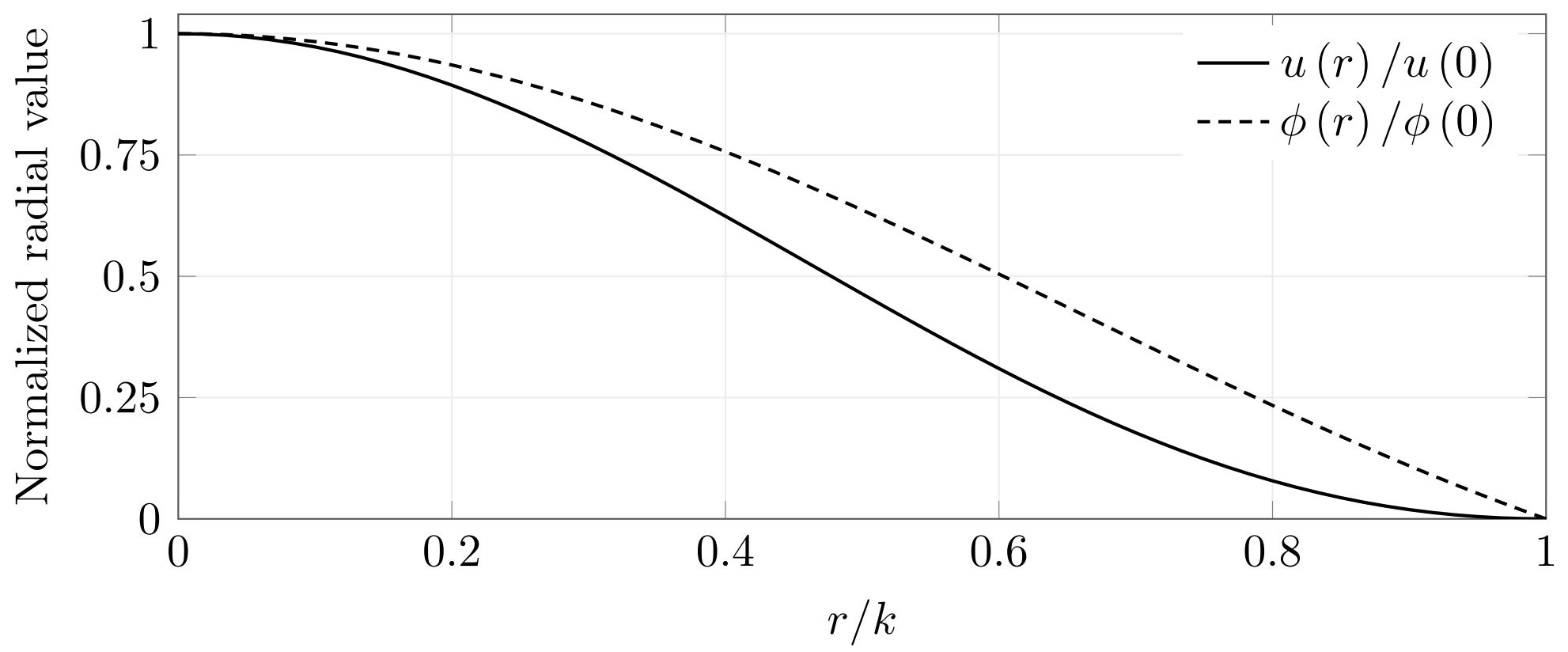
\begin{figure}[htbp]
\centering
\begin{tikzpicture}
\begin{axis}[
 width=.82\linewidth,height=5.8cm,
 xmin=0,xmax=1,ymin=0,ymax=1.04,
 xlabel={$r/k$},ylabel={Normalized radial value},
 xtick={0,.2,.4,.6,.8,1},ytick={0,.25,.5,.75,1},
 grid=major,grid style={gray!15},axis line style={black!65},
 tick label style={font=\small},label style={font=\small},
 legend style={font=\small,draw=none,fill=white,at={(.98,.98)},anchor=north east}]
\addplot[black,thick] coordinates {
(0,1) (0.01666667,0.9992323) (0.03333333,0.9969317) (0.05,0.993106) (0.06666667,0.9877681)
(0.08333333,0.9809358) (0.1,0.9726321) (0.1166667,0.9628849) (0.1333333,0.9517267) (0.15,0.939195)
(0.1666667,0.9253315) (0.1833333,0.9101824) (0.2,0.893798) (0.2166667,0.8762326) (0.2333333,0.8575444)
(0.25,0.8377949) (0.2666667,0.8170489) (0.2833333,0.7953743) (0.3,0.7728419) (0.3166667,0.7495246)
(0.3333333,0.7254978) (0.35,0.7008386) (0.3666667,0.675626) (0.3833333,0.6499399) (0.4,0.6238613)
(0.4166667,0.5974721) (0.4333333,0.5708542) (0.45,0.5440898) (0.4666667,0.5172607) (0.4833333,0.4904481)
(0.5,0.4637323) (0.5166667,0.4371926) (0.5333333,0.4109065) (0.55,0.3849499) (0.5666667,0.3593968)
(0.5833333,0.3343184) (0.6,0.3097839) (0.6166667,0.2858593) (0.6333333,0.2626077) (0.65,0.2400889)
(0.6666667,0.2183592) (0.6833333,0.1974713) (0.7,0.1774741) (0.7166667,0.1584124) (0.7333333,0.140327)
(0.75,0.1232544) (0.7666667,0.1072269) (0.7833333,0.09227225) (0.8,0.07841386) (0.8166667,0.06567054)
(0.8333333,0.05405663) (0.85,0.04358191) (0.8666667,0.03425167) (0.8833333,0.02606671) (0.9,0.0190234)
(0.9166667,0.01311378) (0.9333333,0.008325614) (0.95,0.004642504) (0.9666667,0.002044032) (0.9833333,0.0005058817)
(1,0)
};
\addlegendentry{$u\left(r\right)/u\left(0\right)$}
\addplot[black,thick,densely dashed] coordinates {
(0,1) (0.01666667,0.9995431) (0.03333333,0.9981733) (0.05,0.9958927) (0.06666667,0.9927052)
(0.08333333,0.9886159) (0.1,0.9836316) (0.1166667,0.9777605) (0.1333333,0.9710122) (0.15,0.9633978)
(0.1666667,0.9549297) (0.1833333,0.9456218) (0.2,0.9354893) (0.2166667,0.9245488) (0.2333333,0.9128181)
(0.25,0.9003163) (0.2666667,0.8870638) (0.2833333,0.873082) (0.3,0.8583937) (0.3166667,0.8430225)
(0.3333333,0.8269933) (0.35,0.810332) (0.3666667,0.7930651) (0.3833333,0.7752206) (0.4,0.7568267)
(0.4166667,0.737913) (0.4333333,0.7185093) (0.45,0.6986466) (0.4666667,0.678356) (0.4833333,0.6576696)
(0.5,0.6366198) (0.5166667,0.6152393) (0.5333333,0.5935615) (0.55,0.5716199) (0.5666667,0.5494483)
(0.5833333,0.5270807) (0.6,0.5045512) (0.6166667,0.4818939) (0.6333333,0.459143) (0.65,0.4363326)
(0.6666667,0.4134967) (0.6833333,0.390669) (0.7,0.367883) (0.7166667,0.345172) (0.7333333,0.3225687)
(0.75,0.3001054) (0.7666667,0.2778142) (0.7833333,0.2557263) (0.8,0.2338723) (0.8166667,0.2122824)
(0.8333333,0.1909859) (0.85,0.1700114) (0.8666667,0.1493865) (0.8833333,0.1291382) (0.9,0.1092924)
(0.9166667,0.08987418) (0.9333333,0.07090751) (0.95,0.05241541) (0.9666667,0.03441977) (0.9833333,0.01694141)
(1,0)
};
\addlegendentry{$\phi\left(r\right)/\phi\left(0\right)$}
\end{axis}
\end{tikzpicture}
\caption{The radial buckling and Dirichlet profiles in dimension three,
normalized to have value one at the center. At the boundary $r=k$,
the buckling profile has a double zero,
$u\left(k\right)=u_r\left(k\right)=0$, whereas
$\phi\left(k\right)=0$ and $\phi_r\left(k\right)<0$.
This difference explains why the calibration uses the quotient
$u/\left(-u_r\right)$ to compare the two profiles.
The curves illustrate the explicit definitions and endpoint expansions
in Section~\ref{sec:radial}; no sign assertion is inferred from the plot.}
\label{fig:radial-profiles}
\end{figure}

Let $X=\psi_n$, with $X\left(0\right)=1$, and put
\begin{equation*}
 p=-X_r=\left(r/n\right)\psi_{n+2}\left(r\right),\qquad
 u=X+c,\qquad\phi=X\left(\alpha r\right),\qquad0<r<k.
\end{equation*}
The recurrence and interlacing give $q<k<j_{\nu,2}$, $p>0$, $u>0$,
and $X\left(k\right)=-c$ with $c>0$. The equations are
\begin{equation*}
 X_{rr}+mX_r/r+X=0,\quad
 p_{rr}+mp_r/r+\left(1-m/r^2\right)p=0,\quad
 \phi_{rr}+m\phi_r/r+\lambda\phi=0.
\end{equation*}
Also $\Delta u=-\left(u-c\right)$ and $\Delta^2u+\Delta u=0$ on $B_k$.
The functions $u$ and $\phi$ are the ball buckling and Dirichlet
modes at the common scale $B_k$. Their distinct orders of vanishing
are shown in Figure~\ref{fig:radial-profiles}. The quotient $u/p$
turns these radial data into a profile defined on the whole positive
half-line.

\Needspace{6\baselineskip}
The next lemma makes the radial quotient a global coordinate for the calibration profiles.

\begin{lemma}\label{lem:radial-coordinate}
The function $t=u/p$ is a smooth strictly decreasing bijection
from $\left(0,k\right)$ onto $\left(0,\infty\right)$.
\end{lemma}
\begin{proof}
Set $S=p^2+up_r$. The canonical product gives
\[
 p\left(r\right)=\frac rn\prod_{l\ge1}\left(1-r^2/\kappa_l^2\right),\qquad
 \kappa_l=j_{\nu+1,l},\quad\kappa_1=k,
\]
\[
 \left(\log p\right)_{rr}=-r^{-2}-\sum_{l\ge1}
 \frac{2\left(\kappa_l^2+r^2\right)}{\left(\kappa_l^2-r^2\right)^2}<0.
\]
If $p_r\ge0$, then $S>0$. Otherwise $A=-p_r/p>0$ and strict
logarithmic concavity gives $p\left(v\right)<p\left(r\right)e^{-A\left(r\right)\left(v-r\right)}$ for $v>r$.
Hence $u=\int_r^kp\left(v\right)\dd v<p/A$ and $S=p^2-Aup>0$.
Thus $t_r=-S/p^2<0$. The endpoint limits
$t\sim n\sigma/r$ and $t\sim\left(k-r\right)/2$ prove the assertion.
\end{proof}
We use $s=t\left(r\right)$ as independent calibration variable. Subscripts $r$
denote radial derivatives, and subscripts $s$ derivatives after this change.

The corrections below are chosen so that the radial ball Hessian
lies in the kernel of the energy difference. This determines the
relations proved in Proposition~\ref{prop:kernel}.
Define, as functions of $r$ and equivalently of $s=t\left(r\right)$,
\begin{equation*}
 A=-\frac{p_r}{p},\quad C=-\frac1r,\quad
 f=\frac\phi p,\quad z=f_s,\quad g=f-sz,\quad w=\frac{\phi_r}{p}.
\end{equation*}
For the derivation it is useful to retain $m=n-1$ in the radial equations;
$1\le m\le46$ will be used in the main argument. The corresponding radial equations
give
\begin{align}
 t_r&=tA-1,& C_r&=C^2,&
 A_r&=A^2+mAC+1-mC^2,\label{eq:riccati}\\
 f_r&=w+Af=t_rz,& w_r&=\left(A+mC\right)w-\lambda f,&
 w&=-z-Ag.\label{eq:fw}
\end{align}
The last identity follows by writing $f=g+sz$ and $t_r=sA-1$.

Introduce the scalar $V_0$ and the two correction profiles
\begin{equation}\label{eq:G-b}
 V_0=w^2-mCfw+\lambda f^2,\qquad
 G=\frac{V_0-KA_r}{mC^2},\qquad
 \mathfrak{b}=-KA-mGC-fw.
\end{equation}
Here $V_0$ is a scalar, distinct from the vector-valued perspective map
used later. Since $C\ne0$, these formulas define smooth profiles in the
interior, without an additional differential-equation existence argument.

\Needspace{6\baselineskip}
The differential identities below provide the cancellation needed in the divergence correction.

\begin{lemma}\label{lem:profile-identities}
The profiles satisfy
\begin{align}
 G_r&=\mathfrak{b}+2\left\{A+\left(m-1\right)C\right\}G+2KC,\label{eq:G-ode}\\
 \mathfrak{b}_r+mCG_r&=-2f_rw.\label{eq:b-ode}
\end{align}
\end{lemma}
\begin{proof}
Differentiate using \eqref{eq:riccati}--\eqref{eq:fw}. The terms group as
\begin{align*}
 \left(V_0\right)_r&=\left(2A+mC\right)V_0-mC^2fw,\\
 A_{rr}&=\left(2A+mC\right)A_r+mAC^2-2mC^3.
\end{align*}
Differentiate $V_0=KA_r+mGC^2$ and cancel the common
$K\left(2A+mC\right)A_r$ term. Division by $mC^2$ yields
\[
 G_r=\left(2A+\left(m-2\right)C\right)G-KA+2KC-fw,
\]
which is \eqref{eq:G-ode} after substituting $\mathfrak{b}$. Also
$f_rw-fw_r=V_0$. Differentiating the definition of $\mathfrak{b}$ therefore gives
\[
 \mathfrak{b}_r+mCG_r=-KA_r-mGC^2-f_rw-fw_r
             =-V_0-f_rw-fw_r=-2f_rw.
\]
\end{proof}

Set
\begin{align}
 a&=K-g^2,& d&=2G-sG_s,& E&=\frac{K+\left(m-1\right)G}{m},\label{eq:coefficients}\\
 L_A&=2\mathfrak{b}-s\mathfrak{b}_s-2zg,& L_C&=\mathfrak{b}+G_s,&
 h&=-K-z^2+\lambda f^2+\mathfrak{b}_s.\notag
\end{align}
The symmetric matrix
\begin{equation*}
 M_m=\begin{pmatrix}
 a&d/2&L_A/2\\ d/2&E&L_C/2\\L_A/2&L_C/2&h
 \end{pmatrix}
\end{equation*}
will act on a radial Hessian entry, a tangential trace, and a gradient
magnitude.

\Needspace{6\baselineskip}
The next identities determine its kernel exactly.

\begin{proposition}\label{prop:kernel}
For $0<r<k$,
\begin{equation}\label{eq:kernel}
 L_A=-2aA-dmC,\qquad L_C=-dA-2EmC,\qquad
 h=aA^2+dmAC+Em^2C^2.
\end{equation}
Thus $M_m\left(A,mC,1\right)^{\mathsf T}=0$, and
\begin{equation*}
 \left(X,\tau,P\right)M_m\left(X,\tau,P\right)^{\mathsf T}
 =a\left(X-AP\right)^2+d\left(X-AP\right)\left(\tau-mCP\right)+E\left(\tau-mCP\right)^2.
\end{equation*}
\end{proposition}
\begin{proof}
For the second equation in \eqref{eq:kernel}, use $\left(1-sA\right)G_s=-G_r$:
\begin{align*}
 L_C+dA+2EmC
 &=\mathfrak{b}+\left(1-sA\right)G_s+2AG+2\left\{K+\left(m-1\right)G\right\}C\\
 &=\mathfrak{b}-G_r+2AG+2\left\{K+\left(m-1\right)G\right\}C=0.
\end{align*}
Dividing \eqref{eq:b-ode} by $t_r$ gives
$\mathfrak{b}_s+mCG_s=-2zw$. Consequently
\begin{align*}
 L_A+2aA+dmC
 &=2\left(\mathfrak{b}+KA+mGC\right)-s\left(\mathfrak{b}_s+mCG_s\right)-2g\left(z+Ag\right)\\
 &=-2fw+2szw+2gw=0,
\end{align*}
using $f=g+sz$ and $w=-z-Ag$.

For the last equation, the identities
$\mathfrak{b}_s=-2zw-mCG_s$ and $-z^2-2zw=w^2-A^2g^2$ reduce the difference to
\begin{align*}
 h-aA^2-dmAC-Em^2C^2
 ={}&w^2+\lambda f^2-K\left(1+A^2+mC^2\right)\\
 &-2mACG-m\left(m-1\right)GC^2+mCG_r.
\end{align*}
Substitute \eqref{eq:G-ode}, followed by $\mathfrak{b}=-KA-mGC-fw$.
The result is
\[
 V_0-K\left(A^2+mAC+1-mC^2\right)-mGC^2
 =V_0-KA_r-mGC^2=0.
\]
This proves all three equations and the completed-square expression.
\end{proof}

We record both endpoint expansions, since they will justify the global
Lipschitz extension. These are convergent power series after removal of the
displayed leading factors. Their differentiated versions are therefore valid.
At the radial center,
\begin{align*}
 p&=\frac rn-\frac{r^3}{2n\left(n+2\right)}+O\left(r^5\right),&u&=\sigma-\frac{r^2}{2n}+O\left(r^4\right),\\
 t&=\frac{n\sigma}r+\left(\frac{n\sigma}{2\left(n+2\right)}-\frac12\right)r+O\left(r^3\right),&
 f&=\frac nr+\left(\frac n{2\left(n+2\right)}-\frac\lambda2\right)r+O\left(r^3\right),\\
 w&=-\lambda+O\left(r^2\right),&A&=-r^{-1}+O\left(r\right),\\
 G&=G_0+O\left(r^2\right),&G_0&=\frac{n\lambda-K}{m},\qquad \mathfrak{b}=\frac{2n\lambda}r+O\left(r\right).
\end{align*}
In particular $1/t$ is analytic with a nonzero derivative at zero.
Its analytic inverse gives, as $s\to\infty$,
\begin{align}
 f&=s/\sigma+O\left(s^{-1}\right),& f_s&=\sigma^{-1}+O\left(s^{-2}\right),&g&=O\left(s^{-1}\right),\notag\\
 G&=G_0+O\left(s^{-2}\right),&G_s&=O\left(s^{-3}\right),\label{eq:large-s}\\
 \mathfrak{b}&=\left(2\lambda/\sigma\right)s+O\left(s^{-1}\right),&\mathfrak{b}_s&=2\lambda/\sigma+O\left(s^{-2}\right),&
 \mathfrak{b}-s\mathfrak{b}_s&=O\left(s^{-1}\right).\notag
\end{align}

For the outer endpoint let $\delta=k-r$ and temporarily keep $m$ in
the formulas (with $m=n-1$). The radial equations give
\begin{align}
 p&=c\delta\left(1+a_0\delta+b_0\delta^2+O\left(\delta^3\right)\right),\notag\\
 \phi&=\beta c\delta\left(1+a_0\delta+d_0\delta^2+O\left(\delta^3\right)\right),\label{eq:outer-series}\\
 a_0&=\frac m{2k},\qquad
 b_0=\frac{m\left(m+2\right)/k^2-1}{6},\qquad
 d_0=\frac{m\left(m+1\right)/k^2-\lambda}{6}.\notag
\end{align}
Also $u=c\delta^2/2+ca_0\delta^3/3+O\left(\delta^4\right)$ and
$t=\delta/2+O\left(\delta^2\right)$. With $h_\mathfrak{b}=1-\lambda-m/k^2$,
\begin{align*}
 f&=\beta+\frac{\beta h_\mathfrak{b}}6\delta^2+O\left(\delta^3\right),&
 z&=\frac{2\beta h_\mathfrak{b}}3\delta+O\left(\delta^2\right),\\
 g&=\beta-\frac{\beta h_\mathfrak{b}}6\delta^2+O\left(\delta^3\right).
\end{align*}
To check the removable singularity of $G$, use
\[
 A=\delta^{-1}+a_0+\left(2b_0-a_0^2\right)\delta+O\left(\delta^2\right),\quad
 w=-\beta\left\{\delta^{-1}+a_0+\left(3d_0-a_0^2-b_0\right)\delta+O\left(\delta^2\right)\right\}.
\]
In $V_0-KA_r$ the two pole coefficients vanish; its constant is
$K\left\{6d_0-m\left(m+1\right)/k^2+\lambda\right\}=0$. Thus $G=O\left(\delta\right)$
analytically. The pole and constant of $-KA-fw$ cancel too, so $\mathfrak{b}=O\left(\delta\right)$.
Substitution into the correction ODE then gives
\begin{equation*}
 G=\frac{2K}{3k}\delta+O\left(\delta^2\right),\quad
 d=\frac{2K}{3k}\delta+O\left(\delta^2\right),\quad
 \mathfrak{b}=\frac{2K\left(1-\lambda\right)}3\delta+O\left(\delta^2\right).
\end{equation*}
The equality of the linear factors in \eqref{eq:outer-series} ensures
$f_s\left(0\right)=0$. Extend $f,G$ evenly and $\mathfrak{b}$ oddly to the whole $s$ axis.
These extensions are locally Lipschitz; the even extension of $G$ may
have a corner at zero. With derivatives understood almost everywhere,
the endpoint expansions and interior smoothness give the global bounds
needed in the weak argument:
\begin{equation}\label{eq:bounded-profiles}
 f_s,\quad f-sf_s,\quad \mathfrak{b}_s,\quad \mathfrak{b}-s\mathfrak{b}_s,\quad G,\quad G_s,\quad sG_s
 \quad\text{are essentially bounded on }\mathbb R.
\end{equation}

\subsection{Positivity of the calibration matrix}
\label{sec:scalar-inputs}
We collect the analytic inputs used by the matrix proof. All quantities
come from the radial profiles already defined, and $S=p^2+up_r>0$.
Write $x=r^2$, $b=k^2$, $y=x/b$ and
\[
 F\left(x\right)=\frac{r\phi\left(r\right)}{p\left(r\right)},\qquad
 T\left(x\right)=\frac{ru\left(r\right)}{p\left(r\right)},\qquad
 \mathcal C=\frac{up}{S},\qquad L=\frac{1-y}{1+y}.
\]
The use of $T$ here is unrelated to the tangential Hessian block in
the later pointwise identity.

\Needspace{6\baselineskip}
The following scalar properties supply the inputs for the matrix argument.

\begin{proposition}\label{prop:scalar-inputs}
For every $2\le n\le47$, the following hold on $0<r<k$.
\begin{enumerate}
\item The correction has the expansion with positive weights
\[
 G=\sum_{j\ge1} B_j\left(1-y^j\right),\quad B_j>0,\qquad
 \sum_{j\ge1}jB_j=K/3.
\]
In particular $G_r<0$ and $0<G<K$.
\item $F$ has positive nonconstant Taylor coefficients. The function
$T$ is positive, decreasing and strictly concave as a function of $x$,
has negative nonconstant Taylor coefficients, and satisfies $T\left(b\right)=0$.
\item The radial multiplier satisfies $\mathcal C/r>L$, and
\[
 0<d=2G+\mathcal C G_r<\frac{2K}{3}L.
\]
\end{enumerate}
\end{proposition}
\begin{proof}
In dimension two, use \eqref{pc:G-coeff} and the endpoint slope in
Appendix~\ref{app:planar-coefficients}. In dimensions $n\ge3$, combine Lemma~\ref{an:coefficients} with
Theorems~\ref{generalG:theorem} and~\ref{Glow:theorem} in
Section~\ref{sec:coefficients}, use $G\left(k\right)=0$,
and put $B_j=a_j k^{2j}$. The endpoint slope gives the sum.
The product formula in Section~\ref{sec:radial-bounds} and
Proposition~\ref{prop:radial-series} prove the second assertion. Lemma~\ref{generalClower:main},
Corollary~\ref{generalClower:d} and Corollary~\ref{Cupper47:d}
there prove the third assertion. These scalar proofs do not use matrix positivity.
\end{proof}

We now use the positive Taylor coefficients of $F=r\phi/p$, the
positive weights in the expansion of $G$, and the negative
nonconstant Taylor coefficients of $T=rt$. Put $m=n-1$, $b=k^2$, $y=r^2/b$,
$h=1-y$, and
\[
 L=\frac{1-y}{1+y},\qquad \mathcal C=up/S,\qquad
 a=K-g^2,\quad E=\frac{K+\left(m-1\right)G}{m},\quad
 d=2G+\mathcal CG_r.
\]
For $2\le n\le47$, set
$\mathcal A\left(y\right)=F\left(by\right)/\left(\sqrt K\,k\right)$.
We use this normalized profile throughout the subsection.

\Needspace{6\baselineskip}
We first estimate the radial energy by a rational function of the normalized radius.

\begin{lemma}\label{mgf:rational-energy}
For every $0<y<1$,
\[
 \frac{g}{\sqrt K}<\frac{2\sqrt y}{1+y},\qquad
 \frac aK>L^2.
\]
\end{lemma}
\begin{proof}
Put $\mathcal T\left(y\right)=T\left(by\right)$ and
$D=-\mathcal T'/\mathcal T$. Endpoint regularity gives
$\mathcal A\left(1\right)=1$ and $\mathcal A'\left(1\right)=1/2$.
The positivity of all derivatives of $F$ implies convexity of $\mathcal A$.
Consequently
\[
 \mathcal A'\left(y\right)\le\frac{1-\mathcal A\left(y\right)}{1-y},\qquad
 \mathcal A\left(y\right)\ge\frac{1+y}{2}.
\]
The strictly negative nonconstant coefficients of $T$ imply that
$\mathcal T$ is decreasing and strictly concave in $y$. Since
$\mathcal T\left(1\right)=T\left(b\right)=0$,
\[
 0<D<\frac1{1-y}.
\]
The exact formula $g=2\sqrt x\left(TF_x-FT_x\right)/\left(T-2xT_x\right)$ gives
\[
 \frac{g}{\sqrt K\sqrt y}
 =\frac{2\left(\mathcal A'+D\mathcal A\right)}{1+2yD}
 \le\frac{1+\left(1+y\right)D}{1+2yD}<\frac2{1+y}.
\]
The same formula is strictly positive because
$\mathcal A,\mathcal A',D>0$, so it may be squared.
For the first inequality, substitute the secant bound for $\mathcal A'$
and then the lower tangent bound for $\mathcal A$; its coefficient
$D-1/\left(1-y\right)$ is negative. The last strict inequality is equivalent to
$\left(1-y\right)D<1$. Squaring proves the assertion for $a$.
\end{proof}

\Needspace{6\baselineskip}
This rational estimate already controls the coupled determinant through dimension ten.

\begin{theorem}
\label{mgf:matrix-low}
For every $2\le n\le10$, throughout $0<r<k$,
\[
 \left(K-g^2\right)E-d^2/4>0.
\]
\end{theorem}
\begin{proof}
Proposition~\ref{prop:scalar-inputs} gives
$0<d<\left(2K/3\right)L$. Also $E\ge K/m$, with equality only when $m=1$. Lemma~\ref{mgf:rational-energy}
therefore yields
\[
 \frac{d^2}{4aE}<\frac m9\le1.
\]
The strict bounds for $a$ and $d$ ensure strictness also at $n=2$ and $n=10$.
\end{proof}

\Needspace{6\baselineskip}
To reach higher dimensions, we estimate the correction and its derivative together.

\begin{lemma}
\label{mgf:joint-G}
Suppose $G=\sum_{j\ge1}B_j\left(1-y^j\right)$, $B_j>0$, and
$G_r\left(k\right)=-2K/\left(3k\right)$. If $d\ge0$, then
\[
 \frac{d^2}{4aE}
 \le \frac{mL^2}{3\left(a/K\right)\left\{3+\left(m-1\right)L\right\}}.
\]
\end{lemma}
\begin{proof}
The endpoint derivative gives $\sum jB_j=K/3$. Thus
$w_j=3jB_j/K$ are probability weights. Set
\[
 A_j=\frac{1-y^j}{j},\qquad c=\frac{\mathcal C}r>L,
 \qquad B_j^*=A_j-cy^j.
\]
Then
\[
 \frac GK=\frac13\sum w_j A_j,\quad
 \frac dK=\frac23\sum w_j B_j^*,\quad
 \frac EK=\sum w_j\left\{\frac1m+\frac{m-1}{3m}A_j\right\}.
\]
The geometric-sum inequality
$\left(1+y\right)\left(1-y^j\right)\le j\left(1-y\right)\left(1+y^j\right)$ gives
\[
 \left(B_j^*\right)_+\le\min\left\{A_j,L\right\}.
\]
Here $\left(t\right)_+=\max\left\{t,0\right\}$ denotes the positive part. Writing $e\left(A\right)=1/m+\left(m-1\right)A/\left(3m\right)$, the function $A^2/e\left(A\right)$
is increasing, whereas $L^2/e\left(A\right)$ is decreasing. Hence
$\left(B_j^*\right)_+^2/e\left(A_j\right)\le L^2/e\left(L\right)$. Cauchy--Schwarz yields
\[
 \left(\sum w_jB_j^*\right)^2
 \le\left(\sum w_j\left(B_j^*\right)_+\right)^2
 \le\left(\sum w_je\left(A_j\right)\right)\frac{L^2}{e\left(L\right)}.
\]
The first inequality uses $\sum w_jB_j^*\ge0$. Substitution proves
the result. Convergence follows from $\sum jB_j=K/3$.
\end{proof}

\Needspace{6\baselineskip}
Positive exponential integrals give a stronger curvature estimate, with constants fixed by two exact endpoint moments.

\begin{lemma}
\label{mgf:laplace}
The function $\mathcal A''$ is strictly positive and logarithmically
convex (that is, its logarithm is convex) on $\left[0,1\right]$. Its endpoint moments are
\[
 M_2=\mathcal A''\left(1\right)=\frac{k^2-q^2-m-3}{12},\qquad
 M_3=\mathcal A'''\left(1\right)=\frac{m+3}{8}.
\]
In particular,
\[
 \mathcal A''\left(y\right)\ge
 M_2\exp\left(-\frac{M_3}{M_2}\left(1-y\right)\right).
\]
\end{lemma}
\begin{proof}
The zero-ratio product for $F$ consists, after cancellation of the
first zero, of factors with $a>v>1$ of the form
\[
 \frac{1-y/a}{1-y/v}
 =\frac va+\left(1-\frac va\right)\frac v{v-y}.
\]
Indeed,
\[
 \frac{v}{v-y}=\int_0^\infty e^{yt}v e^{-vt}\dd t,
 \qquad y<v.
\]
Thus each factor is the exponential integral of a positive measure:
a mixture of a unit point mass at zero and an exponential probability
distribution of rate $v$. A finite product has the same form, with
the product represented by convolution of the measures. If $P$ is
such a product, Cauchy--Schwarz applied to its second derivative gives
\[
 \left(P^{\left(3\right)}\right)^2\le P''P^{\left(4\right)}.
\]
This is precisely the differential criterion for logarithmic convexity
of $P''$. The canonical products and their derivatives converge
locally uniformly before the second pole. Passing to the limit
therefore gives the same inequality for $\mathcal A''$; its strict
positivity follows from the positive Taylor coefficients. No limiting
measure is needed for this conclusion.

For an explicit check of the moments, set
$\mathfrak h\left(x\right)=p\left(\sqrt x\right)/\sqrt x$, and normalize
$H\left(y\right)=\mathfrak h\left(by\right)$ and $\Phi\left(y\right)=\phi\left(k\sqrt y\right)$ by separate
constant factors so that $H\left(1\right)=\Phi\left(1\right)=0$ and
$H'\left(1\right)=\Phi'\left(1\right)=1$. The harmless scaling makes $\Phi/H=\mathcal A$.
Their equations are
\[
 4yH''+2\left(n+2\right)H'+bH=0,\qquad
 4y\Phi''+2n\Phi'+q^2\Phi=0.
\]
Writing their Taylor coefficients at $y=1$ as $H_j,\Phi_j$ gives
\[
 H_{j+2}=-\frac{\left(j+1\right)\left(4j+2n+4\right)H_{j+1}+bH_j}
                    {4\left(j+1\right)\left(j+2\right)},
\]
and the same recurrence with $2n+4,b$ replaced by $2n,q^2$.
Division through order three gives
\[
 \mathcal A\left(1+s\right)=1+\frac s2
 +\frac{k^2-q^2-m-3}{24}s^2+\frac{m+3}{48}s^3+O\left(s^4\right).
\]
This proves the exact moments. The tangent inequality for
$\log\mathcal A''$ at $1$ proves the exponential lower bound.
\end{proof}

\Needspace{6\baselineskip}
A lower bound for this curvature strengthens the preceding rational energy estimate.

\begin{lemma}\label{lem:curvature-energy}
For any $0<C\le\mathcal A''\left(y\right)$,
\[
 \frac{g}{\sqrt K\sqrt y}<\frac{2-C\left(1-y\right)^2}{1+y},
 \qquad
 \frac{a/K}{L^2}>1+4yC-yC^2\left(1-y\right)^2.
\]
The right side of the first inequality is positive.
\end{lemma}
\begin{proof}
The nonnegative Taylor coefficients make $\mathcal A''$ increasing, so
$\mathcal A''\ge C$ on $\left[y,1\right]$. Strong convexity gives
\[
 \mathcal A\left(y\right)\ge1-\frac h2+\frac{Ch^2}{2},\qquad
 \mathcal A'\left(y\right)\le\frac{1-\mathcal A\left(y\right)}h-\frac{Ch}{2}.
\]
Substitute in the exact formula for $g$ as in
Lemma~\ref{mgf:rational-energy}. The result is
\[
 \frac{g}{\sqrt K\sqrt y}
 \le\frac{1-2Ch+\left(1+y+Ch^2\right)D}{1+2yD}.
\]
This fractional-linear expression has positive derivative in $D$:
its numerator after differentiation is $h+Ch\left(h+4y\right)>0$.
The strict bound $D<1/h$ now proves the first assertion. Its
right side is positive since $g>0$. Squaring it and subtracting
from one gives the second assertion.
\end{proof}

In dimensions above ten, combine Lemma~\ref{lem:curvature-energy} with
Lemma~\ref{mgf:joint-G} and the scalar exponential comparisons proved
in Section~\ref{sec:bessel-gaps}. The resulting determinant theorem is
given in Theorem~\ref{uniformcoupling:sharp}.

\Needspace{6\baselineskip}
The rational estimate in low dimensions and the Laplace estimate
in the remaining range now give the determinant sign throughout
the required dimension range. We then check the mixed and
trace-free Hessian coefficients.

\begin{theorem}
\label{uniformcoupling:sharp}
For every integer $2\le n\le47$, on the full interval $0<r<k$,
\[
 d^2<4E\left(K-g^2\right).
\]
\end{theorem}
\begin{proof}
For $2\le n\le10$ use Theorem~\ref{mgf:matrix-low}.
For $11\le n\le47$, put $m=n-1$ and $h=1-r^2/k^2$.
Corollary~\ref{Cupper47:d} gives $d>0$.
The rational energy bound gives $a=K-g^2>0$.
Choose $M$ from the corresponding row of the following table:
\[
\begin{array}{c|c}
10\le m\le26&M=\left(m-6\right)/24\\
27\le m\le36&M=\left(3m-24\right)/96\\
37\le m\le41&M=\left(13m-108\right)/432\\
42\le m\le46&M=\left(m-9\right)/36
\end{array}
\]
The gap estimates in Section~\ref{sec:bessel-gaps} prove
$M_2>M>0$. Since $M\mapsto M\exp\left(-M_3h/M\right)$ is increasing,
Lemma~\ref{mgf:laplace} yields
\[
 \mathcal A''\left(1-h\right)>C:=M\exp\left(-\frac{m+3}{8M}h\right).
\]
The strong energy estimate therefore gives
\[
 \frac{a/K}{L^2}>1+4\left(1-h\right)C-\left(1-h\right)C^2h^2,
 \qquad L=\frac h{2-h}.
\]
The elementary exponential comparisons in
Section~\ref{sec:bessel-gaps}, namely
Lemmas~\ref{mgfblock:10-26}, \ref{mgfblockmid:statement}
and~\ref{scalar-exp-matrix:comparison}, prove for these four ranges
\[
 1+4\left(1-h\right)C-\left(1-h\right)C^2h^2
 >\frac{m}{3\left\{3+\left(m-1\right)L\right\}}.
\]
Substitute this strict inequality in Lemma~\ref{mgf:joint-G}.
It gives $d^2/\left(4aE\right)<1$, as required.
\end{proof}

\Needspace{6\baselineskip}
The determinant bound also yields the signs needed for the mixed and trace-free Hessian terms.

\begin{corollary}
\label{uniformcoupling:five}
For every $2\le n\le47$ and $0<r<k$,
\[
 a=K-g^2>0,\quad E>0,\quad aE-d^2/4>0,\quad
 2K-g^2-d>0,\quad K-G>0.
\]
The null direction is given by \eqref{eq:kernel}.
\end{corollary}
\begin{proof}
The determinant estimate has been proved. For the remaining signs,
use the positive weights $B_j$ in the expansion of $G$.
The endpoint derivative gives $\sum jB_j=K/3$, so
\[
 0<G<G_0=\sum B_j<\sum jB_j=K/3.
\]
The strict inequality uses $B_2>0$. Thus $E>0$, $K-G>0$,
and $d=2G+\mathcal CG_r<2G<2K/3$.
The energy bound $g^2<K$ now gives $2K-g^2-d>K/3>0$.
The kernel identities depend only on the defining radial equations.
\end{proof}

The estimates now control every Hessian component. We next turn
this pointwise information into an inequality for arbitrary test functions.

\subsection{The functional comparison}
\label{sec:density}
The profiles are now fixed in a dimension $2\le n\le47$.
We first compute the local energy identity, then justify it across
the exceptional sets and pass to the variational infimum.

For the moment assume $U>0$ and $P=\left|\nabla U\right|>0$ at a point, and
evaluate all profiles at $s=U/P$. In an orthonormal frame with first
vector $e=\nabla U/P$, write
\[
 H=D^2U=\begin{pmatrix}X&\zeta^{\mathsf T}\\\zeta&T\end{pmatrix},
 \qquad \tau=\tr T,\qquad T_0=T-\frac\tau m I_m.
\]
Here $I_m$ is the identity matrix, $\zeta\in\mathbb R^m$, and
$T_0$ is the trace-free part of the tangential block $T$.
For $m=1$, $T_0=0$, so the same calculation includes dimension two.
The symbol $A:B=\operatorname{tr}\left(A^{\mathsf T}B\right)$ denotes the
matrix inner product. Define $T_1\left(H\right)=\left(\tr H\right)I_n-H$, the first Newton tensor. For a smooth
Hessian its divergence vanishes, since
$\partial_i\left(\Delta U\,\delta_{ij}-U_{ij}\right)=0$ (summed over $i$).
Consider
\begin{equation*}
 F=Pf\left(s\right),\qquad J_1=P\mathfrak{b}\left(s\right)\nabla U,\qquad
 J_2=G\left(s\right)T_1\left(H\right)\nabla U.
\end{equation*}

The identities $\nabla P=He$ and
$\nabla s=e-\left(s/P\right)He$ imply
\begin{equation}\label{eq:gradF}
 \nabla F=zPe+gHe,\qquad
 \left|\nabla F\right|^2=\left(zP+gX\right)^2+g^2\left|\zeta\right|^2.
\end{equation}
Direct differentiation of the first vector field gives
\begin{equation*}
 \diver J_1=\mathfrak{b}_sP^2+\left(2\mathfrak{b}-s\mathfrak{b}_s\right)PX+\mathfrak{b}P\tau.
\end{equation*}
For the second, use
\[
 e^{\mathsf T}T_1\left(H\right)e=\tau,\qquad
 e^{\mathsf T}HT_1\left(H\right)e=X\tau-\left|\zeta\right|^2,
\]
and
$T_1\left(H\right):H=2X\tau+\tau^2-\left|T\right|^2-2\left|\zeta\right|^2$.
The result is
\begin{equation}\label{eq:divJ2}
 \diver J_2=dX\tau+G\left(\tau^2-\left|T\right|^2\right)-d\left|\zeta\right|^2+G_sP\tau.
\end{equation}
Using $\left|T\right|^2=\left|T_0\right|^2+\tau^2/m$, combine
\eqref{eq:gradF}--\eqref{eq:divJ2} to obtain
\begin{align}
 &K\left(\left|H\right|^2-P^2\right)-\left(\left|\nabla F\right|^2-\lambda F^2\right)
                      +\diver J_1+\diver J_2\notag\\
 &\quad=aX^2+dX\tau+E\tau^2+L_APX+L_CP\tau+hP^2\notag\\
 &\qquad\quad +\left(2K-g^2-d\right)\left|\zeta\right|^2+\left(K-G\right)\left|T_0\right|^2.\label{eq:density-expanded}
\end{align}
All coefficients are those of \eqref{eq:coefficients},
with $m=n-1$. By Proposition~\ref{prop:kernel}, the right side equals
\begin{align}
 \calD={}&a\left(X-AP\right)^2+d\left(X-AP\right)\left(\tau-mCP\right)+E\left(\tau-mCP\right)^2\notag\\
 &+\left(2K-g^2-d\right)\left|\zeta\right|^2+\left(K-G\right)\left|T_0\right|^2.\label{eq:density-complete}
\end{align}
Corollary~\ref{uniformcoupling:five}
shows $\calD\ge0$. The trace-free term in \eqref{eq:density-complete}
is part of the $n$-dimensional identity and has been retained
throughout the scalar reduction.

It remains to justify the fields globally, including $P=0$. Simply
integrating the pointwise identity off that set would be insufficient:
a smooth function can have a critical set of positive measure.

We begin by extending the perspective maps across $v=0$.
For $v\in\R^n$, $P=\left|v\right|>0$, and $y\in\R$, define
\begin{equation*}
 \Phi\left(y,v\right)=P f\left(y/P\right),\qquad
 B\left(y,v\right)=P \mathfrak{b}\left(y/P\right),\qquad
 \mathcal V\left(y,v\right)=G\left(y/P\right)v.
\end{equation*}
Using the even/odd extensions from Section~\ref{sec:radial}, set
\begin{equation}\label{eq:perspective-limits}
 \Phi\left(y,0\right)=\left|y\right|/\sigma,\qquad
 B\left(y,0\right)=\left(2\lambda/\sigma\right)y,\qquad
 \mathcal V\left(y,0\right)=0.
\end{equation}
These are continuous extensions by \eqref{eq:large-s}. Continuity at
$\left(0,0\right)$ follows from the linear growth bounds for $f,\mathfrak{b}$ and boundedness
of $G$.

Off $v=0$ and the possible corner at $y=0$, differentiation gives
\begin{align*}
 \Phi_y&=f_s,& \nabla_v\Phi&=\left(f-sf_s\right)\frac vP,\\
 B_y&=\mathfrak{b}_s,& \nabla_vB&=\left(\mathfrak{b}-s\mathfrak{b}_s\right)\frac vP,\\
 \mathcal V_y&=G_s\frac vP,&
 D_v\mathcal V&=GI-sG_s\frac{vv^{\mathsf T}}{P^2}.
\end{align*}
Equation \eqref{eq:bounded-profiles} bounds these derivatives uniformly.
To conclude global Lipschitz continuity, integrate along line segments:
the equations $v=0$ and $y=0$ each cut such a segment in at most one
point unless they hold on the entire segment. In the latter case the
explicit restrictions \eqref{eq:perspective-limits}, or the values at
$y=0$, are Lipschitz. Subdivision at the remaining exceptional points
and continuity give the same uniform bound across the whole segment.

Now fix $U\in C_c^\infty\left(\Omega\right)$ and define globally
\begin{equation}\label{eq:global-fields}
 F=\Phi\left(U,\nabla U\right),\quad
 J_1=B\left(U,\nabla U\right)\nabla U,\quad
 J_2=T_1\left(D^2U\right)\mathcal V\left(U,\nabla U\right).
\end{equation}
The pair $\left(U,\nabla U\right)$ is globally Lipschitz. The smooth factor
$T_1\left(D^2U\right)$ and its derivatives are bounded on the support. Hence
$F,J_1,J_2$ are Lipschitz and compactly supported in $\Omega$.
In particular $F\in H_0^1\left(\Omega\right)$: mollification with support smaller
than the distance of $\supp U$ to $\Omega^c$ gives the required
approximation. Also $J_1,J_2\in W^{1,\infty}$, the space of bounded functions
with bounded weak first derivatives, with compact support,
so their weak divergences have zero integral.

The same fields also accommodate sign changes of $U$.
Replacing $U$ by $-U$ changes both $\nabla U$ and $D^2U$ to their
negatives. Evenness of $f$ gives the same $F$. Oddness of $\mathfrak{b}$ gives
$B\left(-U,-\nabla U\right)=-B\left(U,\nabla U\right)$, so $J_1$ is unchanged.
Evenness of $G$ changes $\mathcal V$ to its negative, and
$T_1\left(-H\right)=-T_1\left(H\right)$, so $J_2$ is unchanged as well. Therefore the left
side of \eqref{eq:density-expanded} is invariant, and its nonnegativity
on $\left\{U<0,P>0\right\}$ follows by applying the positive calculation to $-U$.
The set $\left\{U=0,\nabla U\ne0\right\}$ is covered by countably many smooth
local hypersurfaces, by the implicit function theorem, and has measure
zero. The fields themselves are globally Lipschitz across it, so no
surface distribution or unaccounted nodal flux appears.

It remains to verify the inequality on the critical set.
Let $Z=\left\{\nabla U=0\right\}$. We use Sobolev locality: two Sobolev
functions that agree on a measurable set have equal weak gradients
almost everywhere there. The chain rule, locality, and the
approximation facts used here are standard; see \cite[Chapter 2]{Ziemer}.

Apply this fact to each component of $\nabla U$. It gives $D^2U=0$
almost everywhere on $Z$. Moreover $J_1=J_2=0$ pointwise on $Z$,
by \eqref{eq:global-fields} and $\mathcal V\left(U,0\right)=0$. Since these
fields are Lipschitz, every weak first derivative of their components
vanishes almost everywhere on that common zero set. Thus
$\diver J_1=\diver J_2=0$ almost everywhere on $Z$.

The Lipschitz function $F-\left|U\right|/\sigma$ also vanishes on $Z$, so its
weak gradient is zero there. Since $\nabla\left|U\right|=0$ almost everywhere on
$Z$, it follows that $\nabla F=0$ there. The left side of
\eqref{eq:density-expanded} consequently reduces to
\[
 \lambda F^2=\lambda U^2/\sigma^2\ge0
 \quad\hbox{almost everywhere on }Z.
\]
This proof does not differentiate a perspective map at its corner and
does not assume the critical set has measure zero.

The preceding pointwise and regularity arguments establish, almost everywhere in $\Omega$,
\[
 K\left(\left|D^2U\right|^2-\left|\nabla U\right|^2\right)-\left(\left|\nabla F\right|^2-\lambda F^2\right)
                       +\diver J_1+\diver J_2\ge0.
\]
Integrate and use the compact support of both vector fields. For
compactly supported smooth $U$, integration by parts twice gives
\[
 \int U_{ii}U_{jj}\dd x=\int U_{ij}^2\dd x,
 \qquad \int\left(\Delta U\right)^2\dd x=\int\left|D^2U\right|^2\dd x.
\]
We obtain the functional inequality
\begin{equation}\label{eq:target-functional}
 \int\left|\nabla F\right|^2\dd x-\lambda\int F^2\dd x
 \le K\left\{\int\left(\Delta U\right)^2\dd x-\int\left|\nabla U\right|^2\dd x\right\}.
\end{equation}

The function $F$ is not zero when $U$ is not zero. At a point with
$U\ne0$, either $P>0$, in which case $F=Pf\left(\left|U\right|/P\right)>0$, or $P=0$,
in which case $F=\left|U\right|/\sigma>0$. By continuity it is positive on a
set of positive measure. It is therefore an admissible nonzero
Dirichlet trial function.

For $0\ne U\in C_c^\infty\left(\Omega\right)$ put
\[
 \mathcal Q=\frac{\int\left(\Delta U\right)^2\dd x}{\int\left|\nabla U\right|^2\dd x}>0,
 \qquad \rho=\sqrt{\mathcal Q},\qquad
 W\left(x\right)=U\left(x/\rho\right),\quad x\in\rho\Omega.
\]
In dimension $n$ the numerator of the buckling quotient scales by
$\rho^{n-4}$ and its denominator by $\rho^{n-2}$, so $W$ has quotient one.
Apply \eqref{eq:target-functional} to $W$ on $\rho\Omega$. Its right
side is zero, and the Dirichlet variational principle yields
\[
 \lambda_1^{\mathrm D}\left(\rho\Omega\right)
 \le\frac{\int\left|\nabla F_W\right|^2\dd x}{\int F_W^2\dd x}\le\lambda.
\]
Since $\lambda_1^{\mathrm D}\left(\rho\Omega\right)=
\lambda_1^{\mathrm D}\left(\Omega\right)/\mathcal Q$, we obtain
$\mathcal Q\ge\lambda_1^{\mathrm D}\left(\Omega\right)/\lambda$.

For any nonzero $U\in H_0^2\left(\Omega\right)$ choose smooth compactly supported
$U_j\to U$ in $H^2$. Both quadratic integrals in its quotient converge,
and the limiting gradient denominator is positive. The quotients
therefore converge as well. Taking the infimum proves
Theorem~\ref{thm:main} for $2\le n\le47$.

The $n$-dimensional Faber--Krahn inequality gives
$\lambda_1^{\mathrm D}\left(\Omega\right)\ge q^2/R^2$ for $\left|B_R\right|=\left|\Omega\right|$.
This is the standard form-domain version of Faber--Krahn,
valid for arbitrary finite-volume open sets; see \cite{buckeq:BDV}.
The ball formula $\Lambda_1\left(B_R\right)=k^2/R^2$ quoted in the introduction
therefore proves the inequality in Corollary~\ref{cor:ball}.
It also shows that the constant in Theorem~\ref{thm:main} is optimal.
The remaining issue is the exact equality class.

\subsection{Equality and Sobolev capacity}
\label{buckeq:section}

For $s=1,2$, let $\mathcal G_s$ be the Bessel kernel, whose
Fourier transform is $\left(1+\left|\xi\right|^2\right)^{-s/2}$. The Bessel capacity is
\[
 C_{s,2}\left(E\right)=\inf\left\{\left\|h\right\|_{L^2\left(\mathbb R^n\right)}^2:
 h\ge0,\ \mathcal G_s*h\ge1\text{ on }E\right\}.
\]
The convolution is defined for nonnegative $h$, with the value
$+\infty$ allowed. We use only capacity-zero sets, so the normalization
of the Fourier transform is immaterial. A statement holds
\emph{quasi-everywhere}, abbreviated q.e., if its exceptional set has
zero specified capacity. A quasicontinuous representative is one that
becomes continuous after removal of an open set of arbitrarily small
capacity. These representatives exist for Sobolev functions and are
unique up to capacity-zero sets; see \cite{AH}. In particular,
$C_{2,2}\left(E\right)=0$ implies $\left|E\right|=0$ and $C_{1,2}\left(E\right)=0$.

The equality argument depends on the spectral comparison only through
the following hypothesis; the converse assertion holds in every dimension.
Put $q_n=j_{n/2-1,1}$ and $k_n=j_{n/2,1}$.

\Needspace{6\baselineskip}
The following statement separates the capacity argument from the construction of the calibration.

\begin{theorem}
\label{buckeq:capacity}
Fix a dimension $n\ge2$ in which the strong comparison
\begin{equation}\label{buckeq:ratio-hyp}
 \Lambda_1\left(D\right)\ge \frac{k_n^2}{q_n^2}\lambda_1^{\mathrm D}\left(D\right)
 \quad\hbox{for every open }D\subset\mathbb R^n,
 \quad 0<\left|D\right|<\infty
\end{equation}
has been proved.  Let $\Omega$ be such an open set and let
$\left|B_R\right|=\left|\Omega\right|$.  Then
\begin{equation}\label{buckeq:iff}
 \Lambda_1\left(\Omega\right)=\Lambda_1\left(B_R\right)
 \quad\Longleftrightarrow\quad
 \Omega=B_R\left(x_0\right)\setminus E
 \ \hbox{for some }x_0\in\mathbb R^n,
 \quad C_{2,2}\left(E\right)=0.
\end{equation}
Here $E$ is relatively closed in $B_R\left(x_0\right)$, as is forced by openness
of $\Omega$. Consequently this applies in every dimension $2\le n\le47$
by Theorem~\ref{thm:main}. The sufficiency direction in \eqref{buckeq:iff} is valid in every dimension,
independently of \eqref{buckeq:ratio-hyp}.
\end{theorem}

We use the clamped Sobolev trace criterion
\begin{equation}\label{buckeq:synthesis}
 v\in H^2_0\left(D\right)
 \quad\Longleftrightarrow\quad
 \begin{cases}
  \widetilde v=0\quad C_{2,2}\text{-q.e. on }D^c,\\
  \widetilde{\partial_i v}=0\quad C_{1,2}\text{-q.e. on }D^c,
       &1\le i\le n,
 \end{cases}
 \qquad v\in H^2\left(\mathbb R^n\right).
\end{equation}
The representatives in the two lines are the appropriate, possibly
different, quasicontinuous representatives.  This is the $m=p=2$
closed-set synthesis theorem of Hedberg--Wolff
\cite[Theorem 5]{buckeq:HW}; in particular it is valid without boundary
regularity and without an embedding into $C^1$.

\begin{proof}[Proof of Theorem~\ref{buckeq:capacity}]
Suppose first that $\Lambda_1\left(\Omega\right)=k_n^2/R^2$.
The strong comparison and Faber--Krahn give
\[
 \frac{k_n^2}{R^2}=\Lambda_1\left(\Omega\right)
 \ge\frac{k_n^2}{q_n^2}\lambda_1^{\mathrm D}\left(\Omega\right)
 \ge\frac{k_n^2}{q_n^2}\frac{q_n^2}{R^2}.
\]
Thus $\lambda_1^{\mathrm D}\left(\Omega\right)=\lambda_1^{\mathrm D}\left(B_R\right)$.
The Fraenkel asymmetry of $\Omega$ is
\[
 \mathcal A\left(\Omega\right)=\inf_{x_0\in\mathbb R^n}
 \frac{\left|\Omega\mathbin\triangle B_R\left(x_0\right)\right|}{\left|\Omega\right|},
\]
where $\triangle$ denotes symmetric difference and $\left|B_R\right|=\left|\Omega\right|$.
The quantitative Faber--Krahn theorem
\cite[Main Theorem, $q=2$]{buckeq:BDV} implies zero Fraenkel asymmetry.
Hence $\Omega$ agrees almost everywhere with a translate of $B_R$.
For completeness, the centers of balls approaching zero symmetric
difference cannot escape to infinity: intersection of such a ball with
the fixed finite-measure set $\Omega$ would then tend to zero.  A bounded
subsequence of centers converges; continuity of ball translations in
$L^1$ gives one ball $B=B_R\left(x_0\right)$ with $\left|\Omega\mathbin\triangle B\right|=0$.
Since $\Omega$ is open, it follows that $\Omega\subset B$ pointwise.
Indeed, a point of $\Omega\setminus B$, including a point on $\partial B$,
has an open neighborhood whose intersection with $B^c$ has positive
measure.  This would contradict the almost-everywhere equality.

We next justify attainment on $\Omega$, rather than assume it from an
unbounded-domain statement.  Take a minimizing sequence in $H^2_0\left(\Omega\right)$
with $\int\left|\nabla v_j\right|^2\dd x=1$.  Its zero extensions belong to $H^2_0\left(B\right)$.
The identity $\int\left|D^2v_j\right|^2\dd x=\int\left|\Delta v_j\right|^2\dd x$, followed by Poincar\'e
on the bounded ball, gives boundedness in $H^2\left(B\right)$.
Weak compactness, closedness of the zero-extended subspace $H^2_0\left(\Omega\right)$,
and the Rellich compact embedding $H^2\left(B\right)\hookrightarrow H^1\left(B\right)$ \cite{Ziemer} give a nonzero minimizer
$v\in H^2_0\left(\Omega\right)$ with $\int\left|\nabla v\right|^2\dd x=1$.
It is also a first buckling eigenfunction of $B$.

The first buckling eigenspace of a ball is one-dimensional and is
spanned by a smooth radial function $U_B$ that is strictly positive in
the interior; see \cite[Theorem 3.2 and Section 6, $\kappa=0$]{buckeq:DNT}.
Thus $v=cU_B$ in $H^2\left(\mathbb R^n\right)$ for a nonzero constant $c$.
Let $E=B\setminus\Omega$.  By the first line of
\eqref{buckeq:synthesis}, $\widetilde v=0$ $C_{2,2}$-q.e. on $E$.
But the continuous representative $cU_B$ is nonzero at every point
of $E$, and it agrees with $\widetilde v$ outside a $C_{2,2}$-null set.
Consequently $C_{2,2}\left(E\right)=0$. Conversely, suppose $\Omega=B\setminus E$ and $C_{2,2}\left(E\right)=0$.
For any $v\in H^2_0\left(B\right)$, its value trace vanishes $C_{2,2}$-q.e.
on $B^c$ and its derivative traces vanish $C_{1,2}$-q.e. there.
The same statements hold on $B^c\cup E=\Omega^c$ since
$C_{1,2}\left(E\right)=0$ as well.  Formula~\eqref{buckeq:synthesis} yields
$v\in H^2_0\left(\Omega\right)$.  The reverse inclusion is immediate, so that
\begin{equation*}
 H^2_0\left(\Omega\right)=H^2_0\left(B\right)
 \quad\hbox{as zero-extended subspaces of }H^2\left(\mathbb R^n\right).
\end{equation*}
Their buckling infima are equal, and $\left|\Omega\right|=\left|B\right|$ because
$C_{2,2}$-null sets have zero Lebesgue measure.
\end{proof}

\Needspace{6\baselineskip}
Equality in the volume inequality must be distinguished from equality in the spectral ratio.

\begin{remark}\label{rem:equality-distinction}
Theorem~\ref{buckeq:capacity} classifies equality in the volume bound.
The ratio comparison itself also attains equality on a finite disjoint
union of equal-radius balls: both first variational values are the
minimum of the corresponding component values. Its volume bound is
strict when there is more than one component. The theorem makes no
claim to classify equality in the ratio comparison alone.
\end{remark}

A singleton has positive $C_{s,2}$ capacity exactly when $2s>n$
\cite[Chapter 2]{AH}. Applying this fact with $s=2$ to
Theorem~\ref{buckeq:capacity} shows that equality forces a translated
ball in dimensions two and three. In dimensions $n\ge4$, any interior
point $x_1$ can be removed without changing $H_0^2\left(B\right)$ or
the first buckling eigenvalue. This sufficiency assertion does not
require the spectral comparison.

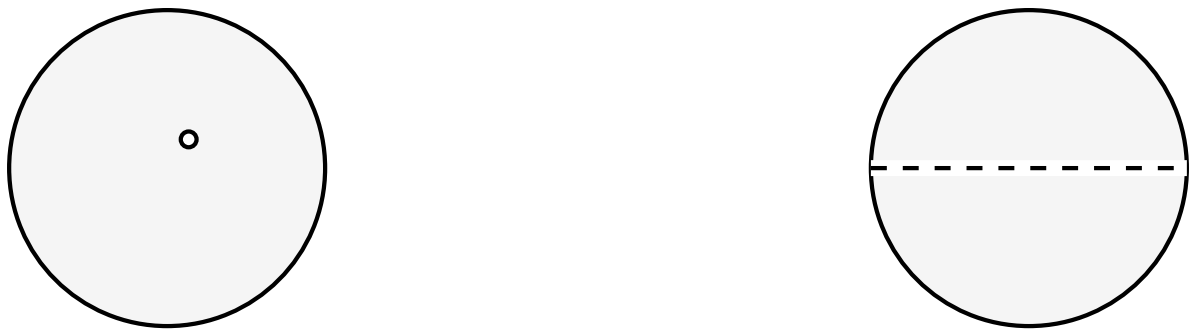
\begin{figure}[tbp]
\centering
\begin{tikzpicture}[scale=.95,every node/.style={font=\small}]
\draw[fill=gray!8,thick] (0,0) circle (1.1);
\draw[fill=white,thick] (.15,.2) circle (.055);
\node at (0,-1.5) {$B\setminus\left\{x_1\right\}$};
\node[text centered,text width=3.8cm] at (0,-2.15) {A point is removable\\when $n\ge4$.};
\draw[fill=gray!8,thick] (6,0) circle (1.1);
\draw[white,line width=3pt] (4.9,0)--(7.1,0);
\draw[dashed,thick] (4.9,0)--(7.1,0);
\node at (6,-1.5) {$B\setminus\left\{x_n=0\right\}$};
\node[text centered,text width=3.8cm] at (6,-2.15) {An internal clamped cut\\has positive capacity.};
\end{tikzpicture}
\caption{Schematic sections of equal-volume sets. Removing a set of
zero volume need not preserve the clamped function space. A point has
zero $H^2$ capacity for $n\ge4$, whereas a hyperplane disk has positive
capacity. The circles represent sections, not a planar puncture example.}
\label{fig:capacity}
\end{figure}

The qualification by capacity cannot be replaced by equality almost
everywhere. For $2\le n\le47$, the set
$D=B_R\setminus\left\{x_n=0\right\}$ has the same volume as $B_R$ but
has strictly larger first buckling eigenvalue. Indeed, the deleted
hyperplane disk has positive $C_{2,2}$ capacity, by the continuous
local hyperplane trace map for $H^2$ \cite{Ziemer}, and
Theorem~\ref{buckeq:capacity} applies. Figure~\ref{fig:capacity}
illustrates this distinction. The equality proof uses the
$C_{2,2}$ value trace as well as the $C_{1,2}$ derivative traces;
Dirichlet rigidity alone does not determine the clamped equality class.

\section{Coefficient signs in dimensions three through forty-seven}\label{sec:coefficients}
We now establish the coefficient estimates used in
Proposition~\ref{prop:scalar-inputs} for $3\le n\le47$.
The three-dimensional profiles admit an explicit trigonometric form.
For higher dimensions, a residue expansion reduces the proof to
comparison of the first retained pole with the signed tail; the same
structure then yields the estimates in dimensions four and five.
Appendix~\ref{app:planar-coefficients} supplies the separate planar argument.

We write $\left[r^j\right]H$ for the coefficient of $r^j$ in the Taylor
series of $H$.

We shall use $157/50<\pi<22/7$. The tangent addition formula gives Machin's identity
$\pi=16\arctan\left(1/5\right)-4\arctan\left(1/239\right)$, and its alternating series give
\begin{align*}
 \pi&>16\left(\frac15-\frac1{3\cdot5^3}\right)-\frac4{239}
       >\frac{157}{50},\\
 \pi&<16\left(\frac15-\frac1{3\cdot5^3}+\frac1{5\cdot5^5}\right)
       -4\left(\frac1{239}-\frac1{3\cdot239^3}\right)<\frac{22}{7}.
\end{align*}

\subsection{The three-dimensional correction coefficients}

For $n=3$ the common profiles become $X\left(r\right)=\sin r/r$ and
$p\left(r\right)=\left(\sin r-r\cos r\right)/r^2$, with $q=\pi$. For this subsection put
\[
 D\left(r\right)=\sin r-r\cos r,\qquad
 \mathcal L\left(r\right)=\left(K-1\right)r^2-K\sin^2r+\alpha^{-2}\sin^2\left(\alpha r\right).
\]
Direct substitution into \eqref{eq:G-b} gives
\begin{equation}\label{n3:explicit-G}
 G=K-\frac{r^2\mathcal L}{2D^2}.
\end{equation}
Write $\kappa_l$ for the positive zeros of $D$, with $\kappa_1=k$.
Since $\tan k=k$ and $\cos k<0$, one has
\[
 c=\left(1+k^2\right)^{-1/2},\qquad K=1+k^{-2}>1.
\]
The following coarse bounds suffice:
\begin{equation}\label{an:coarse}
 4.49<k<4.50,\quad \frac{12}{25}<\lambda<\frac{49}{100}<\frac12,
 \quad \frac{27}{125}<c<\frac{11}{50},\quad K<\frac{21}{20}.
\end{equation}
For completeness the root enclosure can be obtained directly from
the alternating sine and cosine series: if
$S_N\left(x\right)=\sum_{j=0}^N\left(-1\right)^jx^{2j+1}/\left(2j+1\right)!$ and
$C_N\left(x\right)=\sum_{j=0}^N\left(-1\right)^jx^{2j}/\left(2j\right)!$, the two rational inequalities
\[
 S_{11}\left(449/100\right)-\left(449/100\right)C_{12}\left(449/100\right)>0,\qquad
 S_{12}\left(9/2\right)-\left(9/2\right)C_{11}\left(9/2\right)<0
\]
give the assertion by the alternating remainder bound and $D'=r\sin r<0$
on $\left(\pi,3\pi/2\right)$. The remaining bounds in \eqref{an:coarse} follow
from $157/50<\pi<22/7$ and the explicit expressions for $c,K,\lambda$.
The bounds for $\pi$ were established at the start of this section.

\Needspace{6\baselineskip}
These explicit profiles yield the coefficient signs required for the three-dimensional calibration.

\begin{lemma}\label{an:coefficients}
There are $a_j>0$ such that, with convergence of all derivatives on $\left[0,k\right]$,
\[
 G\left(r\right)=G\left(0\right)-\sum_{j\ge1}a_jr^{2j}
      =\sum_{j\ge1}a_j\left(k^{2j}-r^{2j}\right).
\]
In particular $G>0$, $G_r<0$ and $G_{rr}<0$ on $\left(0,k\right)$.
\end{lemma}
\begin{proof}
Use the entire function $\mathcal L$ and \eqref{n3:explicit-G}. The singularities at $0,\pm k$ are removable by
Section~\ref{sec:radial}. All other nonzero zeros of $D$ are the nonzero zeros of
$J_{3/2}$, and are real and simple \cite[Section 10.21(i)]{DLMF}.

At $\kappa=\kappa_l$, $l\ge2$, one has
$D''\left(\kappa\right)/D'\left(\kappa\right)=2/\kappa$. Consequently $r^2/D^2$ has no simple
pole term, and the principal part of $G$ at $\kappa$ is
\[
 -\frac{\left(1+\kappa^2\right)\mathcal L\left(\kappa\right)}{2\kappa^2\left(r-\kappa\right)^2}
 -\frac{\left(1+\kappa^2\right)\mathcal L'\left(\kappa\right)}{2\kappa^2\left(r-\kappa\right)}.
\]
Evenness determines the principal part at $-\kappa$. For $j\ge2$,
paired residues therefore give
\begin{equation}\label{an:residues}
 a_j=\sum_{l\ge2}\frac{1+\kappa_l^2}{\kappa_l^{2j+4}}
 \left(\left(2j+1\right)\mathcal L\left(\kappa_l\right)-\kappa_l \mathcal L'\left(\kappa_l\right)\right).
\end{equation}
Here is a justification of the contour limit. Take the square with
vertical sides $\Re r=\pm N\pi$ and horizontal sides $\Im r=\pm N\pi$.
On a vertical side, with $y=\Im r$,
\[
 \left|D\right|^2=\left(N\pi\right)^2\cosh^2y+\left(\sinh y-y\cosh y\right)^2.
\]
On a horizontal side the inequality
$\left|D\right|\ge \left|r\cos r\right|-\left|\sin r\right|$ gives $\left|D\right|\ge c_1\left|r\right|e^{N\pi}$ for large $N$.
Since $\alpha<1$, these bounds give $G=O\left(N^2\right)$ uniformly on the contour.
The integral of $G\left(r\right)/r^{2j+1}$ is thus $O\left(N^{2-2j}\right)\to0$ for $j\ge2$.
The residue series converges absolutely, its terms being $O\left(\kappa_l^{-2j}\right)$.
No contour claim for $j=1$ is needed.

At a pole put $x=\alpha\kappa$ and
\[
 P=\frac{\kappa^2\left(\kappa^2-k^2\right)}{k^2\left(1+\kappa^2\right)}>0.
\]
Using $\sin^2\kappa=\kappa^2/\left(1+\kappa^2\right)$ gives
\begin{align*}
 \mathcal L\left(\kappa\right)&=P+\lambda^{-1}\sin^2x>0,\\
 5\mathcal L\left(\kappa\right)-\kappa \mathcal L'\left(\kappa\right)
 &=3P+\lambda^{-1}\left(5\sin^2x-x\sin\left(2x\right)\right).
\end{align*}
This last expression is positive for every $l\ge2$, as follows.
For $l=2$, the bounds $9\pi/4<\kappa_2<5\pi/2$ give
$3\pi/2<x<2\pi$, hence $\sin\left(2x\right)<0$.
For $l=3$, $13\pi/4<\kappa_3<7\pi/2$ gives
$P>4$, $9\pi/4<x<11\pi/4$, and $x<7.8$; therefore the expression exceeds
\[
 12-\frac{7.8-2.5}{12/25}=\frac{23}{24}>0.
\]
The bounds on $\kappa_2,\kappa_3$ use $\tan r=1<r$ at the respective
lower endpoints. For $l\ge4$, $\kappa_l>14$: indeed
$\epsilon=9\pi/2-14>0.13$ gives $\tan14=\cot\epsilon<1/\epsilon<14$.
Also
\[
 P>\frac{4\kappa^2}{81}-\frac{21}{20},\qquad \alpha^{-1}<\frac32,
\]
so the expression is greater than
$4\kappa^2/27-63/20-3\kappa/2>0$ for $\kappa\ge14$.
Thus every term in \eqref{an:residues} is positive for $j\ge2$.
Direct expansion at zero gives
\[
 a_1=\frac{K-3\lambda+2\lambda^2}{10}>0,
\]
since $3\lambda-2\lambda^2<1<K$. The nearest possible remaining pole
lies at $\kappa_2>k$, which justifies convergence and differentiation
on $\left[0,k\right]$. Finally $G\left(k\right)=0$ gives the second series representation.
\end{proof}

\subsection{The first pole and the tail in higher dimensions}

This subsection proves the calibration coefficient signs in every integer
dimension $6\le n\le47$. Uniform pole estimates, quantitative
quarter-phase comparisons, and rational polynomial inequalities in a
continuous dimension parameter give the proof.
Put $m=n-1$, $\nu=\left(n-2\right)/2$, $X=\psi_n$, $p=-X_r=r\psi_{n+2}/n$,
$q_l=j_{\nu,l}$, $\kappa_l=j_{\nu+1,l}$, $q=q_1$, $k=\kappa_1$,
$\alpha=q/k$, $\lambda=\alpha^2$, $c=-X\left(k\right)$, and
$K=\lambda p\left(q\right)^2/c^2$. Write
\[
 B=X^2+p^2-mXp/r,\quad Q\left(r\right)=KB\left(r\right)-\lambda B\left(\alpha r\right),\quad
 W\left(r\right)=r^mB\left(r\right),\quad G=K-\frac{r^2Q}{mp^2}.
\]

\Needspace{6\baselineskip}
The first step is to determine the sign of the residue factor at every pole.

\begin{lemma}
\label{generalG:Qpoles}
For every real $n\ge3$ and every $l\ge2$, $Q\left(\kappa_l\right)>0$.
Moreover $q_l/\kappa_l$ strictly increases with $l$.
\end{lemma}
\begin{proof}
The identities
\[
 \left(X^2+p^2\right)_r=-2mp^2/r,\qquad W_r=mr^{m-2}Xp,
\]
follow from $X_r=-p$ and $p_r=X-mp/r$.
One also has
\[
 B\left(r\right)=r^{-m-1}\int_0^r s^m\left\{X\left(s\right)^2+p\left(s\right)^2\right\}\dd s>0.
\]
Indeed the derivative of $r^{m+1}B$ is $r^m\left(X^2+p^2\right)$,
and its value at zero is zero.

Choose a continuous polar phase $X=\rho\cos\theta$,
$p=\rho\sin\theta$. On the quarter-cycle from $q_l$ to
$\kappa_l$ set $\varphi=\theta-\left(l-1\right)\pi-\pi/2$. Then
$0\le\varphi\le\pi/2$ and
\[
 \frac{\mathrm dr_l}{\mathrm d\varphi}=\frac{2r_l}{2r_l+m\sin2\varphi}.
\]
The solutions $r_l\left(\varphi\right)$ strictly increase with $l$, by uniqueness
of this scalar differential equation and the initial values
$r_l\left(0\right)=q_l$. Consequently
\[
 \log\frac{\kappa_l}{q_l}
 =\int_0^{\pi/2}\frac{2\dd \varphi}{2r_l+m\sin2\varphi}
\]
strictly decreases with $l$, proving the zero-ratio assertion.
Also
\[
 \log\frac{W\left(q_l\right)}{W\left(\kappa_l\right)}
 =\int_0^{\pi/2}
 \frac{2m\sin2\varphi}{\left(2r_l+m\sin2\varphi\right)^2}\dd \varphi
\]
strictly decreases with $l$. Its first value is
\[
 R:=\frac{W\left(q\right)}{W\left(k\right)}=K\alpha^{m-2}>1.
\]

The successive local maxima $W\left(q_l\right)$ are nondecreasing. To see this,
put $Y_\nu=\sqrt rJ_\nu$,
$v_\nu=1-\left(\nu^2-1/4\right)/r^2$ and
$\mathcal E_\nu=Y_{\nu,r}^2+v_\nu Y_\nu^2$.
Then $\mathcal E_{\nu,r}=2\left(\nu^2-1/4\right)Y_\nu^2/r^3\ge0$,
and $W\left(q_l\right)=N^2\mathcal E_\nu\left(q_l\right)$ with
$N=2^\nu\Gamma\left(\nu+1\right)$.
The other positive critical points of $W$ are the local minima
$\kappa_l$, while $W\left(0\right)=0$. Hence $W\left(x\right)\le W\left(q_l\right)$ whenever
$0\le x\le q_l$.
Since $\alpha\kappa_l<q_l$ for $l\ge2$, we obtain
\[
 \frac{W\left(\alpha\kappa_l\right)}{W\left(\kappa_l\right)}
 \le\frac{W\left(q_l\right)}{W\left(\kappa_l\right)}<R.
\]
Multiplying by $\lambda\alpha^{-m}W\left(\kappa_l\right)$ gives
$\kappa_l^mQ\left(\kappa_l\right)>0$, as required. The proof includes
$n=3$, for which the local maxima are constant and the strictly
decreasing quarter-cycle ratio still supplies the strict conclusion.
\end{proof}

\Needspace{6\baselineskip}
Before estimating the higher-order residues, we determine the sign of the first coefficient at the center.

\begin{lemma}
\label{generalG:first-coefficient-analytic}
For every real $n\ge6$, with the canonical ball parameters and
$m=n-1$, one has
\[
 K-n\lambda+m\lambda^2>0.
\]
Consequently
\[
 a_1=\frac{K-n\lambda+m\lambda^2}{m\left(n+2\right)}>0.
\]
\end{lemma}
\begin{proof}
Two analytic facts already proved independently of the calibration
coefficient signs will be used. First, the strictly decreasing
quarter-cycle energy gives
\[
 R=\frac{W\left(q\right)}{W\left(k\right)}=K\lambda^{m/2-1}>1.
\]
Second, the zero ratios $q_l/\kappa_l$ strictly increase with $l$.
For $d\ge2$, define the inverse fourth-power zero sum
\[
 \sigma_2^{\left(d\right)}
 :=\sum_{l\ge1}j_{d/2-1,l}^{-4}
 =\frac1{2d^2\left(d+2\right)}.
\]
The value follows by comparing the coefficient of $r^4$ in the
logarithm of the Bessel series with that of the canonical product
\cite[Sections 10.2(i) and 10.21(iii)]{DLMF}.
Since $q_1/\kappa_1=\alpha$ and $q_l/\kappa_l>\alpha$ for $l\ge2$,
the first terms in the following comparison agree and every remaining
term is strictly smaller:
\[
 \lambda^2\sigma_2^{\left(n\right)}
 =\sum_{l\ge1}\frac{\alpha^4}{q_l^4}
 <\sum_{l\ge1}\frac1{\kappa_l^4}
 =\sigma_2^{\left(n+2\right)}.
\]
It follows that
\[
 0<\lambda<\frac n{\sqrt{\left(n+2\right)\left(n+4\right)}}
 <\lambda_H:=\frac{n\left(n+3\right)}{\left(n+2\right)\left(n+4\right)}<1.
\]
The second strict comparison uses
$\left(n+3\right)^2=\left(n+2\right)\left(n+4\right)+1$.

Put $a=m/2$ and $x=1-\lambda$. We claim
\[
 \left(1-x\right)^{-a}>1+2ax.
\]
All coefficients of its binomial series are positive. Define the
increasing polynomial in $x>0$
\[
 Q_a\left(x\right)=-1+\sum_{j=2}^{5}
       \frac{\left(a+1\right)\left(a+2\right)\cdots\left(a+j-1\right)}{j!}\,x^{j-1}.
\]
Since $x>x_H=\left(3n+8\right)/\left(\left(n+2\right)\left(n+4\right)\right)$, it suffices to show
$Q_a\left(x_H\right)>0$. The exact polynomial identity is
\[
 Q_a\left(x_H\right)=\frac{P_8\left(n\right)}{1920\left(n+2\right)^4\left(n+4\right)^4},
\]
where
\begin{align*}
 P_8\left(n\right)={}&591n^8+8850n^7+38876n^6-64850n^5\\
 &-1218075n^4-4848304n^3-9493760n^2\\
 &-9483264n-3870720.
\end{align*}
Its nonzero coefficients have exactly one sign change in descending
order. By Descartes' rule it has at most one positive root. Since
\[
 P_8\left(0\right)<0,\qquad P_8\left(6\right)=1751207184>0,
\]
that root lies below $6$, and $P_8\left(n\right)>0$ for all real $n\ge6$.
Thus $Q_a\left(x\right)>0$. Keeping the binomial series through degree five gives
\[
 \left(1-x\right)^{-a}>1+ax+ax\left\{1+Q_a\left(x\right)\right\}>1+2ax,
\]
proving the claim.

Finally $R>1$ implies $K>\lambda^{1-a}$. Therefore
\begin{align*}
 K-n\lambda+m\lambda^2
 &>\lambda\left\{\lambda^{-a}-n+m\lambda\right\}\\
 &=\lambda\left\{\left(1-x\right)^{-a}-1-2ax\right\}>0.
\end{align*}
This proves the lemma.
\end{proof}

\Needspace{6\baselineskip}
The next bounds locate the first retained pole and keep the denominators in the tail estimate positive.

\begin{lemma}
\label{generalG:geometry-analytic}
For every real $n\ge6$, let $q=j_{n/2-1,1}$,
$k=j_{n/2,1}$, $\kappa_2=j_{n/2,2}$, and set
\[
 \alpha=q/k,\quad A=\left(n^2-1\right)/4,\quad
 L=\frac{\pi}{1-\alpha},\quad D=\alpha^2-A/L^2.
\]
Then
\[
 \kappa_2<2k<L,\qquad D>\frac15.
\]
\end{lemma}
\begin{proof}
For $\eta\ge1/2$, the first Dirichlet eigenvalue of
$-\partial_s^2+\left(\eta^2-1/4\right)/s^2$ on $\left(0,1\right)$ is
$j_{\eta,1}^2$. The one-dimensional Wirtinger inequality gives,
when $\eta>1/2$,
\[
 j_{\eta,1}^2>\eta^2-\frac14+\pi^2.
\]
In particular $k^2>A+\pi^2$.
The function $Y\left(r\right)=\sqrt r J_{n/2}\left(r\right)$ satisfies
\[
 Y''+\left(1-\frac A{r^2}\right)Y=0.
\]
On $r\ge k$, the coefficient is strictly greater than or equal to
$\omega^2=1-A/k^2>0$, with strict inequality inside the interval.
Sturm comparison from the zero $k$ therefore gives
\[
 \kappa_2-k<\frac\pi\omega<k.
\]
This proves $\kappa_2<2k$.

For the canonical radial pair $X,p$ on the quarter-cycle from $q$
to $k$, write $X=\rho\cos\theta$ and $p=\rho\sin\theta$.
Here $\pi/2<\theta<\pi$, and
\[
 \theta_r=1-\frac{n-1}{2r}\sin2\theta>1.
\]
Consequently $k-q<\pi/2$ and
$L=\pi k/\left(k-q\right)>2k$.

Apply the same eigenvalue lower bound to the order $n/2-1$:
\[
 q^2>\frac{n^2-4n+3}{4}+\pi^2.
\]
Since $L>2k$ and $\pi^2>9$, it follows that
\[
 D>\frac{q^2-A/4}{k^2}
   >\frac{3n^2-16n+157}{16k^2}.
\]
The numerator is positive for $n\ge6$.
For the upper bound, use the weighted radial Dirichlet problem with
weight $s^{d-1}$, where $d=n+2$. The trial
$v\left(s\right)=\left(1-s^2\right)^2$ has quotient
\[
 \frac{16\int_0^1s^{d+1}\left(1-s^2\right)^2\dd s}
      {\int_0^1s^{d-1}\left(1-s^2\right)^4\dd s}
 =\frac{d\left(d+8\right)}3.
\]
The equality follows by the substitution $z=s^2$ and integration
of the resulting polynomials. The variational bound is strict since
this polynomial does not satisfy the first eigenfunction equation.
Thus $k^2<\left(n+2\right)\left(n+10\right)/3$, and
\[
 D>\frac{3\left(3n^2-16n+157\right)}{16\left(n+2\right)\left(n+10\right)}>\frac15.
\]
The last comparison is the positive polynomial identity
\[
 29n^2-432n+2035
 =29\left(n-\frac{216}{29}\right)^2+\frac{12359}{29}>0.
\]
This proves the lemma.
\end{proof}

\Needspace{6\baselineskip}
These geometric bounds give a uniform upper estimate for the possible negative tail.

\begin{lemma}\label{GtailV:bound}
In every integer dimension $6\le n\le47$, the constant
\[
 V=\frac{m}{2\sqrt D}\left(\frac{\kappa_2}{L}\right)^4
       \left(\frac1L+\frac1{4\pi}\right),\qquad
 L=\frac\pi{1-\alpha},\quad
 D=\alpha^2-\frac A{L^2},\quad A=\frac{n^2-1}{4},
\]
satisfies $V<1/6$. Here $m=n-1$, $\alpha=q/k$ and $\kappa_2$ is the
second positive zero of $J_{n/2}$. \end{lemma}

\begin{proof}
Put $b=k^2$ and $\delta=k-q$. The following elementary geometric
bounds will be used:
\begin{equation}\label{GtailV:geometry}
 b<\frac{\left(n+2\right)\left(n+10\right)}3,\quad
 \lambda=\frac{q^2}{k^2}>\left(\frac{n-2}{n}\right)^2,\quad
 \frac{\kappa_2}{L}
 <\delta\left(\frac1\pi+\frac1{\sqrt{b-A}}\right).
\end{equation}
The first follows from the weighted radial Dirichlet trial in
Lemma~\ref{generalG:geometry-analytic}. For the second, the one-dimensional Rayleigh formula and
Hardy's inequality show that $j_{\nu,1}^2/\nu^2$ strictly decreases
with $\nu>0$: after division by $\nu^2$, its kinetic-minus-Hardy
part has a positive coefficient which strictly decreases. Apply this
at $\nu=\left(n-2\right)/2$ and $\nu+1$.
For the last inequality, the transformed $J_{n/2}$ equation has
potential $1-A/r^2$. Sturm comparison to the constant potential
$1-A/b>0$ on $\left[k,\kappa_2\right]$ gives
$\kappa_2-k<\pi k/\sqrt{b-A}$, and
$L=\pi k/\delta$ gives the result. The bounds on $b-A$ below ensure
the stated positivity.

The quarter-cycle polar equation gives
\[
 \delta<J\left(m/k\right),\qquad
 J\left(a\right)=\int_0^{\pi/2}\frac{\dd \theta}{1+a\sin\theta\cos\theta}.
\]
This is decreasing in $a$. For $0<a<2$, set $z^2=\left(2-a\right)/\left(2+a\right)$.
Elementary integration and the alternating arctangent estimate give
\[
 J\left(a\right)=\frac{4\arctan z}{\sqrt{4-a^2}}
 <\frac4{2+a}\left(1-\frac{z^2}3+\frac{z^4}5\right).
\]
From the first bound in \eqref{GtailV:geometry}, one has respectively
$m/k>3/4$ for $n\ge6$, $m/k>1$ for $n\ge11$, and $m/k>7/6$ for
$n\ge15$. The cleared polynomial inequalities are
\[
 39n^2-204n-132>0,\quad
 2n^2-18n-17>0,\quad
 59n^2-804n-872>0
\]
on their stated ranges. Each is increasing there and positive at its
left endpoint. Thus
\begin{equation*}
 \delta<\frac{13}{10}\ \left(n\ge6\right),\qquad
 \delta<\frac{11}{9}\ \left(n\ge11\right),\qquad
 \delta<\frac{47}{40}\ \left(n\ge15\right),
\end{equation*}
because the three arctangent upper bounds are respectively
$5168/3993<13/10$, $164/135<11/9$, and
$8024/6859<47/40$.

We also need two direct Barta estimates. The positive comparison
function $w=r^s\left(1-r^2\right)$ in the transformed order-$n/2$ equation has
\[
 -\frac{w''}w+\frac A{r^2}
 =\frac{A-s\left(s-1\right)}{r^2}+\frac{4s+2}{1-r^2}.
\]
Its infimum is the square of the sum of the square roots of the two
positive coefficients. For $s=2$, $n\ge6$, this gives
\[
 b>\frac{n^2+31}{4}+\sqrt{10\left(n^2-9\right)}
 >\frac{n^2+10n+35}{4},\qquad b-A>\frac{5n+18}{2}.
\]
The last strict inequality follows from
$\sqrt{10\left(n^2-9\right)}>5n/2+1$, whose squared difference is
$15n^2/4-5n-91>0$ at and beyond $n=6$.
For $s=4$, $n\ge15$, use
$\sqrt{n^2-49}>n-2$ and $\sqrt2>7/5$ to obtain
\begin{equation*}
 b>\frac{5n^2+84n-53}{20},\qquad
 b-A>\frac{21n-12}{5}.
\end{equation*}
These Barta estimates follow directly from the ground-state
factorization on compactly supported tests and form-domain density;
the comparison function need not have the behavior
$r^{\left(n+1\right)/2}$ of the regular solution at zero.

Write $p_0=50/157>1/\pi$. Since $A/b<1$, the preceding bounds give
\[
 D>\left(\frac{n-2}{n}\right)^2-\delta^2p_0^2.
\]
For three of the four required ranges, the following constant bounds
are immediate by substitution at the left endpoint:
\[
\begin{array}{c|c|c|c|c|c|c}
 n&\delta<&\sqrt D>&k>&\sqrt{b-A}>&m\le&V<\\\hline
 11\le n\le14&11/9&7/10&8&6&13&3/20\\
 15\le n\le29&47/40&39/50&43/4&31/4&28&4/25\\
 30\le n\le47&47/40&17/20&37/2&11&46&3/20
\end{array}
\]
To make the last column fully explicit, its three rational upper
bounds before rounding are
\begin{align*}
 &\frac{65}{7}\left\{\frac{11}{9}\left(p_0+1/6\right)\right\}^4
      p_0\left(\frac{11}{72}+\frac14\right)<\frac3{20},\\
 &\frac{700}{39}\left\{\frac{47}{40}\left(p_0+4/31\right)\right\}^4
      p_0\left(\frac{47}{430}+\frac14\right)<\frac4{25},\\
 &\frac{460}{17}\left\{\frac{47}{40}\left(p_0+1/11\right)\right\}^4
      p_0\left(\frac{47}{740}+\frac14\right)<\frac3{20}.
\end{align*}
Each is a direct rational inequality; clearing the displayed positive
denominators verifies it. They are all less than $1/6$.

It remains to consider the real interval $6\le n\le10$.
The Barta bound gives $k>\left(n+5\right)/2$ and
$\sqrt{b-A}>\left(n+13\right)/4$: for the latter square the inequality
$n^2-14n+25<0$, which holds throughout $\left[6,10\right]$.
Moreover
\[
 \sqrt D>\frac{n+7}{25}.
\]
Indeed $f\left(n\right)=\sqrt{\left(1-2/n\right)^2-\left(65/157\right)^2}$ is strictly concave for
$n\ge6$, as follows by differentiating $f=\sqrt{w^2-c^2}$ with
$w=1-2/n$, $w'>0$, $w''<0$. The direct rational comparisons
$f\left(6\right)>13/25$ and $f\left(10\right)>17/25$ give the stated lower chord.
Consequently $V<W\left(n\right)$, where
\[
 W\left(n\right)=\frac{25\left(n-1\right)}{2\left(n+7\right)}
 \left\{\frac{13}{10}\left(p_0+\frac4{n+13}\right)\right\}^4
 p_0\left\{\frac{13}{5\left(n+5\right)}+\frac14\right\}.
\]
This rational function strictly decreases for every $n\ge6$.
For a short check, put
\[
 \mathcal D=\left(n-1\right)\left(n+7\right)\left(n+13\right)\left\{p_0\left(n+13\right)+4\right\}
                   \left(n+5\right)\left(n+77/5\right)>0.
\]
With $v=n-6$, direct logarithmic differentiation gives
\[
 -\mathcal D\frac{W'}W
 =\frac{11209208+53618064v+10245920v^2+645360v^3+13160v^4}{785}>0.
\]
Finally direct substitution gives
\[
 \frac16-W\left(6\right)
 =\frac{183971814833017}{1640916934670745204}>\frac1{9000}>0.
\]
Thus $V<W\left(6\right)<1/6$ also in this last range, completing the proof.
\end{proof}

\begin{figure}[tbp]
\centering
\begin{tikzpicture}[>=stealth,every node/.style={font=\small}]
\draw[->] (0,0)--(11,0) node[right] {$r$};
\foreach \x in {1,2.8,4.5,6.6,8.2,9.6}
  \draw (\x,-.1)--(\x,.1);
\node[below] at (1,-.15) {$k$};
\node[below] at (2.8,-.15) {$\kappa_2$};
\node[below] at (4.5,-.15) {$2k$};
\node[below] at (6.6,-.15) {$L$};
\node[below] at (8.2,-.15) {$\kappa_l$};
\node[below] at (9.6,-.15) {$\cdots$};
\draw[fill=white] (1,0) circle (.055);
\fill (2.8,0) circle (.055);
\node[text width=2cm,text centered] at (1,1.0) {removable\\singularity};
\node[text width=3cm,text centered] at (3.5,1.15) {first retained pole\\positive term};
\node[text width=3.7cm,text centered] at (8.6,1.0) {possible negative tail\\estimated in total};
\draw[->] (1,.65)--(1,.17);
\draw[->] (3.1,.7)--(2.82,.17);
\draw[densely dotted] (6.6,.2)--(6.6,1.5);
\end{tikzpicture}
\caption{The pole decomposition, schematically and not to scale.
The inequalities $\kappa_2<2k<L$ separate the first retained pole
from the region where the residue contribution may change sign.
The argument compares a lower bound for the first contribution with
an upper bound for the entire adverse tail.}
\label{fig:poles}
\end{figure}
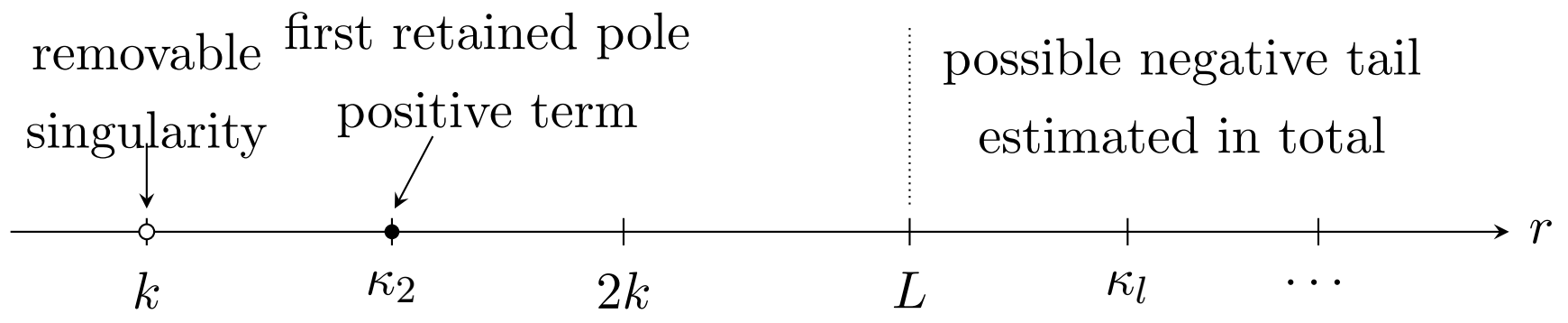

Figure~\ref{fig:poles} shows the separation of the first retained pole
from the tail.

\Needspace{6\baselineskip}
The next criterion reduces positivity of all higher coefficients to
domination of one infinite tail.

\begin{lemma}
\label{generalG:tailcriterion}
Let $n\ge4$, and define
\[
 L=\frac\pi{1-\alpha},\quad A=\frac{m\left(m+2\right)}4,\quad
 D=\alpha^2-A/L^2,\quad
 V=\frac{m}{2\sqrt D}\left(\frac{\kappa_2}{L}\right)^4
                   \left(\frac1L+\frac1{4\pi}\right),
\]
\[
 T_l=3Q\left(\kappa_l\right)+\frac{m\alpha}{\kappa_l}
 X\left(\alpha\kappa_l\right)p\left(\alpha\kappa_l\right),\quad
 U_l=\frac{T_l}{KX\left(\kappa_l\right)^2}.
\]
Suppose $D>0$, $\kappa_2<L$, and $U_2>V$. Then for every $j\ge2$,
\[
 \sum_{l\ge2}\frac{T_l}{X\left(\kappa_l\right)^2\kappa_l^{2j}}>0.
\]
\end{lemma}
\begin{proof}
If $\kappa_l<L$, separation of consecutive positive zeros by more
than $\pi$ gives
$\kappa_{l-1}<\kappa_l-\pi<\alpha\kappa_l<q_l$.
On $\left(\kappa_{l-1},q_l\right)$ the factors $X,p$ have the same sign.
Lemma~\ref{generalG:Qpoles} therefore gives $T_l>0$ in this range.

For the remaining poles put $x=\alpha\kappa_l$,
$Y=\sqrt xJ_{\nu+1}\left(x\right)$, $v=1-A/x^2$, and
$\mathcal E=Y_r^2+vY^2$. Since $x\ge\alpha L>\sqrt A$,
this is a positive quadratic form. The identities
\[
 X\left(x\right)=Nx^{-m/2}\left\{Y_r+mY/\left(2x\right)\right\},\qquad
 p\left(x\right)=Nx^{-m/2}Y
\]
give
\[
 \kappa_l^mT_l=3K\kappa_l^mX\left(\kappa_l\right)^2
 -3N^2\lambda\alpha^{-m}
 \left\{Y_r^2+\left(1-5m^2/\left(12x^2\right)\right)Y^2-mYY_r/\left(3x\right)\right\}.
\]
Relative to $\mathcal E$, the $Y^2$ diagonal correction is
$-m\left(m-3\right)/\left(6x^2\right)\le0$. Dropping that nonpositive correction shows
that the largest relative eigenvalue is at most
\[
 1+\frac{m}{6\sqrt{x^2-A}}.
\]
Moreover $\mathcal E\left(x\right)<\mathcal E\left(\kappa_l\right)$ by the Sonin
derivative, and $N^2\mathcal E\left(\kappa_l\right)
=\kappa_l^mX\left(\kappa_l\right)^2$, because $J_{\nu+1}\left(\kappa_l\right)=0$.
Using $R>1$ from Lemma~\ref{generalG:Qpoles}, this yields
\[
 U_l> -\frac{m}{2R\sqrt{\alpha^2\kappa_l^2-A}}
 >-\frac{m}{2\sqrt D\,\kappa_l}\qquad\left(\kappa_l\ge L\right).
\]
The spacing bound supplies
\[
 \sum_{\kappa_l\ge L}\kappa_l^{-5}
 \le L^{-5}+\frac{L^{-4}}{4\pi}.
\]
Indeed the first term is at most $L^{-5}$, and comparison of each
subsequent term with the integral of the decreasing function $r^{-5}$
over the preceding gap gives the second term.
For every $j\ge2$ and $\kappa_l\ge L>\kappa_2$,
\[
 \frac{\kappa_2^{2j}}{\kappa_l^{2j+1}}
 \le\frac{\kappa_2^4}{\kappa_l^5}.
\]
Keep the positive $l=2$ term, discard all the other positive terms
below $L$, and use the tail bound. The resulting normalized sum is
strictly larger than $U_2-V>0$, proving the claim.
\end{proof}

For the following phase comparison, let $n\ge4$ and assume
$\kappa_2<L=\pi/\left(1-\alpha\right)$. Write
$m=n-1$, $A=\left(n^2-1\right)/4$, $A_\nu=\left(\left(n-2\right)^2-1\right)/4$, and retain
the notation $q=q_1$, $k=\kappa_1$, $\alpha=q/k$, $\kappa=\kappa_2$,
$x=\alpha\kappa$ and $\delta=k-q$. In particular
$k<x<q_2<\kappa$. The first inequality follows from
$\kappa<L$, the zero spacing $\kappa-k>\pi$, and
$\alpha\kappa>\kappa-\pi$.

\Needspace{6\baselineskip}
A quantitative phase comparison supplies the lower bound for the first retained pole.

\begin{lemma}
\label{Gfirstpole:phase}
Let $n\ge4$ and $\kappa<L$. Suppose that
$K_0>k$, $C_0>\kappa$, $d_0>\delta$ and
$B_0<b-A$, where $b=k^2$. Put $p_-=157/50$ and $p_+=22/7$, and define
\[
 s=\frac{A_\nu}{C_0^2},\qquad
 g=p_-\left(1+\frac s2\right),\qquad
 I=\frac{2p_-g}{\left(2K_0+2m/p_-\right)\left(2C_0+2m/p_-\right)},
\]
\[
 D_k=2K_0+\frac{mp_+}4,\quad
 D_c=2C_0+\frac{mp_+}4,\quad
 E=\frac{4mg\left(D_k+D_c\right)}{D_k^2D_c^2},\quad
 c_0=1-\frac{2d_0^2p_+^2}{3B_0},
\]
\[
 S=\frac{mI c_0}{1+mI},\qquad
 u_0=\frac S{1-S},\qquad
 \mathcal U=\frac{u_0+3E}{1+E}.
\]
If $0<s<1$, $0<c_0<1$ and $mI<1$, then $U_2>\mathcal U$.
\end{lemma}
\begin{proof}
Let $r_l\left(\varphi\right)$ be the quarter-phase trajectories from $q_l$ to
$\kappa_l$, as in Lemma~\ref{generalG:Qpoles}. They satisfy
\[
 r_l'=\frac{2r_l}{2r_l+m\sin2\varphi},\qquad
 0\le\varphi\le\frac\pi2.
\]
The right side is increasing in $r_l$. Thus
$r_2-r_1\ge q_2-q_1$. The transformed $J_\nu$ equation on
$\left[q_1,q_2\right]$ has potential at most $1-A_\nu/C_0^2$.
Sturm comparison and $\left(1-s\right)^{-1/2}>1+s/2$ give
\[
 q_2-q_1>\frac\pi{\sqrt{1-A_\nu/C_0^2}}>g.
\]
Let $R_l=W\left(q_l\right)/W\left(\kappa_l\right)$. The exact quarter-phase identities give
\begin{align*}
 I_*:=\log\frac{q_2/\kappa}{q/k}
 &=\int_0^{\pi/2}
 \frac{4\left(r_2-r_1\right)}{\left(2r_1+m\sin2\varphi\right)
                         \left(2r_2+m\sin2\varphi\right)}\dd \varphi,\\
 E_*:=\log\frac{R_1}{R_2}
 &=\int_0^{\pi/2}2m\sin2\varphi
 \left\{\frac1{\left(2r_1+m\sin2\varphi\right)^2}
       -\frac1{\left(2r_2+m\sin2\varphi\right)^2}\right\}\dd \varphi.
\end{align*}
For the first integral, bound $r_1<K_0$, $r_2<C_0$, and
$r_2-r_1>g$. The function
\[
 t\longmapsto\frac1{\left(2K_0+mt\right)\left(2C_0+mt\right)}
\]
is convex for $t\ge0$.
The average of $\sin2\varphi$ on this interval is $2/\pi$.
Jensen's inequality, followed by $p_-<\pi$, yields $I_*>I$.

For the second integral use
\[
 \frac1{a^2}-\frac1{b^2}
   =\frac{\left(b-a\right)\left(a+b\right)}{a^2b^2},\qquad b-a=2\left(r_2-r_1\right)>2g.
\]
The function $\left(a+b\right)/\left(a^2b^2\right)=1/\left(ab^2\right)+1/\left(a^2b\right)$ is decreasing
in each positive variable. After their upper substitution, it is
convex in $t=\sin2\varphi$: each summand is a product of positive
inverse powers of positive affine functions, with nonnegative second
derivative. The measure $\sin2\varphi\dd \varphi$ has total mass one
and mean value $\pi/4$ for $t$. Jensen's inequality and $\pi<p_+$
therefore give $E_*>E$.

Let $\theta$ be the polar phase of $\left(X,p\right)$, and set
$\vartheta=3\pi/2-\theta\left(x\right)\in\left(0,\pi/2\right)$ and $d=q_2-x$.
Since $q_2/x=e^{I_*}$, one has $d>xI$. Run the phase equation
backwards from $q_2$ to $x$. At distance $t$ from $q_2$, its phase
deficit satisfies
\[
 \vartheta_t=1-\frac{m\sin2\vartheta}{2\left(q_2-t\right)}
       \ge1-\frac m x\vartheta.
\]
Integration and $e^{-z}\le\left(1+z\right)^{-1}$ show
\[
 \vartheta\ge\frac{x}{m}\left(1-e^{-md/x}\right)
       \ge\frac{d}{1+md/x}>\frac{xI}{1+mI}.
\]
On the other hand $\vartheta\le d$. The quarter-cycle lengths strictly
increase with $l$, so $\kappa-q_2>k-q=\delta$. Hence
\[
 \vartheta\le d<\delta\frac{\kappa-k}{k}
        <\frac{d_0\pi}{\sqrt{b-A}}
        <\frac{d_0p_+}{\sqrt{B_0}}.
\]
Here the second strict bound is the Sturm spacing estimate for the
order-$n/2$ equation. Using $\sin t\ge t-t^3/6$ gives
\[
 \frac{m\sin2\vartheta}{2x}
 \ge\frac{m\vartheta}{x}\left(1-\frac{2\vartheta^2}{3}\right)
 >\frac{mI}{1+mI}c_0=S.
\]
The positive quantity
\[
 u=\frac{xW'\left(x\right)}{W\left(x\right)}
   =\frac{m\sin2\vartheta/\left(2x\right)}{1-m\sin2\vartheta/\left(2x\right)}
\]
therefore exceeds $u_0$. Its denominator is positive because
$W\left(x\right)>0$.

Finally put $\sigma=W\left(x\right)/\left(R_1W\left(\kappa\right)\right)$. Since
$W\left(x\right)<W\left(q_2\right)=R_2W\left(\kappa\right)$, one has
$0<\sigma<e^{-E_*}<\left(1+E\right)^{-1}$. The exact expression for the
first pole is
\[
 U_2=3-\sigma\left(3-u\right).
\]
The hypotheses imply $0<S<1/2$ and $u_0<1$, hence
$\mathcal U<3$. If $u\ge3$ the conclusion is immediate; if $u<3$,
the preceding upper bound for $\sigma$ gives
$U_2>\left(u_0+3E\right)/\left(1+E\right)$, as claimed.
\end{proof}

\Needspace{6\baselineskip}
We next replace the phase parameters by elementary bounds that hold on whole dimension intervals.

\begin{lemma}
\label{Gfirstpole:envelopes}
The following choices satisfy every hypothesis of
Lemma~\ref{Gfirstpole:phase} throughout the indicated real intervals:
\[
\begin{array}{c|c|c|c|c}
 n&K_0&C_0-K_0&d_0&B_0\\\hline
 {\left[6,10\right]}&3\left(n+5\right)/5&19/5&13/10&\left(5n+18\right)/2\\
 {\left[11,14\right]}&3\left(n+5\right)/5&4&11/9&\left(5n+18\right)/2\\
 {\left[15,29\right]}&\left(7n+38\right)/12&23/5&47/40&\left(21n-12\right)/5\\
 {\left[30,47\right]}&\left(7n+38\right)/12&5&47/40&\left(21n-12\right)/5
\end{array}
\]
In addition $V<\mathcal V$, where
\[
 \mathcal V=\frac{m}{2D_0}
 \left\{d_0\left(50/157+1/H_0\right)\right\}^4\frac{50}{157}
                       \left(\frac{d_0}{k_0}+\frac14\right)
\]
with the following choices, except that on $\left[11,14\right]$ one can take
the simpler bound $\mathcal V=3/20$:
\[
\begin{array}{c|c|c|c}
 n&D_0&k_0&H_0\\\hline
 {\left[6,10\right]}&\left(n+7\right)/25&\left(25n+133\right)/50&\left(7n+104\right)/30\\
 {\left[15,29\right]}&\left(7n+1065\right)/1500&\left(n+6\right)/2&\left(2n+40\right)/9\\
 {\left[30,47\right]}&\left(3n+1355\right)/1700&\left(n+6\right)/2&\left(3n+110\right)/18
\end{array}
\]
\end{lemma}
\begin{proof}
Lemma~\ref{generalG:geometry-analytic} gives $\kappa<L$ for every real
$n\ge6$. The Barta and quarter-phase bounds in the proof of
Lemma~\ref{GtailV:bound} give $b-A>B_0$ and $\delta<d_0$ on
these real intervals. Its Dirichlet trial gives
$b<b_+=\left(n+2\right)\left(n+10\right)/3$. The two upper bounds for $k$ follow from
\[
 \left\{\frac{3\left(n+5\right)}5\right\}^2-b_+
   =\frac{2n^2-30n+175}{75}>0,
 \qquad
 \left\{\frac{7n+38}{12}\right\}^2-b_+
   =\frac{\left(n-22\right)^2}{144}\ge0.
\]
To bound the second zero, let
$f\left(t\right)=t+\pi t/\sqrt{t^2-A}$. It is increasing for $t\ge k$.
Indeed it suffices that $\pi^2A^2<\left(b-A\right)^3$; for $6\le n\le14$
this follows from
\[
 \pi^2A^2<\frac{5n^4}{8}
 \le\frac{70n^3}{8}<\left(\frac{5n}{2}\right)^3<\left(b-A\right)^3,
\]
and for $15\le n\le47$ from
$5n^4/8\le235n^3/8<\left(4n\right)^3<\left(b-A\right)^3$.
Sturm comparison gives $\kappa<f\left(k\right)<f\left(\sqrt{b_+}\right)$.
The function
\[
 R\left(n\right)=\frac{n^2+12n+20}{n^2+48n+83}
\]
is increasing, since the numerator of its derivative is
$36n^2+126n+36$. Thus
\[
 f\left(\sqrt{b_+}\right)<K_0+2p_+\sqrt{R\left(n\right)}<C_0.
\]
For the four intervals, the exact positive differences at their right
endpoints between $\left(C_0-K_0\right)^2$ and $4p_+^2R\left(n\right)$ are respectively
\[
 \frac{37269}{270725},\quad \frac{720}{15533},\quad
 \frac{126453}{236425},\quad \frac{279}{379}.
\]

Clearly $0<s<1$. Also $0<c_0<1$: its smallest coarse comparison
is the low-dimensional endpoint $n=6$, where
$2\left(13/10\right)^2p_+^2/\left(3\cdot24\right)<1$.
For the remaining hypothesis, $g<3\pi/2$ and
$K_0\ge\left(n+5\right)/2$, $C_0>\left(n+11\right)/2$ imply
\[
 mI<\frac{30\left(n-1\right)}{\left(n+5\right)\left(n+11\right)}<1,
\]
because $\left(n+5\right)\left(n+11\right)-30\left(n-1\right)=\left(n-7\right)^2+36>0$.

For the upper bound on $V$, it remains to check
$\sqrt D>D_0$, $k>k_0$, and $\sqrt{b-A}>H_0$.
The Hardy estimate and the quarter-phase bound give
\[
 \sqrt D>\sqrt{\left(1-2/n\right)^2-\left(d_0\cdot50/157\right)^2}.
\]
For fixed $d_0$ the right side is concave in $n$ where it is
positive, so its endpoint lower chords apply. The squared positive
differences at $n=6,10$ from $\left(13/25\right)^2$ and
$\left(17/25\right)^2$ are, respectively,
\[
 \frac{365746}{138650625},\qquad \frac{95414}{15405625}.
\]
For the high intervals, use endpoints $15,30,47$ and chord values
$39/50,17/20,22/25$. The squared positive differences are
\[
 \frac{5951131}{2218410000},\quad
 \frac{76159}{8873640},\quad
 \frac{1239214471}{544496410000}.
\]
The printed $D_0$ are exactly these chords, using $\left[15,30\right]$ for
the third row. The Barta lower bound for $b$ gives $k>k_0$;
in the low interval the squared difference is
$\left(2093-200n\right)/1250>0$, and in the high intervals it is
$\left(24n-233\right)/20>0$. Finally the differences $B_0-H_0^2$ are concave
quadratics. At the endpoints of $\left[6,10\right]$ they equal $71/225,9/25$;
at those of $\left[15,30\right]$ they equal $43/405,58/405$; and at those of
$\left[30,47\right]$ they equal $58/405,179/324$. Hence all are positive.
Substituting these bounds into the definition of $V$ proves
$V<\mathcal V$. The bound $V<3/20$ on $\left[11,14\right]$ was proved directly
in Lemma~\ref{GtailV:bound}.
\end{proof}

\Needspace{6\baselineskip}
We first compare the lower and upper bounds on the two lower-dimensional parameter intervals.

\begin{lemma}\label{Gfirstpole:low-algebra}
For the parameters and scalar functions $\mathcal U,\mathcal V$ in
Lemma~\ref{Gfirstpole:envelopes}, one has
\[
 \mathcal U>\mathcal V
 \quad\text{throughout }\left[6,10\right]\text{ and throughout }\left[11,14\right].
\]
\end{lemma}

\begin{proof}
The following polynomial recipe specifies both exact expansions used
in the proof. Put
\[
 \mathsf p=157/50,\quad \mathsf P=22/7,\quad
 \mathsf r=50/157,\quad m=n-1,\quad K_0=3\left(n+5\right)/5.
\]
On the first interval take $C_0=K_0+19/5$ and $d_0=13/10$; on the
second take $C_0=K_0+4$ and $d_0=11/9$. Define the polynomials
\begin{align*}
 H&=2C_0^2+\left\{\left(n-2\right)^2-1\right\}/4,\\
 I_N&=\mathsf p^4H,&
 I_D&=4C_0^2\left(\mathsf pK_0+m\right)\left(\mathsf pC_0+m\right),\\
 D_k&=2K_0+m\mathsf P/4,&D_c&=2C_0+m\mathsf P/4,\\
 E_N&=2m\mathsf pH\left(D_k+D_c\right),&
 E_D&=C_0^2D_k^2D_c^2,\\
 C_D&=3\left(5n+18\right),&C_N&=C_D-4d_0^2\mathsf P^2,\\
 u_N&=mI_NC_N,&u_D&=I_DC_D+mI_N\left(C_D-C_N\right),\\
 U_N&=u_NE_D+3E_Nu_D,&U_D&=u_D\left(E_D+E_N\right).
\end{align*}
The formulas in Lemma~\ref{Gfirstpole:phase} give exactly
$I=I_N/I_D$, $E=E_N/E_D$, $u_0=u_N/u_D$, and
$\mathcal U=U_N/U_D$. All displayed denominators are positive.

On $\left[6,10\right]$, let
\begin{align*}
 V_N={}&\frac{25m d_0^4\mathsf r}{8}
 \left\{\mathsf r\left(7n+104\right)+30\right\}^4\left(25n+133+200d_0\right),\\
 V_D={}&\left(n+7\right)\left(7n+104\right)^4\left(25n+133\right).
\end{align*}
Then $\mathcal V=V_N/V_D$ exactly. On $\left[11,14\right]$, take simply
$V_N=3$, $V_D=20$. In either case put
\[
 Q\left(n\right)=U_NV_D-V_NU_D.
\]
Thus it suffices to prove that this explicitly defined polynomial is
positive on the corresponding interval.

For the first interval write
$Q\left(6+v\right)=10^{25}\sum_{j=0}^{17}q_jv^j$.
Multiplication of the displayed polynomial factors gives the following
strict rational lower bounds. Powers of ten in this table are exact.
\[
\begin{array}{c|rrrrrrrrr}
 j&0&1&2&3&4&5&6&7&8\\\hline
 q_j>&1&3&3&1&3/5&13/100&11/500&7/2500&13/50000
\end{array}
\]
\[
\begin{array}{c|rrrrrrrrr}
 j&9&10&11&12&13&14&15&16&17\\\hline
 q_j>&19\cdot10^{-6}&10^{-6}&4\cdot10^{-8}&10^{-9}&
 3\cdot10^{-11}&4\cdot10^{-13}&2\cdot10^{-15}&
 -7\cdot10^{-18}&-2\cdot10^{-19}
\end{array}
\]
Every term through degree fourteen is nonnegative for $v\ge0$,
with a positive constant. On $0\le v\le4$, the remaining terms obey
\[
 q_{15}v^{15}+q_{16}v^{16}+q_{17}v^{17}
 \ge v^{15}\left\{2\cdot10^{-15}-28\cdot10^{-18}-32\cdot10^{-19}\right\}
 \ge0.
\]
Consequently $Q>0$ on the entire first interval.

For the second interval write
$Q\left(11+v\right)=10^{16}\sum_{j=0}^{11}q_jv^j$.
The same factor recipe gives these strict lower bounds:
\[
\begin{array}{c|rrrrrr}
 j&0&1&2&3&4&5\\\hline
 q_j>&5&3&4/5&11/100&9/1000&3\cdot10^{-4}
\end{array}
\]
\[
\begin{array}{c|rrrrrr}
 j&6&7&8&9&10&11\\\hline
 q_j>&-9\cdot10^{-6}&-18\cdot10^{-7}&-10^{-7}&
 -3\cdot10^{-9}&-5\cdot10^{-11}&-4\cdot10^{-13}
\end{array}
\]
For $0\le v\le3$, group all terms from degree five onwards. Their
coefficient after factoring $v^5$ is greater than
\begin{align*}
 3\cdot10^{-4}
 &-3\left(9\cdot10^{-6}\right)-9\left(18\cdot10^{-7}\right)-27\cdot10^{-7}\\
 &-81\left(3\cdot10^{-9}\right)-243\left(5\cdot10^{-11}\right)-729\left(4\cdot10^{-13}\right)
 >\frac1{4000}>0.
\end{align*}
The lower-degree terms are nonnegative with a positive constant.
Thus $Q>0$ on the second interval as well. This proves the lemma.
\end{proof}

\Needspace{6\baselineskip}
The same comparison on the two remaining intervals completes the parameter estimates through dimension forty-seven.

\begin{lemma}\label{Gfirstpole:high-algebra}
Use the quantities $\mathcal U,\mathcal V$ in
Lemmas~\ref{Gfirstpole:phase}--\ref{Gfirstpole:envelopes}.
Then
\[
 \mathcal U>\mathcal V
 \qquad\left(15\le n\le29\ \hbox{or}\ 30\le n\le47\right).
\]
\end{lemma}

\begin{proof}
We give factored polynomials and small coefficient bounds for the
whole intervals. Put
\[
 \begin{gathered}
 p=\frac{157}{50},\quad z=\frac1p,\quad P=\frac{22}{7},\quad
 \delta=\frac{47}{40},\quad m=n-1,\\
 K_0=\frac{7n+38}{12},\quad B=\frac{21n-12}{5},
 \quad b_*=\frac{2\delta^2P^2}{3}.
 \end{gathered}
\]
Choose the remaining affine polynomials by the following table:
\[
\begin{array}{c|c|c|c|c}
 \hbox{interval}& C_0-K_0&D_0&k_0&H_0\\ \hline
 {\left[15,29\right]}&23/5&\left(7n+1065\right)/1500&\left(n+6\right)/2&\left(2n+40\right)/9\\
 {\left[30,47\right]}&5&\left(3n+1355\right)/1700&\left(n+6\right)/2&\left(3n+110\right)/18
\end{array}
\]
Define
\begin{gather*}
 A_\nu=\frac{\left(n-2\right)^2-1}{4},\qquad
 R=2C_0^2+A_\nu,\qquad
 B_k=K_0+\frac mp,\quad B_c=C_0+\frac mp,\\
 D_k=2K_0+\frac{mP}{4},\qquad
 D_c=2C_0+\frac{mP}{4},\\
 I_N=p^2R,\qquad I_D=4C_0^2B_kB_c,\qquad
 E_N=2mpR\left(D_k+D_c\right),\qquad E_D=C_0^2D_k^2D_c^2,\\
 J=mI_N\left(B-b_*\right),\qquad Q=I_DB+mI_Nb_*.
\end{gather*}
Every denominator and every displayed factor is positive on the
stated intervals. Direct substitution in the phase formulas gives
\[
 I=\frac{I_N}{I_D},\qquad E=\frac{E_N}{E_D},\qquad
 u_0=\frac JQ,\qquad
 \mathcal U=\frac{JE_D+3E_NQ}{Q\left(E_D+E_N\right)}.
\]
The following additional polynomials remove a common positive
factor $m$:
\begin{gather*}
 Z=Q\left(E_D+E_N\right),\\
 H=pR\left\{p\left(B-b_*\right)E_D+6\left(D_k+D_c\right)Q\right\},\\
 N=\delta^4z\left(zH_0+1\right)^4\left(4\delta+k_0\right),\qquad
 T=8D_0H_0^4k_0,\qquad
 \Pi=HT-NZ.
\end{gather*}
Indeed $JE_D+3E_NQ=mH$, so that the exact identity is
\begin{equation}\label{Gfirstpole:high-clearing}
 \mathcal U-\mathcal V=\frac{m\Pi\left(n\right)}{ZT}.
\end{equation}
In particular $\Pi$ is a polynomial of degree $16$ with completely
specified rational coefficients.

Here are the two endpoint checks used only to normalize this
polynomial; all quantities in them are the rational functions just
displayed:
\[
\begin{array}{c|c|c}
 n&\mathcal U\left(n\right)>&\mathcal V\left(n\right)<\\ \hline
 15&1/6&2/25\\
 30&7/50&1/10
\end{array}
\]
They follow by multiplication of positive rational denominators.
Thus $\Pi\left(15\right)>0$ in the first block, and $\Pi\left(30\right)>0$ in the second.

For a block $\left[\ell,r\right]$, set $t=\left(n-\ell\right)/\left(r-\ell\right)$ and write
\[
 \frac{\Pi\left(\ell+\left(r-\ell\right)t\right)}{\Pi\left(\ell\right)}
   =\sum_{j=0}^{16}q_jt^j,\qquad
 q_j=\frac{\left(r-\ell\right)^j\Pi^{\left(j\right)}\left(\ell\right)}{j!\,\Pi\left(\ell\right)}.
\]
The following are strict lower bounds for $q_j$; in both rows $q_0=1$.
They are obtained simply by expanding the factored polynomial $\Pi$
above and using the displayed coefficient formula.
\[
\begin{array}{c|rrrrrrrr}
 &q_1&q_2&q_3&q_4&q_5&q_6&q_7&q_8\\ \hline
 {\left[15,29\right]}&7&29&64&97&103&79&43&16\\
 {\left[30,47\right]}&4&10&12&7&-1&-5&-5&-3
\end{array}
\]
{\setlength{\arraycolsep}{4pt}
\[
\begin{array}{c|rrrrrrrr}
 &q_9&q_{10}&q_{11}&q_{12}&q_{13}&q_{14}&q_{15}&q_{16}\\ \hline
 {\left[15,29\right]}&3&-284/1000&-474/1000&-178/1000&
                   -39/1000&-6/1000&-1/1000&-1/1000\\
 {\left[30,47\right]}&-2&-1&-1&-1&-1&-1&-1&-1
\end{array}
\]
}

For the first interval and $0\le t\le1$, discard the nine positive
nonconstant terms and use $t^j\le1$ on the remaining terms. This gives
\[
 \frac{\Pi\left(15+14t\right)}{\Pi\left(15\right)}
 \ge1-\frac{284+474+178+39+6+1+1}{1000}
 =\frac{17}{1000}>0.
\]
For the second interval, $t^j\ge t^5$ when $1\le j\le4$, whereas
$t^j\le t^5$ when $j\ge5$. The positive lower bounds sum to
$4+10+12+7=33$; the absolute values of all negative lower bounds
sum to $1+5+5+3+2+7=23$. Therefore
\[
 \frac{\Pi\left(30+17t\right)}{\Pi\left(30\right)}
 \ge1+\left(33-23\right)t^5\ge1>0.
\]
Equation~\eqref{Gfirstpole:high-clearing} now proves both continuous
blocks.
\end{proof}

\Needspace{6\baselineskip}
The preceding estimates now compare the first pole with the entire possible negative tail.

\begin{lemma}
\label{generalG:parameters}
For every integer $6\le n\le47$, one has $U_2>V$.
\end{lemma}
\begin{proof}
Choose the row of Lemma~\ref{Gfirstpole:envelopes} containing $n$.
Lemma~\ref{Gfirstpole:phase} gives $U_2>\mathcal U$ and
Lemma~\ref{Gfirstpole:envelopes} gives $V<\mathcal V$.
Lemmas~\ref{Gfirstpole:low-algebra} and~\ref{Gfirstpole:high-algebra}
give $\mathcal U>\mathcal V$ in all four rows. Hence $U_2>V$.
\end{proof}

\Needspace{6\baselineskip}
The residue formula and this comparison establish all coefficient signs in the higher-dimensional range.

\begin{theorem}
\label{generalG:theorem}
For every integer $6\le n\le47$ the calibration has a convergent
expansion on $0\le r<\kappa_2$,
\[
 G\left(r\right)=\frac{n\lambda-K}{m}-\sum_{j\ge1}a_jr^{2j},\qquad a_j>0.
\]
In particular on $0<r<k$,
\[
 0<G<K,\qquad G_r<0,\qquad G_{rr}<0,
 \qquad G_r\left(k\right)=-\frac{2K}{3k}.
\]
\end{theorem}
\begin{proof}
The points $0,\pm k$ are removable. At $k$ one has
$Q=Q_r=0$, $Q_{rr}=2mKc^2/k^2$, and hence $G\left(k\right)=0$.
For $j\ge2$ the paired residue identity is
\[
 a_j=\sum_{l\ge2}\frac{2}{mX\left(\kappa_l\right)^2\kappa_l^{2j}}
 \left\{\left(2j-1\right)Q\left(\kappa_l\right)
      +\frac{m\alpha}{\kappa_l}X\left(\alpha\kappa_l\right)p\left(\alpha\kappa_l\right)\right\}.
\]
To check the residue formula,
the double- and simple-pole coefficients of $K-G$ at $r=\kappa$ are
\[
 A_\kappa=\frac{\kappa^2Q\left(\kappa\right)}{mX\left(\kappa\right)^2},\qquad
 D_\kappa=\frac{\kappa^2Q_r\left(\kappa\right)+\left(m+2\right)\kappa Q\left(\kappa\right)}
                         {mX\left(\kappa\right)^2}.
\]
The paired residue at $\pm\kappa$ is
$2\left\{\left(2j+1\right)A_\kappa/\kappa^{2j+2}-D_\kappa/\kappa^{2j+1}\right\}$.
Substitution gives the stated formula, using
$p_r\left(\kappa\right)=X\left(\kappa\right)$,
$p_{rr}\left(\kappa\right)/p_r\left(\kappa\right)=-m/\kappa$ and
\[
 Q_r\left(\kappa\right)=-mQ\left(\kappa\right)/\kappa
                  -m\alpha X\left(\alpha\kappa\right)p\left(\alpha\kappa\right)/\kappa^2.
\]
To justify the identity, integrate $G\left(r\right)/r^{2j+1}$ around squares
whose vertical sides are
\[
 \Re r=\pm\left(N+\left(n+1\right)/4\right)\pi
\]
and whose horizontal sides have the same absolute height.
The Bessel asymptotics in \cite[Section 10.17]{DLMF} give $\left|p\left(r\right)\right|\ge c_0\left|r\right|^{-m/2}
e^{\left|\Im r\right|}$ and $G=O\left(\left|r\right|^2\right)$ there. The contour integral tends
to zero for $j\ge2$. The residue series converges absolutely
because $Q\left(\kappa\right)=O\left(\kappa^{-m}\right)$,
$X\left(\kappa\right)^2\asymp\kappa^{-m}$ and
$X\left(\alpha\kappa\right)p\left(\alpha\kappa\right)=O\left(\kappa^{-m}\right)$.
Every brace is $T_l+2\left(j-2\right)Q\left(\kappa_l\right)$.
Lemma~\ref{generalG:parameters} gives $U_2>V$.
Together with Lemmas~\ref{generalG:Qpoles},
\ref{generalG:geometry-analytic}, and~\ref{generalG:tailcriterion},
this proves $a_j>0$. The center expansion supplies
\[
 a_1=\frac{K-n\lambda+m\lambda^2}{m\left(n+2\right)}>0
\]
by Lemma~\ref{generalG:first-coefficient-analytic}, independently of
Lemma~\ref{generalG:parameters}. The positivity of the $a_j$ and $G\left(k\right)=0$ give
$G>0$, $G_r<0$, $G_{rr}<0$ on $\left(0,k\right)$.

The endpoint expansion in Section~\ref{sec:radial}, derived directly
from the radial equations, gives $G_r\left(k\right)=-2K/\left(3k\right)$.
Since $k<\kappa_2$, the coefficient series may be differentiated at $k$,
and therefore
\[
 \sum_{j\ge1}j a_jk^{2j}=-\frac{kG_r\left(k\right)}2=\frac K3.
\]
Together with $G\left(k\right)=0$, this gives, for $0<r<k$,
\[
 0<G\left(r\right)<\sum_{j\ge1}a_jk^{2j}<\frac K3<K.
\]
The penultimate inequality is strict because $a_2>0$.
\end{proof}

\subsection{The coefficients in dimensions four and five}

We now establish the coefficient signs in dimensions four and five
using the same residue formula and tail criterion. The phase lemma
applies once the condition $\kappa_2<L$ has been checked.

\Needspace{6\baselineskip}
For dimensions four and five, the same argument uses the following explicit root bounds.

\begin{lemma}\label{Glow:roots}
Write $q=j_{n/2-1,1}$, $k=j_{n/2,1}$ and $A=\left(n^2-1\right)/4$. The data
\[
\begin{array}{c|c|c|c|c|c}
 n&q_-&k_-&k_+&C_0&d_D\\\hline
 4&3.83&5.13&5.14&9&18/25\\
 5&4.49&5.76&5.77&10&3/4
\end{array}
\]
satisfy $q>q_-$, $k_-<k<k_+$, $q<3.84$ for $n=4$, and $q<4.50$
for $n=5$. With $L=\pi/\left(1-q/k\right)$ and $D=\left(q/k\right)^2-A/L^2$,
\[
 \kappa_2<C_0<L,\qquad D>d_D^2>0.
\]
\end{lemma}
\begin{proof}
For $d\ge3$ put
\[
 P_{d,N}\left(r\right)=\sum_{j=0}^N\frac{\left(-1\right)^jr^{2j}}{4^j j!\left(d/2\right)_j}.
\]
Here $\left(a\right)_j=a\left(a+1\right)\cdots\left(a+j-1\right)$ for $j\ge1$ and $\left(a\right)_0=1$
is the rising factorial. The decreasing alternating tail gives $P_{d,9}<\psi_d<P_{d,8}$
at all arguments in the following exact rational table:
\[
\begin{array}{c|c|c|c|c}
 d&a&b&P_{d,9}\left(a\right)>&P_{d,8}\left(b\right)<\\\hline
 4&3.83&3.84&1/5000&-1/5000\\
 5&4.49&4.50&1/5000&-1/5000\\
 6&5.13&5.14&1/5000&-1/5000\\
 7&5.76&5.77&1/5000&-1/5000
\end{array}
\]
Each comparison is specified by substituting into the displayed finite
sum and clearing its positive denominator. The initial sign is positive,
so these positive-to-negative crossings have odd index. Consecutive
positive zeros are separated by more than $\pi$; a third zero is greater
than $2\pi$, whereas all the table endpoints are less than $6<2\pi$.
Thus the crossings identify first zeros. Use $p_-=157/50<\pi<p_+=22/7$,
as established at the start of Section~\ref{sec:coefficients}.

Put $B_0=k_-^2-A$ and $L_0=p_-/\left(1-q_-/k_+\right)$. Sturm comparison gives
\[
 \kappa_2<k+\frac{\pi k}{\sqrt{k^2-A}}
 <k_++\frac{p_+k_-}{\sqrt{k_-^2-A}}<C_0.
\]
The last comparison follows from the following positive values of
$\left(C_0-k_+\right)^2B_0-p_+^2k_-^2$, for $n=4,5$ respectively:
\[
 \frac{93456794469}{1225000000},\qquad
 \frac{24281435853}{153125000}.
\]
Also $L>L_0$, whose respective values are $40349/3275$ and $90589/6400$;
each exceeds $C_0$. Finally $D>\left(q_-/k_+\right)^2-A/L_0^2>d_D^2$; the two
positive margins in the last comparison are
\[
 \frac{24670408577}{2035052251250},\qquad
 \frac{1718746495}{131301870736}.
\]
This proves all assertions.
\end{proof}

\Needspace{6\baselineskip}
These root bounds retain the strict first-pole dominance in dimensions four and five.

\begin{lemma}\label{Glow:firstpole}
For both dimensions, $U_2>V$ in Lemma~\ref{generalG:tailcriterion}.
\end{lemma}
\begin{proof}
Apply Lemma~\ref{Gfirstpole:phase} with $K_0=k_+$, $C_0$ from
Lemma~\ref{Glow:roots}, $d_0=k_+-q_-$ and $B_0=k_-^2-A$.
Substitution in its displayed rational formulas gives
\[
 mI>6/25,\qquad c_0>49/100,\qquad E>1/61.
\]
For these data, $0<s<1$, $K_0\ge\left(n+5\right)/2$ and
$C_0>\left(n+11\right)/2$. Since $g<3\pi/2$, the formula for $I$ gives
$mI<30\left(n-1\right)/\left[\left(n+5\right)\left(n+11\right)\right]<1$.
Also $c_0<1$ by its definition, so every hypothesis of the phase lemma
holds. In particular $0<S<1/2$, and the phase estimate gives
\[
 U_2>\mathcal U>108/713>3/25.
\]
Here $S=mIc_0/\left(1+mI\right)$, $u_0=S/\left(1-S\right)$ and
$\mathcal U=\left(u_0+3E\right)/\left(1+E\right)$ increase in the three lower inputs;
for the displayed data $u_0<1<3$.

The tail is smaller than
\[
 \mathcal W_n=\frac{m}{2d_D}\left(\frac{C_0}{L_0}\right)^4
                      \left(\frac1{L_0}+\frac1{4p_-}\right).
\]
The exact positive differences $11/100-\mathcal W_n$ for $n=4,5$ are
\[
 \frac{312567793029496310535353}{21389167149811660795349800},\qquad
 \frac{6206044125316845926731439}{610066710957425872061194900}.
\]
Thus $U_2>3/25>11/100>V$, as required.
\end{proof}

\Needspace{6\baselineskip}
The coefficient signs in these two dimensions follow from the same residue structure.

\begin{theorem}\label{Glow:theorem}
The conclusions of Theorem~\ref{generalG:theorem}, including every
coefficient sign, $0<G<K$ and $G_r\left(k\right)=-2K/\left(3k\right)$, also hold for $n=4,5$.
\end{theorem}
\begin{proof}
The residue identity, removable singularities, contour limit and
absolute convergence in that theorem apply unchanged. Lemmas
\ref{generalG:Qpoles}, \ref{generalG:tailcriterion} and
\ref{Glow:firstpole} give $a_j>0$ for all $j\ge2$.
For $a_1$, the independent quarter-cycle comparison gives
$K>\lambda^{1-m/2}$. The root bounds give $\lambda<9/16$ for $n=4$
and $\lambda<16/25$ for $n=5$. Hence
\[
 \begin{array}{ll}
 K>\lambda^{-1/2}>4/3\ge4\lambda-3\lambda^2,&n=4,\\
 K>\lambda^{-1}>25/16\ge5\lambda-4\lambda^2,&n=5.
 \end{array}
\]
Thus $a_1=\left(K-n\lambda+m\lambda^2\right)/\left[m\left(n+2\right)\right]>0$ in both cases.
The endpoint expansion and the positive-coefficient argument at the end
of that proof give all remaining conclusions, including $G<K$.
\end{proof}

\section{Radial bounds and determinant estimates}\label{sec:multipliers}
The coefficient signs from Section~\ref{sec:coefficients} and
Appendix~\ref{app:planar-coefficients} reduce the matrix problem to bounds for the radial multiplier and the normalized
Dirichlet profile. We first prove the product and multiplier estimates,
then establish the scalar exponential inequalities used in the
determinant argument.

\subsection{Product expansions and radial bounds}\label{sec:radial-bounds}

For $n\ge2$ let $x=r^2$, $b=k^2$, and $F\left(x\right)=r\phi\left(r\right)/p\left(r\right)$.
The zero-ratio comparison proved in Appendix~\ref{app:planar-coefficients}
for $n=2$, and in Lemma~\ref{generalG:Qpoles} for $n\ge3$, yields
\[
 F\left(x\right)=n\prod_{l\ge2}
 \frac{1-x/\left(q_l/\alpha\right)^2}{1-x/\kappa_l^2},\qquad
 \left(q_l/\alpha\right)^2>\kappa_l^2.
\]
Each factor has strictly positive Taylor coefficients after its
constant term. The products converge uniformly on compact subsets
of $\left|x\right|<\kappa_2^2$. Thus $F$ has strictly positive Taylor coefficients there. In
particular it is absolutely monotone: every derivative is nonnegative
on $\left[0,b\right]$.

We next examine the coefficients of the radial parameter.
Put $x=r^2$, $b=k^2$, $p=rh\left(x\right)$ and $T\left(x\right)=rt=u/h$.
Let $\kappa_l$ be the positive zeros of $p/r$, with $\kappa_1=k$.

\Needspace{6\baselineskip}
The sign of this series will control the radial multiplier through the logarithmic derivative of $T$.

\begin{proposition}\label{prop:radial-series}
In every dimension $n\ge2$,
\[
 T\left(x\right)=n\left(1+c\right)-\sum_{j\ge1}s_jx^j,\qquad s_j>0,
 \qquad \left|x\right|<\kappa_2^2.
\]
Consequently $T'<0$ and $T''<0$ on $0<x<b$.
In fact every positive-order derivative of $T$ is strictly negative
on this interval.
\end{proposition}
\begin{proof}
Write
\[
 h\left(x\right)=\frac1n\left(1-x/b\right)P\left(x\right),\qquad
 P\left(x\right)=\prod_{l\ge2}\left(1-x/\kappa_l^2\right).
\]
The radial equation gives $X=nh+2xh'$. Evaluating at $b$ proves
\[
 nc=2P\left(b\right).
\]
Hence
\[
 T=n+2x\frac{h'}h+\frac ch.
\]
For the coefficient computation normalize $a_l=b/\kappa_l^2$,
so $a_1=1$ and $0<a_l<1$ for $l\ge2$, and put
$P_0=\prod_{l\ge2}\left(1-a_l\right)>0$.
Let $p_j=\sum_{l\ge1}a_l^j$. Define the complete homogeneous
symmetric functions $h_j$ by
\[
 \prod_{l\ge1}\left(1-a_l z\right)^{-1}=\sum_{j\ge0}h_jz^j.
\]
Thus $h_j$ is the sum of all degree-$j$ monomials with nondecreasing
indices, with $h_0=1$. The absolutely convergent product expansions give,
for every $j\ge1$,
\[
 b^j\left[x^j\right]T=2\left(P_0h_j-p_j\right).
\]
Since $a_1=1$,
\[
 h_j\left(a_1,a_2,\ldots\right)=\sum_{i=0}^j h_i\left(a_2,a_3,\ldots\right)
 <\sum_{i\ge0}h_i\left(a_2,a_3,\ldots\right)=\frac1{P_0}.
\]
Also $p_j>1$, so every positive-order coefficient of $T$ is
strictly negative. The singularity at $x=b$ is removable because
$u$ has a double zero while $h$ has a simple zero there; the next
possible poles are at $\kappa_2^2$. Thus the Taylor expansion has
the stated disk of convergence. The derivative signs follow termwise.
\end{proof}

Let $x=r^2$, $b=k^2$, $p=rh\left(x\right)$, and write the Bessel product as
\[
 h\left(x\right)=\frac1n\left(1-x/b\right)P\left(x\right),\qquad
 P\left(x\right)=\prod_{l\ge2}\left(1-x/\kappa_l^2\right).
\]
Set $T\left(x\right)=rt=u/h$ and $\mathcal C=up/S$.

\Needspace{6\baselineskip}
The product representation yields the lower multiplier bound used in the coupling estimate.

\begin{lemma}\label{generalClower:main}
On $0<r<k$, in every dimension $n\ge2$,
\[
 \frac{\mathcal C}{r}>\frac{k^2-r^2}{k^2+r^2}.
\]
Moreover $T\left(x\right)/\left(b-x\right)$ has nonnegative Taylor coefficients, and its
derivative is strictly positive on $\left[0,b\right]$.
\end{lemma}
\begin{proof}
Since $u\left(r\right)=\int_r^kp\left(s\right)\dd s$, substitution of the product gives
\[
 R\left(x\right):=\frac{T\left(x\right)}{b-x}
 =\frac12\int_0^1\left(1-v\right)
 \frac{P\left(x+\left(b-x\right)v\right)}{P\left(x\right)}\dd v.
\]
Each factor of the quotient is
\[
 1-\frac{v\left(b-x\right)}{\kappa_l^2-x}
 =\left(1-v\right)+\frac{v\left(\kappa_l^2-b\right)}{\kappa_l^2-x}.
\]
It has nonnegative Taylor coefficients and strictly positive
$x$ derivative when $0<v<1$. Uniform convergence on compact
subintervals of $x<\kappa_2^2$ justifies multiplication and
differentiation of the products, and then integration. Thus $R'>0$.

Direct radial differentiation gives
\[
 \frac{\mathcal C}{r}=
 \frac1{1-2xT'/T},\qquad
 -\frac{2xT'}T=\frac{2x}{b-x}-\frac{2xR'}R
 <\frac{2x}{b-x}.
\]
The first denominator is positive because $S>0$. Taking reciprocals
proves the strict lower bound.
\end{proof}

\Needspace{6\baselineskip}
The lower multiplier bound gives the upper estimate for the coupling coefficient.

\begin{corollary}
\label{generalClower:d}
Suppose $G=\sum_{j\ge1}B_j\left(1-\left(r/k\right)^{2j}\right)$ with $B_j>0$ and
$G_r\left(k\right)=-2K/\left(3k\right)$. Then, in every such dimension,
\[
 d=2G+\mathcal CG_r
 <\frac{2K}{3}\frac{1-\left(r/k\right)^2}{1+\left(r/k\right)^2}.
\]
\end{corollary}
\begin{proof}
Put $y=r^2/k^2$. The endpoint derivative gives
$\sum_{j\ge1}jB_j=K/3$, while the preceding lemma gives
\[
 \frac d2<\sum_{j\ge1}jB_j
 \left\{\frac{1-y^j}{j}-\frac{1-y}{1+y}y^j\right\}.
\]
For every integer $j\ge1$,
\[
 \left(1+y\right)\left(1-y^j\right)\le j\left(1-y\right)\left(1+y^j\right).
\]
Indeed summing $y^i+y^{j-i}\le1+y^j$ over $1\le i<j$ gives
the equivalent geometric-sum inequality. Hence each brace is at
most $\left(1-y\right)/\left(1+y\right)$. Substitution of the endpoint sum proves the claim.
\end{proof}

Let $x=r^2$, $b=k^2$, $y=x/b$, and use
$R\left(x\right)=T\left(x\right)/\left(b-x\right)$ from the lower multiplier bound.

\Needspace{6\baselineskip}
Logarithmic convexity of $R$ gives a complementary upper multiplier bound through dimension twenty-two.

\begin{proposition}\label{generalCupper:main}
In every dimension $2\le n\le22$, throughout $0<r<k$,
\[
 0<\mathcal C=\frac{up}{S}<\frac{k^2-r^2}{r}.
\]
More generally, for arbitrary $n\ge2$, whenever the displayed denominator
is positive,
\[
 \frac{\mathcal C}{r}<
 \frac{1-y}{1+y-\frac{n+2}{6}y\left(1-y\right)}.
\]
\end{proposition}
\begin{proof}
We justify the logarithmic convexity used here. In the product-integral
formula for $R$, each factor, for $a=\kappa_l^2>b$ and $0<v<1$, is
\[
 h_v\left(x\right)=\left(1-v\right)+\frac{v\left(a-b\right)}{a-x}.
\]
It is strictly logarithmically convex: writing $A=1-v>0$,
$B=v\left(a-b\right)>0$ and $t=a-x>0$ gives
\[
 h_vh_v''-\left(h_v'\right)^2=\frac{2AB}{t^3}+\frac{B^2}{t^4}>0.
\]
Products preserve logarithmic convexity, and so does integration
against the positive measure $\left(1-v\right)\dd v/2$, by H\"older's inequality.
Uniform convergence of the product and its first two derivatives on
$\left[0,b\right]$ gives logarithmic convexity for the infinite product-integral.
For the strict endpoint comparison below, one may keep one factor:
its strict logarithmic convexity for $0<v<1$ makes H\"older's bound
strict when the two $x$ arguments differ. Thus $\left(\log R\right)'$ strictly
increases. For the endpoint evaluation, put
\[
 \Sigma\left(x\right)=\sum_{l\ge2}\left(\kappa_l^2-x\right)^{-1}.
\]
The function $h\left(x\right)=p\left(\sqrt x\right)/\sqrt x$ solves
$4xh''+2\left(n+2\right)h'+h=0$. Since $h\left(b\right)=0$, its product gives
\[
 -2\Sigma\left(b\right)=\frac{h''\left(b\right)}{h'\left(b\right)}=-\frac{n+2}{2b},
 \qquad\Sigma\left(b\right)=\frac{n+2}{4b}.
\]
At the endpoint the product factors have value $1$ and their
$x$ derivatives are $v/\left(\kappa_l^2-b\right)$. Hence
\[
 R\left(b\right)=\frac14,\qquad
 R'\left(b\right)=\frac12\Sigma\left(b\right)\int_0^1v\left(1-v\right)\dd v
       =\frac{n+2}{48b}.
\]
Consequently, for every $x<b$,
\[
 \frac{R'\left(x\right)}{R\left(x\right)}<\frac{R'\left(b\right)}{R\left(b\right)}=\frac{n+2}{12b}.
\]
In the exact multiplier identity this gives
\[
 \frac r{\mathcal C}=1-\frac{2xT'}T
 =1+\frac{2x}{b-x}-\frac{2xR'}R
 >\frac{1+y-\frac{n+2}{6}y\left(1-y\right)}{1-y}.
\]
Taking reciprocals proves the second assertion where its denominator
is positive. If $n\le22$, then
\[
 1-\frac{n+2}{6}y\left(1-y\right)\ge1-4y\left(1-y\right)\ge0.
\]
Thus the strict preceding inequality gives
$r/\mathcal C>y/\left(1-y\right)>0$, including when $n=22$, $y=1/2$.
It follows that $\mathcal C/r<\left(1-y\right)/y$, as asserted.
\end{proof}

For $23\le n\le47$, we control the radial multiplier by the first
three coefficients of a positive power series. Normalizing their sum
to one turns the required estimate into an elementary moment inequality.

Put $x=r^2$, $b=k^2$, $y=x/b$, $\varepsilon=1-y$,
$T\left(x\right)=rt\left(r\right)$, $\rho=-t/t_r=\mathcal C$, and
\[
 R\left(x\right)=T\left(x\right)/\left(b-x\right),\qquad Q\left(y\right)=4R\left(by\right)=\sum_{j\ge0}\nu_jy^j.
\]
The product-integral formula for $R$ proves $\nu_j\ge0$, and the
endpoint identity $R\left(b\right)=1/4$ gives $\sum\nu_j=1$.
The exact multiplier formula is
\[
 \frac\rho r=
 \left(1+\frac{2y}{\varepsilon}-2y\frac{Q'\left(y\right)}{Q\left(y\right)}\right)^{-1}.
\]
Consequently $\rho/r<\varepsilon/y$ is equivalent to
\begin{equation}\label{Cupper47:moment-target}
 \sum_{j\ge0}\nu_j f_j\left(y\right)<0,\qquad
 f_j\left(y\right)=\left\{2j\left(1-y\right)-1\right\}y^j.
\end{equation}

\Needspace{6\baselineskip}
Three weighted coefficient sums suffice for the upper multiplier bound.

\begin{lemma}\label{lem:three-masses}
For every real parameter $23\le n\le47$ in the radial formulas,
\[
 \nu_0>\frac15,\qquad \nu_1>\frac3{25},\qquad
 \nu_2>\frac{89}{1000}.
\]
\end{lemma}
\begin{proof}
Write $d=n+2$, $h=n+4$ and $c=-X\left(k\right)>0$, and expand
$T\left(x\right)=T_0+T_1x+T_2x^2+O\left(x^3\right)$ at the center. Then
\[
 \nu_0=4T_0/b,\quad
 \nu_1=4\left(T_0/b+T_1\right),\quad
 \nu_2=4\left(T_0/b+T_1+bT_2\right),
\]
\[
 T_0=n\left(1+c\right),\quad T_1=\left(nc-2\right)/\left(2d\right),\quad
 T_2=\left\{nc\left(n+6\right)-8\right\}/\left(8d^2h\right).
\]
All three expressions increase with $c$. For real $n$, the radial
problem is the one-dimensional problem with weight $r^{n+1}$.
The trial function $v\left(r\right)=\left(1-r^2\right)^2$ gives
\[
 b<\frac{16\int_0^1r^{n+3}\left(1-r^2\right)^2\dd r}
          {\int_0^1r^{n+1}\left(1-r^2\right)^4\dd r}
   =\frac{\left(n+2\right)\left(n+10\right)}3.
\]
The inequality is strict because the trial function has a double zero
at $1$, whereas the first radial eigenfunction has a simple zero.
After setting $c=0$, each lower expression for $\nu_0,\nu_1,\nu_2$
decreases with $b$. Substituting this upper bound therefore gives
\[
 \nu_0>\frac{12n}{d\left(n+10\right)},\quad
 \nu_1>\frac{8\left(n-5\right)}{d\left(n+10\right)},\quad
 \nu_2>\frac{8\left(n-5\right)}{d\left(n+10\right)}-\frac{4\left(n+10\right)}{3dh}.
\]
For the first two desired comparisons, clear their positive
denominators. The resulting inequalities are
\[
 n^2-48n+20<0,\qquad -3n^2+164n-1060>0.
\]
The first convex polynomial is negative at both $23,47$, and the
second concave polynomial is positive there. The third comparison
reduces to
\[
 P_3\left(n\right)=-267n^3+15728n^2-122156n-901360>0.
\]
It is concave on $\left[23,47\right]$, since
$P_3''=31456-1602n<0$, and
$P_3\left(23\right)=1360575$, $P_3\left(47\right)=379719$.
\end{proof}

\Needspace{6\baselineskip}
We bound the remaining mass by elementary estimates for the coefficient tail.

\begin{lemma}\label{lem:mass-tail}
For the intervals in the next table, $f_j\left(y\right)<\beta_i$ holds for
all integers $j\ge0$ and all $y$ up to the indicated right endpoint.
The last row holds on the whole interval $0<y<1$.
\[
\begin{array}{c|c|c}
 b_i&r_i&\beta_i\\\hline
11/20&33/50&29/100\\
3/5&19/30&31/100\\
13/20&61/100&33/100\\
7/10&59/100&7/20\\
3/4&57/100&37/100\\
4/5&11/20&39/100\\
1&1/2&9/20
\end{array}
\]
\end{lemma}
\begin{proof}
Set $\ell=-\log y$ and $r\left(y\right)=\ell/\left[2\left(1-y\right)\right]$. Maximizing in the
continuous variable $u=j\ell\ge0$ gives
\[
 f_j\left(y\right)=\left\{u/r\left(y\right)-1\right\}e^{-u}
 \le\frac{e^{-1-r\left(y\right)}}{r\left(y\right)}.
\]
The last expression decreases as $r\left(y\right)$ increases, and $r\left(y\right)$
strictly decreases with $y$. At $y=b_i<1$, put
$\xi=\left(1-b_i\right)/\left(1+b_i\right)$. The elementary logarithm series gives
\[
 r\left(b_i\right)\ge\frac{1+\xi^2/3}{1+b_i}>r_i.
\]
Every last rational comparison is immediate from the displayed small
fractions. For $b_i=1$, use $r\left(y\right)>1/2$.
To check the exponential bounds without any decimal constant, let
$E_6\left(v\right)=\sum_{l=0}^6v^l/l!$. In each of the seven rows, direct
rational multiplication gives
\[
 \beta_i r_i E_6\left(1+r_i\right)>1+\frac3{1000}>1.
\]
Since $e^{1+r_i}>E_6\left(1+r_i\right)$, this proves all the bounds.
\end{proof}

\Needspace{6\baselineskip}
Combining the three masses with the tail bounds gives the required upper multiplier estimate.

\begin{theorem}\label{thm:multiplier-upper}
For every integer $2\le n\le47$ and every $0<r<k$,
\[
 0<\rho<\frac{k^2-r^2}{r}.
\]
\end{theorem}
\begin{proof}
Proposition~\ref{generalCupper:main} proves the result for
$n\le22$. For $23\le n\le47$ use Lemmas~\ref{lem:three-masses}
and~\ref{lem:mass-tail}.
If $0<y\le1/2$, then $f_1\le1/8$ and $f_j\le1/4$ for all
$j\ge2$. For the latter assertion the maximum is at $y=1/2$,
where $\left(j-1\right)2^{-j}\le1/4$; the case $j=2$ also has its maximum
there. Thus
\[
 \sum\nu_jf_j\le\frac14-\frac54\nu_0-\frac18\nu_1
                         <-\frac3{200}<0.
\]

On each remaining interval $\left[a_i,b_i\right]$, the first two polynomials
$f_1,f_2$ decrease for $y\ge1/2$. Hence
\[
 \sum\nu_jf_j
 \le \beta_i-\left(1+\beta_i\right)\nu_0
       -\left(\beta_i-f_1\left(a_i\right)\right)\nu_1
       -\left(\beta_i-f_2\left(a_i\right)\right)\nu_2.
\]
All the subtracted coefficients are positive. Substituting the three
mass lower bounds gives the following strictly negative upper bounds:
\[
\begin{array}{c|c|c|c}
{\left[a_i,b_i\right]}&f_1\left(a_i\right)&f_2\left(a_i\right)&\hbox{upper bound}\\\hline
{\left[1/2,11/20\right]}&0&1/4&-159/25000\\
{\left[11/20,3/5\right]}&-11/200&121/500&-463/250000\\
{\left[3/5,13/20\right]}&-3/25&27/125&-73/500000\\
{\left[13/20,7/10\right]}&-39/200&169/1000&-1509/1000000\\
{\left[7/10,3/4\right]}&-7/25&49/500&-97/15625\\
{\left[3/4,4/5\right]}&-3/8&0&-1451/100000\\
{\left[4/5,1\right)}&-12/25&-16/125&-1521/500000
\end{array}
\]
This proves
\eqref{Cupper47:moment-target} on all of $\left(0,1\right)$.
\end{proof}

\Needspace{6\baselineskip}
Theorem~\ref{thm:multiplier-upper} supplies the strict positivity of the coupling coefficient.

\begin{corollary}
\label{Cupper47:d}
In every $2\le n\le47$, $d>0$ throughout $0<r<k$, using only the
positive weights in the established expansion of $G$.
\end{corollary}
\begin{proof}
Write $G=\sum_{j\ge1}a_j\left(b^j-r^{2j}\right)$ with $a_j>0$. Then
\[
 d=2\sum_{j\ge1}a_j
       \left\{b^j-r^{2j}-j\rho r^{2j-1}\right\}>0,
\]
since $b^j-r^{2j}\ge j\left(b-r^2\right)r^{2j-2}>j\rho r^{2j-1}$.
\end{proof}

\subsection{Bessel gaps and exponential estimates}\label{sec:bessel-gaps}

For $m=n-1$, put $\Delta=k^2-q^2$. The following lower bounds
are sufficient for the scalar matrix proof:
\[
\begin{array}{c|c|c}
 n&\Delta/m& M_2=\left(\Delta-m-3\right)/12\\\hline
 11\le n\le27&>3/2&>\left(m-6\right)/24\\
 28\le n\le37&>11/8&>\left(3m-24\right)/96\\
 38\le n\le42&>49/36&>\left(13m-108\right)/432\\
 43\le n\le47&>4/3&>\left(m-9\right)/36
\end{array}
\]
Only the four right endpoints of the dimensional ranges need to be
estimated. Indeed the variational secant argument in Section~\ref{sec:dimension-threshold} for
$L\left(t\right)=j_{\sqrt t,1}^2$ gives that
\[
 \frac{j_{\nu+1,1}^2-j_{\nu,1}^2}{2\nu+1}
\]
is nonincreasing for $\nu\ge1/2$. Here $\nu=n/2-1$.

\[
\begin{array}{c|cc}
n&q&k\\\hline
27&\left(17.250,17.251\right)&\left(18.351,18.352\right)\\
37&\left(22.714,22.716\right)&\left(23.797,23.799\right)\\
42&\left(25.4171,25.4172\right)&\left(26.4936,26.4937\right)\\
47&\left(28.1043,28.1044\right)&\left(29.1755,29.1756\right)
\end{array}
\]
For $P_{d,N}\left(r\right)=\sum_{j=0}^Nt_j$, with
\[
 t_0=1,\qquad t_j=-\frac{r^2}{2j\left(d+2j-2\right)}t_{j-1},
\]
the alternating remainder gives $P_{d,61}<\psi_d<P_{d,60}$
at these arguments: $r_+^2<2\cdot61\left(d+120\right)$ ensures decreasing
tail magnitudes from the required index onward; monotonicity of the
early terms is not needed. The following are explicit finite rational
comparisons; terminating decimals denote exact rationals:
\[
\begin{array}{c|cc|c|c}
d&r_-&r_+&P_{d,61}\left(r_-\right)>\varepsilon&
                         P_{d,60}\left(r_+\right)<-\varepsilon\\\hline
27&17.250&17.251&10^{-7}&10^{-7}\\
29&18.351&18.352&10^{-8}&10^{-8}\\
37&22.714&22.716&6\cdot10^{-8}&6\cdot10^{-8}\\
39&23.797&23.799&3\cdot10^{-8}&3\cdot10^{-8}\\
42&25.4171&25.4172&10^{-9}&10^{-9}\\
44&26.4936&26.4937&8\cdot10^{-10}&8\cdot10^{-10}\\
47&28.1043&28.1044&10^{-10}&10^{-10}\\
49&29.1755&29.1756&10^{-10}&10^{-10}
\end{array}
\]
All entries mean precisely the displayed 61- or 62-term rational sum.
For an integer-only check, write $r=u/v$, put
$D_N=\prod_{j=1}^N2j\left(d+2j-2\right)$, and use
\[
 D_Nv^{2N}P_{d,N}\left(u/v\right)
 =\sum_{j=0}^N\left(-1\right)^j u^{2j}v^{2\left(N-j\right)}
                  \prod_{i=j+1}^N2i\left(d+2i-2\right).
\]

Here is an elementary identification of the roots as the first ones.
Let $\nu=d/2-1$, $a_\nu=\sqrt{\nu^2-1/4}$, and
$Y=r^{\nu+1/2}\psi_d\left(r\right)$. Then
\[
 Y''+V_\nu Y=0,\qquad V_\nu=1-a_\nu^2/r^2.
\]
The regular solution and its derivative are positive up to $a_\nu$.
On $\left[a_\nu,r_+\right]$ one has $V_\nu<\omega^2$, with
\[
 \begin{array}{c|cccc}
 d&27,29&37,39&42,44&47,49\\\hline
 \omega&3/4&2/3&5/8&3/5
 \end{array}
\]
These assertions are the rational comparisons
$\left(1-\omega^2\right)r_+^2<a_\nu^2$.
Duhamel comparison with frequency $\omega$ shows that the first
zero is after $a_\nu+\pi/\left(2\omega\right)$. Since $V_\nu<1$,
successive zeros have distance greater than $\pi$. A second zero
before $r_+$ would therefore exceed
\[
 a_\nu+\pi\left(1+\frac1{2\omega}\right)
 >\nu-\frac1{40}+\frac{314}{100}
                      \left(1+\frac1{2\omega}\right)>r_+.
\]
The inequalities hold for every row above; $a_\nu>\nu-1/40$
follows by squaring. This contradiction, together with the Taylor
sign change, proves each stated first-root bracket.

Finally the margins obtained from these brackets are
\[
\begin{array}{c|c}
n&k_-^2-q_+^2-\theta\left(n-1\right)\\\hline
27,\quad\theta=3/2&811/5000\\
37,\quad\theta=11/8&780553/10^6\\
42,\quad\theta=49/36&4006663/56250000\\
47,\quad\theta=4/3&5750267/\left(3\cdot10^8\right)
\end{array}
\]
They are strictly positive. The secant monotonicity proves all four
dimensional ranges in the first table.

The preceding estimates complete the Bessel-function part of the
determinant argument. The remaining scalar exponential inequalities
are proved by rational bounds and polynomial positivity. Each
coefficient comparison below applies on a whole interval.

\Needspace{6\baselineskip}
For the first dimension range above ten, the following exponential comparison closes the determinant estimate.

\begin{lemma}
\label{mgfblock:10-26}
Let $10\le m\le26$, put $M=\left(m-6\right)/24$, and for $0\le h\le1$ set
\[
 c=M\exp\left(-\frac{m+3}{8M}h\right).
\]
Then
\[
 1+4\left(1-h\right)c-\left(1-h\right)c^2h^2
 >\frac{m}{3\left\{3+\left(m-1\right)h/\left(2-h\right)\right\}}.
\]
\end{lemma}
\begin{proof}
Write $t=m-6$, so $4\le t\le20$, and abbreviate
\[
 R=\frac{\left(t+6\right)\left(2-h\right)}{3\left\{6+\left(t+2\right)h\right\}},\qquad
 H=\frac{t-3}{2\left(t+3\right)},\qquad z=\frac{3\left(t+9\right)}t h.
\]
Thus $c=\left(t/24\right)e^{-z}$. Since $0<c\le5/6$,
\[
 1+4\left(1-h\right)c-\left(1-h\right)c^2h^2
 =1+\left(1-h\right)c\left(4-ch^2\right)\ge1.
\]
If $h\ge H$, then $R\le1$; the desired inequality is strict for
$h<1$ by the positive extra term, and for $h=1$ because
$R=m/\left(3\left(m+2\right)\right)<1$.

It remains to consider $0\le h\le H<1/2$. In this range
$ch^2<5/24<1/4$, so
\[
 1+\left(1-h\right)c\left(4-ch^2\right)>1+\frac{15}4\left(1-h\right)c.
\]
Moreover
\[
 z\le\frac{3\left(t+9\right)\left(t-3\right)}{2t\left(t+3\right)}<2,
\]
since the positive numerator of the last difference is
$t^2-6t+81=\left(t-3\right)^2+72$.
For $0\le z<2$ one has the elementary lower bound
\[
 e^{-z}\ge\frac{2-z}{2+z}.
\]
Indeed $\log\left(\left(2+z\right)/\left(2-z\right)\right)-z$ vanishes at zero and has derivative
$z^2/\left(4-z^2\right)\ge0$.
Consequently it suffices to show
\[
 1+\frac{5t}{32}\left(1-h\right)
        \frac{2t-3\left(t+9\right)h}{2t+3\left(t+9\right)h}>R.
\]
All denominators are positive. Upon multiplication this becomes
$Q\left(t,h\right)>0$, where
\begin{align*}
 Q\left(t,h\right)={}&3\left\{6+\left(t+2\right)h\right\}\left\{2t+3\left(t+9\right)h\right\}\\
 &+\frac{15t}{32}\left(1-h\right)\left\{6+\left(t+2\right)h\right\}\left\{2t-3\left(t+9\right)h\right\}\\
 &-\left(t+6\right)\left(2-h\right)\left\{2t+3\left(t+9\right)h\right\}.
\end{align*}
Set $h=Hu$, $0\le u\le1$. The following cubic identity is exact:
\[
 \left(t+3\right)^3Q\left(t,Hu\right)=\sum_{j=0}^3B_j\left(t\right)\binom{3}{j}
                         u^j\left(1-u\right)^{3-j},
\]
where
\begin{align*}
 B_0\left(t\right)&=\frac{t\left(t+3\right)^3\left(13t+96\right)}8,\\
 B_1\left(t\right)&=\frac{\left(t+3\right)^2}{96}
       \left\{t^3\left(15t-52\right)+702t^2+10269t-7776\right\},\\
 B_2\left(t\right)&=\frac{t+3}{384}
       \left\{t^4\left(45t-131\right)+975t^3+39627t^2
                                 +136620t-93312\right\},\\
 B_3\left(t\right)&=\frac{15t\left(t+5\right)\left(t+6\right)\left(t+9\right)}{256}
                                  \left\{\left(t-3\right)^2+72\right\}.
\end{align*}
All four are positive for $t\ge4$. For the middle two, use
\[
 15t-52\ge8,\quad 10269t-7776\ge33300,
\]
\[
 45t-131\ge49,\quad 136620t-93312\ge453168.
\]
The Bernstein basis functions $\binom dj u^j\left(1-u\right)^{d-j}$ on $\left[0,1\right]$ are nonnegative and sum to one.
Thus $Q\left(t,Hu\right)>0$, completing the proof.
\end{proof}

\Needspace{6\baselineskip}
The next comparison treats the two middle-dimensional ranges with the same functional inequality.

\begin{lemma}
\label{mgfblockmid:statement}
Let either
\[
 27\le m\le36,\qquad M=\frac{m-8}{32},
 \qquad\text{or}\qquad
 37\le m\le41,\qquad M=\frac{13m-108}{432}.
\]
For $0\le h\le1$ put
\[
 c=M\exp\!\left(-\frac{m+3}{8M}h\right).
\]
Then
\begin{equation}
\label{mgfblockmid:goal}
 1+4\left(1-h\right)c-\left(1-h\right)c^2h^2
 >\frac{m\left(2-h\right)}{3\left\{6+\left(m-4\right)h\right\}}
 =\frac{m}{3\left\{3+\left(m-1\right)h/\left(2-h\right)\right\}}.
\end{equation}
The two ranges are continuous real intervals.
\end{lemma}

We record the short polynomial checks explicitly.  If $P\left(h\right)=\sum_{i=0}^d
p_i h^i$, its Bernstein coefficients on $\left[0,H\right]$ are
\begin{equation}
\label{mgfblockmid:bernstein}
 b_j\left(P;H\right)=\sum_{i=0}^j p_iH^i
                  \frac{\binom ji}{\binom di},\qquad 0\le j\le d.
\end{equation}
The binomial theorem gives
$P\left(Hu\right)=\sum_{j=0}^d b_j\left(P;H\right)\binom dj u^j\left(1-u\right)^{d-j}$.
Consequently, positive $b_j$ imply $P>0$ on $\left[0,H\right]$.

Define the following factored polynomials:
\begin{align*}
 D_1\left(h\right)={}&239\left(1-h\right)\left(19+44h\right)\left(6+23h\right)^2\left(19-15h\right)^4\\
 &\hspace{12mm}-12160\left(2-h\right)\left(6-4h\right)\left(19+15h\right)^4,\\
 P_1\left(h\right)={}&\left(-36+96h\right)\left(19+15h\right)^4
       +\frac{4541}{640}\left(6+23h\right)\left(1-h\right)\left(19-15h\right)^4,\\
 D_2\left(h\right)={}&\frac{2587}{7200}\left(1-h\right)
       \left(1+\frac{7938}{4849}h\right)
       \left(6+33h\right)^2\left(373-270h\right)^4\\
 &\hspace{12mm}-\left(2-h\right)\left(6-4h\right)\left(373+270h\right)^4,\\
 P_2\left(h\right)={}&\left(-56+136h\right)\left(373+270h\right)^4
       +\frac{74227}{7200}\left(6+33h\right)\left(1-h\right)\left(373-270h\right)^4.
\end{align*}
Here are strictly positive integer lower bounds obtained from
\eqref{mgfblockmid:bernstein}; the entries run from $j=0$ to $j=d$.
Thus each row is a direct rational polynomial identity and inequality.
\[
\begin{array}{c|c|c|c|l}
 P&H&d&S& \text{strict lower bounds for }b_j\left(P;H\right)/S\\ \hline
 D_1&9/22&8&10^8&
 22,66,107,136,144,131,99,55,2\\
 D_2&8/19&8&10^9&
 18,112,214,302,358,377,359,313,246\\
 P_2&7/17&6&10^9&
 113,210,146,65,60,203,558
\end{array}
\]
In particular, $D_1>0$ on $\left[0,9/22\right]$, $D_2>0$ on $\left[0,8/19\right]$,
and $P_2>0$ on $\left[0,7/17\right]$.

For completeness, the closer polynomial $P_1$ has a convenient
calculus proof of positivity on $\left[0,3/8\right]$.  Its expanded coefficients
in ascending order are
\[
 \begin{gathered}
 \frac{274065063}{320},\quad
 -\frac{2627387963}{640},\quad
 -\frac{18045957883}{640},\\
 \frac{9764885517}{64},\quad
 -\frac{1368225945}{16},\quad
 \frac{6761625525}{128},\quad
 -\frac{1057485375}{128}.
 \end{gathered}
\]
They give, throughout $0\le h\le3/8$,
\begin{align*}
 P_1'''\left(h\right)
 &>10^6\left\{915-2053h-992h^3\right\}\\
 &\ge10^6\left\{915-2053\frac38
                          -992\left(\frac38\right)^3\right\}>0.
\end{align*}
Direct substitution of the rational points gives
\begin{gather*}
 P_1'\left(0\right)<0,\quad P_1'\left(1/8\right)<0,\quad P_1\left(1/8\right)>181000,
 \quad P_1''\left(1/8\right)>44000000,\\
 P_1\left(3/16\right)>7300,\qquad \left|P_1'\left(3/16\right)\right|<530000.
\end{gather*}
Since $P_1'$ is convex, it is negative on $\left[0,1/8\right]$, where
$P_1\ge P_1\left(1/8\right)>0$.  On $\left[1/8,3/8\right]$ the increasing function $P_1''$
is greater than $44000000$.  Taylor's integral formula, in either
direction from $3/16$, therefore gives
\[
 P_1\left(h\right)>7300-530000\left|h-\frac3{16}\right|
              +22000000\left(h-\frac3{16}\right)^2
 \ge7300-\frac{530000^2}{88000000}>4100.
\]
This proves $P_1>0$ on the required interval.

\begin{proof}[Proof of Lemma~\ref{mgfblockmid:statement}]
Write
\[
 R_m\left(h\right)=\frac{m\left(2-h\right)}{3\left\{6+\left(m-4\right)h\right\}},\qquad
 a=\frac{m+3}{8M},\qquad c=Me^{-ah}.
\]
First suppose $27\le m\le36$.  The elementary maximum of
$h^2e^{-ah}$ on $h\ge0$, together with $e^2>7$, gives
\[
 ch^2\le\frac{4M}{a^2e^2}
 <\frac{4\left(7/8\right)}{\left(39/7\right)^2\,7}
 =\frac{49}{3042}<\frac1{60}.
\]
Consequently the left side of \eqref{mgfblockmid:goal} is greater than
$1+\kappa\left(1-h\right)c$ for $h<1$, where $\kappa=239/60$.
At $h=1$ the desired inequality holds directly.  Also $R_m\left(h\right)\le1$
when $h\ge\left(m-9\right)/\left(2m-6\right)$; thus $h\ge9/22$ is already settled.
For $0\le h\le9/22$ consider
\[
 \Phi_m\left(h\right)=1+\kappa\left(1-h\right)c-R_m\left(h\right).
\]
For fixed $h$, differentiation with respect to the real parameter $m$
gives
\[
 c_m=\frac1{32}e^{-4h}\left(1+z\right)e^{-z},\qquad
 z=\frac{44h}{m-8}.
\]
The function $\left(1+z\right)e^{-z}$ is decreasing on $z\ge0$.  Since $z$
decreases with $m$, $c_m$ increases with $m$.  Moreover
\[
 \partial_mR_m\left(h\right)
 =\frac{\left(2-h\right)\left(6-4h\right)}{3\left\{6+\left(m-4\right)h\right\}^2}
\]
decreases with $m$ on this interval.  It follows that
\[
 \partial_m\Phi_m\ge
 \frac{239}{1920}\left(1-h\right)\left(1+\frac{44h}{19}\right)e^{-120h/19}
 -\frac{\left(2-h\right)\left(6-4h\right)}{3\left(6+23h\right)^2}.
\]
For $0\le x<8$ we shall use the elementary inequality
\begin{equation}
\label{mgfblockmid:pade}
 e^{-x}\ge\left(\frac{1-x/8}{1+x/8}\right)^4.
\end{equation}
Indeed, if $u=x/8$, then
$\log\left(\left(1-u\right)/\left(1+u\right)\right)\le-2u$, as follows by differentiation.
Using \eqref{mgfblockmid:pade} in the preceding estimate and multiplying
by the positive denominator $36480\left(6+23h\right)^2\left(19+15h\right)^4$ proves
$\partial_m\Phi_m>0$, because its numerator is $D_1\left(h\right)>0$.
Hence $\Phi_m\ge\Phi_{27}$.
For $h\ge3/8$, $R_{27}\left(h\right)\le1$ and $\Phi_{27}>0$.
For $0\le h\le3/8$, applying \eqref{mgfblockmid:pade} once more shows
\[
 3\left(6+23h\right)\left(19+15h\right)^4\Phi_{27}\left(h\right)\ge P_1\left(h\right)>0.
\]
This completes the first range.

Now suppose $37\le m\le41$.  We have $M<1$, $a\ge2376/425>11/2$,
and therefore
\[
 ch^2\le\frac{4M}{a^2e^2}<\frac{16}{847}<\frac1{50}.
\]
Use $\kappa=199/50$ and the same definition of $\Phi_m$.
Here $R_m\left(h\right)\le1$ for $h\ge8/19$, and for smaller $h$
\[
 c_m=\frac{13}{432}e^{-54h/13}\left(1+z\right)e^{-z},\qquad
 z=\frac{7938h}{13\left(13m-108\right)}.
\]
Exactly the same monotonicity argument yields
\[
 \partial_m\Phi_m\ge
 \frac{199}{50}\frac{13}{432}\left(1-h\right)
      \left(1+\frac{7938h}{4849}\right)e^{-2160h/373}
 -\frac{\left(2-h\right)\left(6-4h\right)}{3\left(6+33h\right)^2}.
\]
Use \eqref{mgfblockmid:pade}, now with
$x=2160h/373$. Multiplication by
$3\left(6+33h\right)^2\left(373+270h\right)^4$ gives the lower bound
$D_2\left(h\right)>0$. Thus $\Phi_m\ge\Phi_{37}$.
The latter is positive if $h\ge7/17$, because $R_{37}\left(h\right)\le1$.
For $0\le h\le7/17$, the same exponential bound gives
\[
 3\left(6+33h\right)\left(373+270h\right)^4\Phi_{37}\left(h\right)\ge P_2\left(h\right)>0.
\]
This proves the second range and the lemma.
\end{proof}

\Needspace{6\baselineskip}
For the highest five dimensions, the determinant estimate reduces
to the following exponential comparison. Its parameter interval
covers $43\le n\le47$ after setting $m=n-1$.

\begin{lemma}\label{scalar-exp-matrix:comparison}
Let $42\le m\le46$, $u=m-9$, and put
\[
 \mu=\frac{u}{36},\qquad
 \alpha=\frac92+\frac{54}{u},\qquad
 c\left(h\right)=\mu e^{-\alpha h},\qquad
 \Theta_m\left(h\right)=\frac{m\left(2-h\right)}{3\left\{6+\left(m-4\right)h\right\}}.
\]
For every $0<h<1$,
\begin{equation}\label{scalar-exp-matrix:goal}
 1+4\left(1-h\right)c\left(h\right)-\left(1-h\right)c\left(h\right)^2h^2>\Theta_m\left(h\right).
\end{equation}
At $h=0$ the two sides agree.
\end{lemma}

\begin{proof}
Write the difference in \eqref{scalar-exp-matrix:goal} as $\mathcal F\left(u,h\right)$.
Since $\Theta_m\left(h\right)\le1$ exactly when
\[
 h\ge h_0\left(u\right):=\frac{u}{2\left(u+6\right)},
\]
and $c\left(h\right)h^2\le37/36<4$, this part of the interval is immediate.
It remains to consider
\[
 0<h\le h_0\left(u\right)\le\frac{37}{86}<\frac{13}{30}<\frac12.
\]
We first prove that $\mathcal F$ increases with $u$ on the larger rectangle
$33\le u\le37$, $0<h\le37/86$.

Use the elementary Pad\'e inequality
\begin{equation}\label{scalar-exp-matrix:pade}
 e^{-x}>\frac{S\left(x\right)}{N\left(x\right)},\qquad
 N\left(x\right)=1+\frac x2+\frac{x^2}{10}+\frac{x^3}{120},\quad
 S\left(x\right)=1-\frac x2+\frac{x^2}{10}-\frac{x^3}{120}
 \quad\left(0<x\le3\right).
\end{equation}
Indeed $S>0$ on this interval, and, for the unnormalized polynomials
$\widehat N=120N$, $\widehat S=120S$,
\[
 \frac{\mathrm d}{\dd x}\log\frac{\widehat N}{\widehat S}-1
   =\frac{x^6}{\widehat N\widehat S}>0.
\]
Their quotient is one at zero, proving \eqref{scalar-exp-matrix:pade}.

Direct differentiation gives
\begin{align*}
 \mathcal F_u={}&
 \frac{\left(1-h\right)\left(4-2ch^2\right)e^{-\alpha h}}{36}
       \left(1+\frac{54h}{u}\right)
 -\frac{\left(2-h\right)\left(6-4h\right)}{3\left\{6+\left(u+5\right)h\right\}^2}.
\end{align*}
Using $c\le37/36$, $\alpha\le135/22$, and $u\ge33$ shows that
\[
 \mathcal F_u\ge
 \frac{\left(1-h\right)\left(4-37h^2/18\right)\left(1+54h/37\right)e^{-135h/22}}{36}
 -\frac{\left(2-h\right)\left(6-4h\right)}{3\left(6+38h\right)^2}.
\]
All factors used in this comparison are positive. Apply
\eqref{scalar-exp-matrix:pade} with $x=135h/22<3$ and multiply by
$36N\left(135h/22\right)\left(6+38h\right)^2$. The resulting numerator is
\begin{align*}
 \mathcal Q\left(h\right)={}&
 \left(1-h\right)\left(4-37h^2/18\right)\left(1+54h/37\right)\left(6+38h\right)^2 S\left(135h/22\right)\\
 &\quad-12\left(2-h\right)\left(6-4h\right)N\left(135h/22\right).
\end{align*}
It has zero constant coefficient. The nine coefficients of
$\mathcal Q\left(h\right)/h$ have the following strict lower bounds, in increasing
order of degree:
\begin{equation*}
 \left(1174,\ 997,\ -13353,\ 11623,\ 19996,\ -37930,\ 9293,\ 16303,\ -8342\right).
\end{equation*}
These bounds follow by multiplication of the displayed cubic polynomials;
all entries are integers. Consequently
\begin{align*}
 \frac{\mathcal Q\left(h\right)}h>{}&
 R\left(h\right)+h^4\left(19996-37930h\right)
        +h^6\left(9293+16303h-8342h^2\right),\\
 R\left(h\right)={}&1174+997h-13353h^2+11623h^3.
\end{align*}
The last two terms are nonnegative for $0\le h\le13/30$ and
strictly positive when $h>0$.
On this same interval the cubic $R$ has Bernstein coefficients
\[
 1174,\qquad \frac{118621}{90},\qquad
 \frac{187867}{300},\qquad \frac{1198921}{27000},
\]
which are positive. Thus $\mathcal F_u>0$ throughout the stated rectangle.

It remains to prove $\mathcal F\left(33,h\right)>0$ for $0<h\le11/26$.
Here $\mu=11/12$, $\alpha=135/22$, and, with
$B=26/11$, $D=19/3$,
\[
 \Theta_{42}\left(h\right)-1=4\mu\frac{1-Bh}{1+Dh}.
\]
The map $z\mapsto4z-h^2z^2$ is increasing for the relevant positive $z$.
Hence \eqref{scalar-exp-matrix:pade} reduces the claim to
$\mathcal P\left(h\right)>0$, where
\begin{align*}
 \mathcal P\left(h\right)={}&
 \left(1-h\right)\left(1+Dh\right)\left\{S\left(\alpha h\right)N\left(\alpha h\right)
                    -\tfrac\mu4 h^2S\left(\alpha h\right)^2\right\}\\
 &\quad -\left(1-Bh\right)N\left(\alpha h\right)^2.
\end{align*}
Again the constant coefficient vanishes. The ten coefficients of
$\mathcal P\left(h\right)/h$, in increasing order, have strict lower bounds
\begin{equation*}
 \left(\begin{gathered}
 \frac{1560}{1000},-\frac{10886}{1000},
 \frac{3236}{1000},\frac{57078}{1000},\frac{36105}{1000},\\
 63,-79,78,-26,0
 \end{gathered}\right).
\end{equation*}
The pairings $h^5\left(63-79h\right)>0$ and $h^7\left(78-26h\right)>0$ for $h<1/2$ give
\[
 \frac{\mathcal P\left(h\right)}h>q\left(h\right):=
 \frac{1560-10886h+3236h^2+57078h^3+36105h^4}{1000}.
\]
This quartic is positive for every $h\ge0$. To see this without root
localization, note that $q''\left(h\right)>6$ there, whereas at $a=13/60$,
\[
 q\left(a\right)=\frac{11580751}{864000000}>\frac1{100},\qquad
 q'\left(a\right)=\frac{17057}{720000}<\frac1{30}.
\]
Strong convexity therefore gives, with $z=h-a$,
\[
 q\left(h\right)>\frac1{100}+q'\left(a\right)z+3z^2
 \ge\frac1{100}-\frac1{10800}>0.
\]
This proves the base case. Increasing $u$ proves the desired inequality
when $h\le h_0\left(33\right)$. If $h>h_0\left(33\right)$, the already established trivial
part gives $\mathcal F\left(33,h\right)>0$ as well; the same monotonicity applies
up to $37/86$. Thus the argument covers every required pair $\left(u,h\right)$.
\end{proof}

\section{The exact endpoint dimension threshold}\label{sec:dimension-threshold}

This section concerns the calibration used in this paper, rather than the
truth of the buckling inequality in arbitrary dimension. It proves an
obstruction in every dimension $n\ge48$. The endpoint test is positive in dimension
$47$; full-interval positivity is supplied by
Theorem~\ref{uniformcoupling:sharp}. The conclusion for all higher dimensions follows
from a variational concavity argument at the end of this section.

\begin{figure}[htbp]
\centering
\begin{tikzpicture}
\begin{axis}[width=.82\textwidth,height=5.2cm,
 xlabel={Dimension $n$},ylabel={$3\left(k^2-q^2\right)-4\left(n-1\right)$},
 xmin=40,xmax=55,xtick={40,43,47,48,51,55},
 grid=major,grid style={gray!18},axis line style={gray!60},
 tick label style={font=\small},label style={font=\small}]
\addplot[black,dashed,domain=40:55]{0};
\addplot[black,thick,mark=*,mark size=1.5pt] coordinates {(40,5.0419572933) (41,4.3466880162) (42,3.6469187258) (43,2.9428145242) (44,2.2345309624) (45,1.5222147869) (46,0.8060046138) (47,0.0860315381) (48,-0.6375803145) (49,-1.3647132861) (50,-2.0952557299) (51,-2.8291015946) (52,-3.5661500442) (53,-4.3063051105) (54,-5.0494753746) (55,-5.7955736743)};
\addplot[only marks,mark=square*,mark size=2.8pt] coordinates
 {(47,0.0860315381) (48,-0.6375803145)};
\end{axis}
\end{tikzpicture}
\caption{The endpoint coefficient changes sign between dimensions
$47$ and $48$. The plotted values illustrate the obstruction;
the sign assertions are proved below by rational root bounds and
variational concavity. They are not inferred from the plot.}
\label{fig:endpoint}
\end{figure}
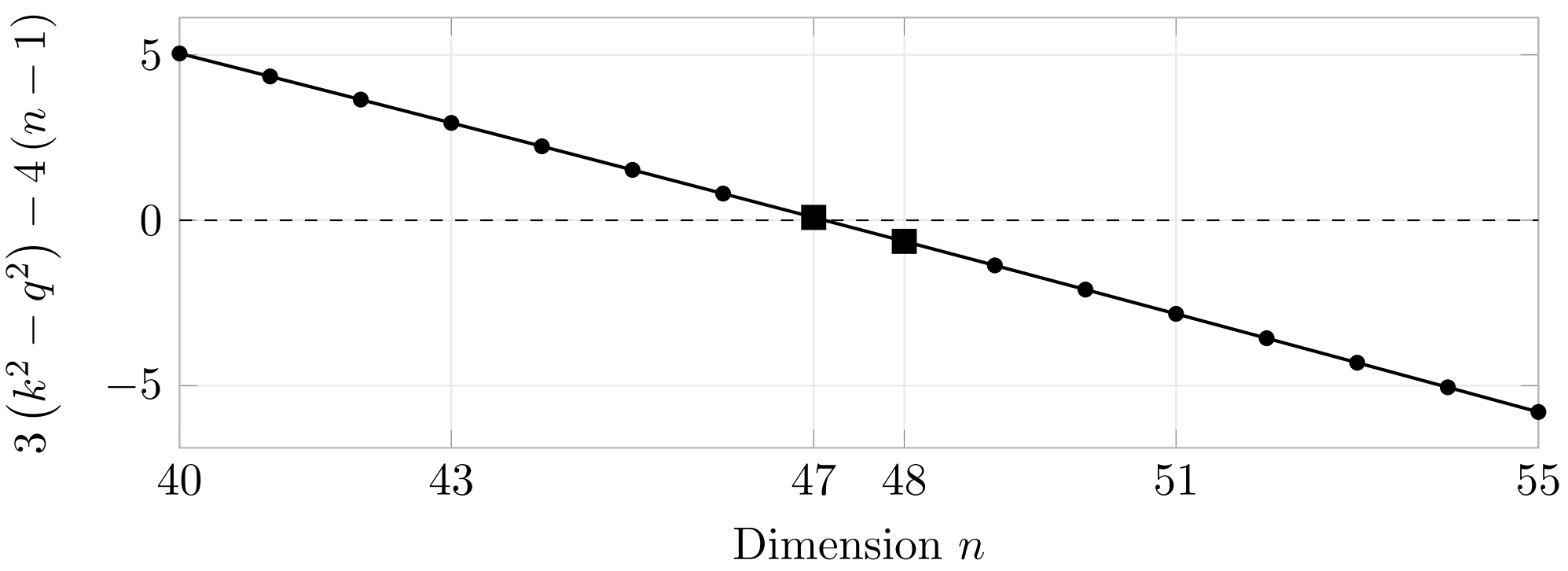

Figure~\ref{fig:endpoint} illustrates the endpoint sign. The proof
below supplies its exact value and propagates the sign to all larger
dimensions.

Retain the profiles of Section~\ref{sec:radial}, now in any fixed
dimension $n\ge2$, and put $m=n-1$. The defining ODEs and endpoint
expansions do not require positivity of the calibration.

\Needspace{6\baselineskip}
The leading term of the determinant gives a necessary condition for the calibration to remain positive semidefinite.

\begin{proposition}\label{n48:endpoint}
For every fixed dimension, as $\delta=k-r\downarrow0$,
\begin{equation}\label{n48:det-expansion}
 aE-\frac{d^2}{4}
 =\frac{K^2\delta^2}{9mk^2}
   \left\{3\left(k^2-q^2\right)-4m\right\}+O\left(\delta^3\right).
\end{equation}
Consequently $3\left(k^2-q^2\right)\ge4m$ is necessary for the positive
semidefiniteness of this calibration matrix throughout $0<r<k$.
\end{proposition}

\begin{proof}
The endpoint expansions in Section~\ref{sec:radial} follow solely
from the radial equations and hold in every fixed dimension. With
$h_{\mathfrak b}=1-\lambda-m/k^2$, they give
\[
 g=\beta-\frac{\beta h_{\mathfrak b}}6\delta^2+O\left(\delta^3\right),
 \qquad
 a=\frac{Kh_{\mathfrak b}}3\delta^2+O\left(\delta^3\right).
\]
To recover the slope of the correction, substitute
$G=\gamma\delta+O\left(\delta^2\right)$ into \eqref{eq:G-ode}.
The constant terms give $-\gamma=2\gamma-2K/k$, and hence
\[
 G=\frac{2K}{3k}\delta+O\left(\delta^2\right),\qquad
 d=\frac{2K}{3k}\delta+O\left(\delta^2\right),\qquad
 E=\frac Km+O\left(\delta\right).
\]
Consequently
\[
 aE-\frac{d^2}{4}
 =K^2\left\{\frac{h_{\mathfrak b}}{3m}-\frac1{9k^2}\right\}
       \delta^2+O\left(\delta^3\right).
\]
Substituting $h_{\mathfrak b}=1-q^2/k^2-m/k^2$ proves
\eqref{n48:det-expansion}. Its leading prefactor is positive,
so a negative bracket precludes positive semidefiniteness near $k$.
\end{proof}

\begin{figure}[htbp]
\centering
\begin{minipage}{.49\linewidth}
\centering
\begin{tikzpicture}
\begin{axis}[
 width=.94\linewidth,height=5.7cm,
 xmin=0,xmax=.1,ymin=-.015,ymax=.54,
 xlabel={$\varepsilon=1-r/k$},ylabel={$D_n\left(\varepsilon\right)$},
 xtick={0,.05,.1},
 grid=major,grid style={gray!15},axis line style={black!65},
 tick label style={font=\small},label style={font=\small},
 legend style={font=\small,draw=none,fill=white,at={(.03,.97)},anchor=north west}]
\addplot[gray!55,thin,forget plot] coordinates {(0,0) (.1,0)};
\addplot[black,thick] coordinates {
(0,0.0002078056) (0.002,0.01236846) (0.004,0.02451347) (0.006,0.03663831) (0.008,0.04873836)
(0.01,0.06080887) (0.012,0.07284507) (0.014,0.08484208) (0.016,0.09679501) (0.018,0.1086989)
(0.02,0.1205488) (0.022,0.1323398) (0.024,0.1440668) (0.026,0.155725) (0.028,0.1673095)
(0.03,0.1788153) (0.032,0.1902377) (0.034,0.2015719) (0.036,0.2128133) (0.038,0.2239573)
(0.04,0.2349994) (0.042,0.2459353) (0.044,0.2567606) (0.046,0.2674713) (0.048,0.2780634)
(0.05,0.288533) (0.052,0.2988764) (0.054,0.3090899) (0.056,0.3191703) (0.058,0.3291142)
(0.06,0.3389185) (0.062,0.3485804) (0.064,0.358097) (0.066,0.3674658) (0.068,0.3766844)
(0.07,0.3857505) (0.072,0.394662) (0.074,0.4034171) (0.076,0.4120141) (0.078,0.4204514)
(0.08,0.4287277) (0.082,0.4368416) (0.084,0.4447923) (0.086,0.4525788) (0.088,0.4602005)
(0.09,0.4676567) (0.092,0.4749471) (0.094,0.4820714) (0.096,0.4890296) (0.098,0.4958217)
(0.1,0.5024479)
};
\addlegendentry{$n=47$}
\addplot[black,thick,densely dashed] coordinates {
(0,-0.001507282) (0.002,0.01079243) (0.004,0.02307895) (0.006,0.03534754) (0.008,0.04759338)
(0.01,0.05981153) (0.012,0.07199697) (0.014,0.08414463) (0.016,0.09624937) (0.018,0.108306)
(0.02,0.1203094) (0.022,0.1322543) (0.024,0.1441355) (0.026,0.1559479) (0.028,0.1676863)
(0.03,0.1793456) (0.032,0.1909208) (0.034,0.2024069) (0.036,0.2137991) (0.038,0.2250925)
(0.04,0.2362824) (0.042,0.2473642) (0.044,0.2583336) (0.046,0.2691862) (0.048,0.2799178)
(0.05,0.2905243) (0.052,0.3010019) (0.054,0.3113468) (0.056,0.3215555) (0.058,0.3316246)
(0.06,0.3415508) (0.062,0.3513312) (0.064,0.3609627) (0.066,0.3704428) (0.068,0.3797689)
(0.07,0.3889388) (0.072,0.3979502) (0.074,0.4068013) (0.076,0.4154902) (0.078,0.4240154)
(0.08,0.4323754) (0.082,0.4405691) (0.084,0.4485954) (0.086,0.4564534) (0.088,0.4641423)
(0.09,0.4716618) (0.092,0.4790113) (0.094,0.4861907) (0.096,0.4931999) (0.098,0.500039)
(0.1,0.5067082)
};
\addlegendentry{$n=48$}
\end{axis}
\end{tikzpicture}
\end{minipage}\hfill
\begin{minipage}{.49\linewidth}
\centering
\begin{tikzpicture}
\begin{axis}[
 width=.94\linewidth,height=5.7cm,
 xmin=0,xmax=5,ymin=-1.8,ymax=3.7,
 xlabel={$10^4\varepsilon$},ylabel={$10^3D_n\left(\varepsilon\right)$},
 xtick={0,1,2,3,4,5},ytick={-1,0,1,2,3},
 grid=major,grid style={gray!15},axis line style={black!65},
 tick label style={font=\small},label style={font=\small},
 legend style={font=\small,draw=none,fill=white,at={(.03,.97)},anchor=north west}]
\addplot[gray!55,thin,forget plot] coordinates {(0,0) (5,0)};
\addplot[black,thick] coordinates {
(0,0.2078056) (0.1,0.2686405) (0.2,0.329475) (0.3,0.3903092) (0.4,0.4511432)
(0.5,0.5119769) (0.6,0.5728103) (0.7,0.6336434) (0.8,0.6944762) (0.9,0.7553087)
(1,0.816141) (1.1,0.876973) (1.2,0.9378047) (1.3,0.9986361) (1.4,1.059467)
(1.5,1.120298) (1.6,1.181129) (1.7,1.241959) (1.8,1.302789) (1.9,1.363618)
(2,1.424448) (2.1,1.485277) (2.2,1.546106) (2.3,1.606934) (2.4,1.667762)
(2.5,1.72859) (2.6,1.789418) (2.7,1.850245) (2.8,1.911072) (2.9,1.971899)
(3,2.032725) (3.1,2.093552) (3.2,2.154377) (3.3,2.215203) (3.4,2.276028)
(3.5,2.336853) (3.6,2.397678) (3.7,2.458502) (3.8,2.519326) (3.9,2.58015)
(4,2.640973) (4.1,2.701796) (4.2,2.762619) (4.3,2.823442) (4.4,2.884264)
(4.5,2.945086) (4.6,3.005908) (4.7,3.066729) (4.8,3.12755) (4.9,3.188371)
(5,3.249191)
};
\addlegendentry{$n=47$}
\addplot[black,thick,densely dashed] coordinates {
(0,-1.507282) (0.1,-1.445758) (0.2,-1.384235) (0.3,-1.322711) (0.4,-1.261188)
(0.5,-1.199666) (0.6,-1.138143) (0.7,-1.07662) (0.8,-1.015098) (0.9,-0.9535762)
(1,-0.8920544) (1.1,-0.8305329) (1.2,-0.7690115) (1.3,-0.7074904) (1.4,-0.6459695)
(1.5,-0.5844488) (1.6,-0.5229284) (1.7,-0.4614081) (1.8,-0.3998881) (1.9,-0.3383684)
(2,-0.2768488) (2.1,-0.2153295) (2.2,-0.1538104) (2.3,-0.09229151) (2.4,-0.03077287)
(2.5,0.03074555) (2.6,0.09226374) (2.7,0.1537817) (2.8,0.2152994) (2.9,0.2768169)
(3,0.3383342) (3.1,0.3998512) (3.2,0.461368) (3.3,0.5228846) (3.4,0.5844009)
(3.5,0.645917) (3.6,0.7074329) (3.7,0.7689485) (3.8,0.8304639) (3.9,0.8919791)
(4,0.953494) (4.1,1.015009) (4.2,1.076523) (4.3,1.138037) (4.4,1.199551)
(4.5,1.261065) (4.6,1.322579) (4.7,1.384092) (4.8,1.445605) (4.9,1.507118)
(5,1.56863)
};
\addlegendentry{$n=48$}
\end{axis}
\end{tikzpicture}
\end{minipage}
\caption{The calibration determinant near the outer endpoint. Here
$\varepsilon=1-r/k$ and
$D_n\left(\varepsilon\right)=\left(aE-d^2/4\right)/\left(K^2\varepsilon^2\right)$,
with all profiles evaluated at $r=k\left(1-\varepsilon\right)$.
The left panel shows $0\le\varepsilon\le0.1$; the right panel magnifies
$0\le\varepsilon\le5\times10^{-4}$ and rescales both axes as indicated.
The values at zero are the continuous limits supplied by
\eqref{n48:det-expansion}. The negative boundary interval in dimension
$48$ is visible only at the finer scale.
These curves are illustrations; the strict signs are established
by the determinant estimates and the exact endpoint analysis.}
\label{fig:determinant-boundary}
\end{figure}
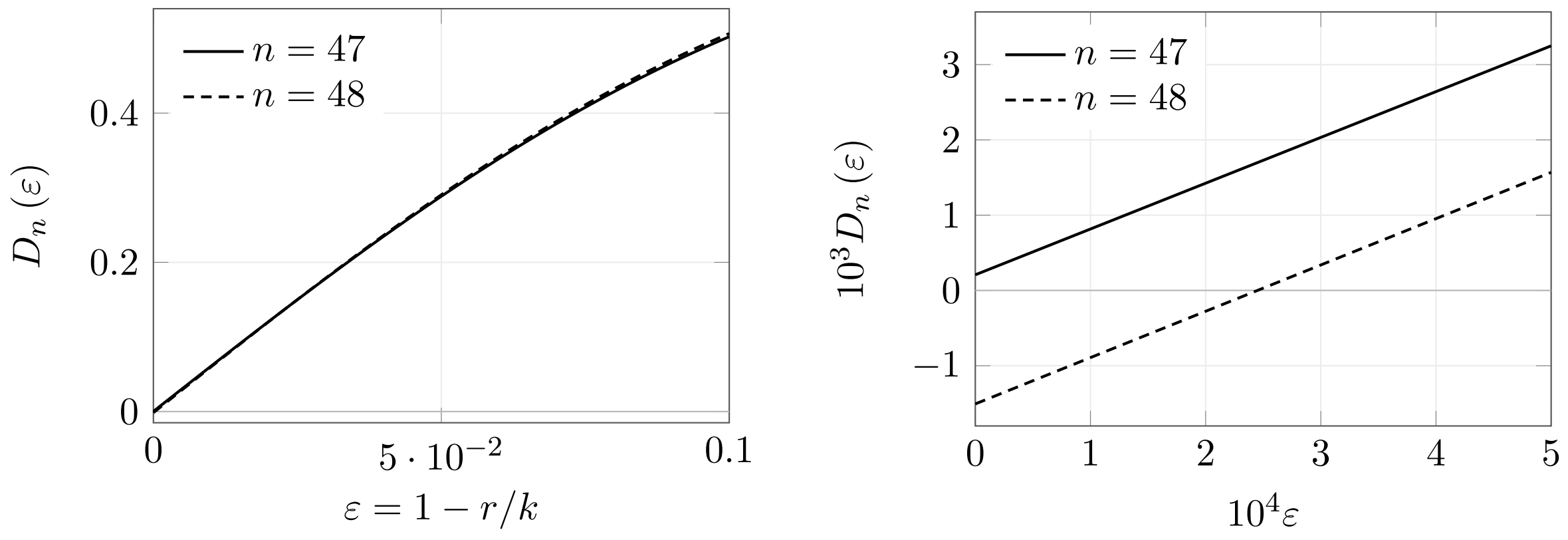

Figure~\ref{fig:determinant-boundary} resolves the endpoint behavior
of the determinant itself. Its limiting values are given exactly by
\eqref{n48:det-expansion}; the next step proves their signs without
using the plotted approximations.

We give finite Taylor comparisons and an elementary zero-counting
argument, so the root identification does not rely on a numerical root
finder. Define the rational polynomial
\[
 P_{n,N}\left(r\right)=\sum_{j=0}^{N}t_j,\qquad
 t_0=1,\quad
 t_j=-\frac{r^2}{2j\left(n+2j-2\right)}t_{j-1}.
\]
For all arguments below, $r<30$ and $n\ge47$. The magnitudes of the
terms decrease strictly from index $61$ onward, since their successive
ratio is at most $900/\left[124\left(47+122\right)\right]<1$. The alternating tail therefore
gives
\begin{equation}\label{n48:taylor-bounds}
 P_{n,61}\left(r\right)<\psi_n\left(r\right)<P_{n,60}\left(r\right).
\end{equation}
The following are finite rational inequalities. All terminating decimals
in the table denote exact rational numbers; the recurrence defining
$P_{n,N}$ specifies the arithmetic completely.
\[
\begin{array}{c|c|c|c}
 n&r&N&\text{finite rational comparison}\\\hline
48&28.640&61&P_{48,61}\left(28.640\right)>1/10^9\\
48&28.641&60&P_{48,60}\left(28.641\right)<-6/10^9\\
50&29.710&61&P_{50,61}\left(29.710\right)>26/10^{10}\\
50&29.711&60&P_{50,60}\left(29.711\right)<-25/10^{10}\\\hline
47&28.1043&61&P_{47,61}\left(28.1043\right)>1/10^{11}\\
47&28.1044&60&P_{47,60}\left(28.1044\right)<-1/10^{11}\\
49&29.1755&61&P_{49,61}\left(29.1755\right)>1/10^{11}\\
49&29.1756&60&P_{49,60}\left(29.1756\right)<-1/10^{11}
\end{array}
\]
For an integer-arithmetic verification of the table, let $r=v/w$ in positive integers and $N=60$ or
$61$, multiply $P_{n,N}\left(v/w\right)-\eta$ by the positive integer
\[
 M=w^{2N}\prod_{i=1}^{N}2i\left(n+2i-2\right)\cdot\operatorname{den}\left(\eta\right).
\]
The required sign is the sign of the integer
\[
 \operatorname{den}\left(\eta\right)
 \sum_{j=0}^{N}\left(-1\right)^jv^{2j}w^{2N-2j}
                  \prod_{i=j+1}^{N}2i\left(n+2i-2\right)
 -\operatorname{num}\left(\eta\right)w^{2N}
                  \prod_{i=1}^{N}2i\left(n+2i-2\right).
\]

To identify the roots in these sign-changing intervals as first roots,
put $\nu=n/2-1$, $a_\nu=\sqrt{\nu^2-1/4}$, and
\[
 Y\left(r\right)=r^{\nu+1/2}\psi_n\left(r\right),\qquad
 Y''+V_\nu\left(r\right)Y=0,\qquad
 V_\nu\left(r\right)=1-\frac{a_\nu^2}{r^2}.
\]
The regular expansion at zero gives $Y,Y'>0$ initially. As long as
$0<r\le a_\nu$ and $Y>0$, one has $Y''\ge0$. Hence
$Y,Y'>0$ throughout this interval, including at $a_\nu$.
Let $B$ be the upper argument in the corresponding pair of rows.
For each of the four orders $\nu=22.5,23,23.5,24$, direct rational
comparison gives
\[
 16B^2<25\left(\nu^2-1/4\right),\qquad
 V_\nu\left(r\right)<9/25\quad\left(a_\nu\le r\le B\right).
\]
Comparison with frequency $3/5$ shows that the first zero, if it lies
before $B$, is strictly greater than $a_\nu+5\pi/6$. One can see this
without invoking a zero theorem: before its first zero, the Duhamel
formula about $a_\nu$ is the positive homogeneous solution
\[
 Y\left(a_\nu\right)\cos\left(3t/5\right)+\tfrac53Y'\left(a_\nu\right)\sin\left(3t/5\right)
\]
plus the integral with kernel $\tfrac53\sin\left(3\left(t-v\right)/5\right)$ and forcing
$\left(9/25-V_\nu\right)Y$. Every term is nonnegative, and the homogeneous
solution is strictly positive, for $0\le t\le5\pi/6$.
Similarly, successive zeros have distance strictly greater than $\pi$:
after a zero, choose the sign to make the derivative positive and use
Duhamel comparison with frequency $1$; the forcing is
$\left(1-V_\nu\right)Y>0$ before the next zero. A next zero at distance at most
$\pi$ contradicts the positive sine formula and integral.
Consequently a second zero before $B$ would exceed
\[
 a_\nu+\frac{11\pi}{6}
 >\nu-\frac1{40}+\frac{11\left(3.14\right)}6
 >\nu+5.7316>B.
\]
The table and \eqref{n48:taylor-bounds} exhibit a zero before $B$,
whereas this argument permits at most one zero there. Thus the first-root
brackets are
\begin{align}
 n=48:\quad&28.640<q<28.641,\qquad29.710<k<29.711,
 \label{n48:root-brackets}\\
 n=47:\quad&28.1043<q<28.1044,\qquad29.1755<k<29.1756.
\notag\end{align}

\Needspace{6\baselineskip}
The brackets in dimension forty-eight force a negative leading determinant coefficient.

\begin{corollary}\label{cor:dimension48}
In dimension $48$ the present calibration matrix is not positive
semidefinite on the whole radial interval. In fact its radial--trace
principal determinant is negative for all $r<k$ sufficiently close to
$k$.
\end{corollary}
\begin{proof}
Here $m=47$. The brackets \eqref{n48:root-brackets} give the strict
rational inequality
\[
 3\left(k^2-q^2\right)-188
 <3\left\{\left(29.711\right)^2-\left(28.640\right)^2\right\}-188
 =-\frac{518237}{10^6}<0.
\]
Apply Proposition~\ref{n48:endpoint}; its leading factor
$K^2/\left(9mk^2\right)$ is strictly positive.
\end{proof}

For comparison, in dimension $47$ the equally explicit lower bound is
\[
 3\left(k^2-q^2\right)-184
 >3\left\{\left(29.1755\right)^2-\left(28.1044\right)^2\right\}-184
 =\frac{5750267}{10^8}>0.
\]
Thus the boundary principal determinant is positive sufficiently close
to $k$ in that dimension. This last statement is local: it neither
proves full-interval positivity in dimension $47$ nor establishes that
$48$ is the first dimension in which an interior obstruction occurs.

The conclusion for all higher dimensions follows from a variational
concavity argument, rather than from a numerical trend in Bessel zeros.
For real $\nu\ge1/2$, let $j_\nu$ denote the first positive zero of the
regular solution of order $\nu$. Put
\[
 L\left(t\right)=j_{\sqrt t}^{\,2},\qquad t\ge1/4.
\]
Then
\begin{equation}\label{n48:rayleigh}
 L\left(t\right)=\inf_{0\ne v\in H_0^1\left(0,1\right)}
 \frac{\displaystyle\int_0^1
     \left\{\left|v'\right|^2+\frac{t-1/4}{r^2}\left|v\right|^2\right\}\dd r}
      {\displaystyle\int_0^1\left|v\right|^2\dd r}.
\end{equation}
All forms have the same domain: the elementary Hardy inequality
$\int_0^1 v^2/r^2\dd r\le4\int_0^1\left|v'\right|^2\dd r$ holds on $H_0^1\left(0,1\right)$.
For a direct verification of \eqref{n48:rayleigh}, the positive function
\[
 w\left(r\right)=r^{\nu+1/2}\psi_{2\nu+2}\left(j_\nu r\right)
\]
satisfies the eigenvalue equation, vanishes at $1$, and belongs to
$H_0^1\left(0,1\right)$. For $v\in C_c^\infty\left(0,1\right)$, integration by parts gives
\[
 \int_0^1\left\{\left|v'\right|^2+
      \left(\frac{\nu^2-1/4}{r^2}-j_\nu^2\right)v^2\right\}\dd r
 =\int_0^1w^2\left|\left(\frac vw\right)'\right|^2\dd r\ge0.
\]
Approximation, the Hardy bound, and evaluation on $w$ prove
\eqref{n48:rayleigh}. For Hardy's inequality and the density statement, see \cite{Ziemer}.

Because the right-hand side of \eqref{n48:rayleigh} is an infimum of
affine functions of $t$, $L$ is concave. Thus its secant slope
\[
 R\left(\nu\right)=\frac{L\left(\left(\nu+1\right)^2\right)-L\left(\nu^2\right)}{\left(\nu+1\right)^2-\nu^2}
       =\frac{j_{\nu+1}^2-j_\nu^2}{2\nu+1}
\]
is nonincreasing in $\nu$. Here one uses the elementary concave-function
fact that a secant slope cannot increase when either endpoint of the
secant interval is moved to the right (while preserving their order).
This also handles overlapping intervals by first moving the right
endpoint, then the left endpoint; for disjoint intervals the usual
three-point secant inequalities apply.

For integer dimension $n\ge3$, $\nu=n/2-1$ and
$R\left(\nu\right)=\left(k^2-q^2\right)/\left(n-1\right)$. The rational comparisons already proved show
\[
 R\left(22.5\right)>\frac43,\qquad R\left(23\right)<\frac43.
\]
It follows that
\begin{equation*}
 \begin{cases}
 3\left(k^2-q^2\right)-4\left(n-1\right)>0,&3\le n\le47,\\
 3\left(k^2-q^2\right)-4\left(n-1\right)<0,&n\ge48.
 \end{cases}
\end{equation*}
In particular the present pointwise calibration fails in every dimension
$n\ge48$, already in a left neighborhood of its outer endpoint.
This endpoint argument alone gives only a local statement in dimensions
$3$ through $47$. For these dimensions, full-interval
positivity is supplied by Theorem~\ref{uniformcoupling:sharp}.

\appendix
\section{The planar coefficient estimates}\label{app:planar-coefficients}

The dimension-independent algebra in Sections~\ref{sec:radial}
and~\ref{sec:density} also applies when $n=2$. We therefore need only
the planar coefficient signs. The perspective map and the correction
profile are those of Tang and Zhang \cite{TZ}; the argument below
estimates the residue coefficients by comparing the first pole with
the remaining signed tail.

\subsection{Zero geometry and the first pole}

Put
\[
 q=j_{0,1},\quad k=j_{1,1},\quad \alpha=q/k,\quad
 \lambda=\alpha^2,\quad c=-J_0\left(k\right),\quad\sigma=1+c,\quad
 \beta=\alpha J_1\left(q\right)/c,\quad K=\beta^2.
\]
For $0<r<k$ write
\[
 X=J_0\left(r\right),\quad p=J_1\left(r\right),\quad u=X+c,\quad
 \phi=J_0\left(\alpha r\right),\quad S=p^2+up_r,\quad s=u/p,
\]
\[
 f\left(s\left(r\right)\right)=\phi/p,\quad z=f_s,\quad g=f-sz,\quad
 \rho=up/S.
\]
Interlacing gives $p,u,\phi>0$. The equations are
$X_r=-p$, $p_r=X-p/r$, and $\phi_{rr}+\phi_r/r+\lambda\phi=0$.
Also $u\left(r\right)=\int_r^kp\left(v\right)\dd v$.

Write $q_l=j_{0,l}$ and $\kappa_l=j_{1,l}$. The ratios
$q_l/\kappa_l$ strictly increase. Here is a proof that remains valid
in dimension two and does not use a higher-dimensional Sonin sign.
On the quarter-cycle from $q_l$ to $\kappa_l$, polar coordinates
$X=R\cos\theta$, $p=R\sin\theta$ give
\[
 \frac{\mathrm dr_l}{\mathrm d\varphi}=\frac{2r_l}{2r_l+\sin2\varphi},
 \qquad 0\le\varphi\le\pi/2.
\]
The trajectories strictly increase with $l$, while
\[
 \log\frac{\kappa_l}{q_l}
 =\int_0^{\pi/2}\frac{2\dd \varphi}{2r_l+\sin2\varphi}
\]
strictly decreases. This proves the ratio assertion.

Consecutive positive zeros of $J_1$ have distance greater than $\pi$:
apply Wirtinger's inequality between zeros to
$Y=\sqrt rJ_1\left(r\right)$, which solves
$Y''+\left(1-3/\left(4r^2\right)\right)Y=0$. The same variational equation on $\left(0,1\right)$
gives $k^2>3/4+\pi^2$. Sturm comparison on $\left[k,\kappa_2\right]$ then gives
\[
 \kappa_2-k<\frac{\pi k}{\sqrt{k^2-3/4}}<k.
\]
The quarter-phase speed is strictly greater than one, so
$k-q<\pi/2$. Thus
\[
 \kappa_2<2k<\frac\pi{1-\alpha},\qquad
 k<\alpha\kappa_2<q_2.
\]
For the last lower bound use
$\alpha\kappa_2>\kappa_2-\pi>k$; the upper bound is the
zero-ratio assertion just proved.

The following rational root bounds give the constants needed for the pole estimates.
For $\nu=0,1$ set
\[
 S_{\nu,N}\left(v\right)=\sum_{j=0}^N
 \frac{\left(-1\right)^j\left(v/2\right)^{2j+\nu}}{j!\left(j+\nu\right)!}.
\]
At the arguments below the alternating tail decreases from index two
onwards. Hence $S_{\nu,7}<J_\nu<S_{\nu,6}$, and direct rational
multiplication gives just four comparisons:
\[
\begin{array}{c|c|c|c|c}
 \nu&v_-&v_+&S_{\nu,7}\left(v_-\right)>&S_{\nu,6}\left(v_+\right)<\\\hline
 0&12/5&49/20&1/500&-1/50\\
 1&19/5&77/20&1/100&-1/200
\end{array}
\]
The series also gives $J_1\left(r\right)>r/2-r^3/16\ge0$ for
$0<r\le\sqrt8$. Zero spacing greater than $\pi$ therefore makes
the indicated $J_1$ zero the first one; interlacing identifies the
$J_0$ zero. Thus
\begin{equation*}
 \frac{12}5<q<\frac{49}{20},\qquad
 \frac{19}5<k<\frac{77}{20},\qquad
 \frac{48}{77}<\alpha<\frac{49}{76}.
\end{equation*}

The bounds $157/50<\pi<22/7$ used below were established
at the start of Section~\ref{sec:coefficients}.

For the pole estimates, define the auxiliary quantities
\[
 \begin{gathered}
 \mathcal B\left(r\right)=r\left(X^2+p^2\right)-Xp,\qquad
 \mathcal L\left(r\right)=K\mathcal B\left(r\right)-\alpha\mathcal B\left(\alpha r\right),\\
 Q\left(r\right)=\mathcal L\left(r\right)/r,\qquad G=K-r^2Q/p^2.
 \end{gathered}
\]
Differentiation and integration give
\begin{gather*}
 \mathcal B'=Xp/r,\qquad
 2\mathcal B\left(r\right)-p\left(r\right)^2/r=3\int_0^rp\left(v\right)^2/v^2\dd v,\\
 \mathcal B\left(q_{l+1}\right)-\mathcal B\left(q_l\right)
       =-\tfrac12\int_{q_l}^{q_{l+1}}X\left(v\right)^2/v^2\dd v.
\end{gather*}
The successive minima $\mathcal B\left(\kappa_l\right)$ increase strictly and
the successive maxima $\mathcal B\left(q_l\right)$ decrease strictly. In particular
\begin{equation*}
 0<\mathcal B\left(\kappa_l\right)\le\mathcal B\left(v\right)\le\mathcal B\left(q_l\right)
 \quad\left(v\ge q_l\right),\qquad
 K\mathcal B\left(k\right)=\alpha\mathcal B\left(q\right).
\end{equation*}
For every $l\ge2$, these identities imply
$Q\left(\kappa_l\right)>0$. Also
\[
 R_1:=\frac{\mathcal B\left(q\right)}{\mathcal B\left(k\right)}=K/\alpha>1.
\]

The energy ratio also determines the sign of the first coefficient at the center.
The quarter-phase formula for the energy ratio gives
\[
 \log R_1=\int_0^{\pi/2}
       \frac{2\sin2\varphi}{\left(2r_1+\sin2\varphi\right)^2}\dd \varphi.
\]
The measure $\sin2\varphi\dd \varphi$ has mass one and mean $\pi/4$
for $\sin2\varphi$. Jensen's inequality and $r_1<77/20$ give
\[
 \log R_1>\frac2{\left(77/10+\pi/4\right)^2}
 >\frac2{\left(297/35\right)^2}>\frac1{37}.
\]
Consequently $K>\left(38/37\right)\alpha$. The function
$38/37-2t+t^3$ decreases on $0<t\le49/76$, and
\[
 \frac{38}{37}-2\frac{49}{76}+\left(\frac{49}{76}\right)^3
       =\frac{90325}{16242112}>0.
\]
The center coefficient therefore has the sign
\begin{equation}\label{pc:new-a1}
 a_1:=\frac{K-2\lambda+\lambda^2}{4}>0.
\end{equation}

The higher coefficients require a quantitative lower bound for the first retained pole.
Put $\kappa=\kappa_2$, $x=\alpha\kappa$, $\delta=k-q$,
$A=3/4$, and $L=\pi/\left(1-\alpha\right)$. The root bounds and Sturm comparison give
\begin{gather*}
 \delta<29/20,\qquad k^2-A>\left(37/10\right)^2,\\
 \kappa<k+\frac{\pi k}{\sqrt{k^2-A}}
 <\frac{77}{20}+\frac{22}{7}\frac{38}{37}<\frac{71}{10},
 \qquad L>\frac{12089}{1450}>\frac{25}{3}.
\end{gather*}
Here $t/\sqrt{t^2-A}$ decreases, which explains the substitutions
in the second line. It follows that
\begin{equation}\label{pc:new-D}
 D:=\alpha^2-A/L^2>
 \left(\frac{48}{77}\right)^2-\frac{27}{2500}
 =\frac38+\frac{82959}{29645000}>\frac38.
\end{equation}

We next obtain two phase gaps without a $q_2$ or $\kappa_2$ root table.
The transformed $J_0$ potential on $\left[q,q_2\right]$ is
$1+1/\left(4r^2\right)<601/576$. Sturm comparison gives
\[
 q_2-q>\frac{24\pi}{\sqrt{601}}>3.
\]
Let $r_1,r_2$ be the quarter-phase trajectories from $q,k$ and
$q_2,\kappa$, respectively. Their differential equation has an
increasing right side, so $r_2-r_1\ge q_2-q>3$. Set
\[
 I_*:=\log\frac{q_2/\kappa}{q/k},\qquad
 E_*:=\log\frac{R_1}{R_2},\qquad
 R_2=\mathcal B\left(q_2\right)/\mathcal B\left(\kappa\right).
\]
Differentiating these ratios along the quarter-phase gives
\begin{align*}
 I_*&=\int_0^{\pi/2}
   \frac{4\left(r_2-r_1\right)}{\left(2r_1+\sin2\varphi\right)\left(2r_2+\sin2\varphi\right)}\dd \varphi,\\
 E_*&=\int_0^{\pi/2}2\sin2\varphi
   \left\{\frac1{\left(2r_1+\sin2\varphi\right)^2}
           -\frac1{\left(2r_2+\sin2\varphi\right)^2}\right\}\dd \varphi.
\end{align*}
Use $r_1<77/20$, $r_2<71/10$, and the gap $>3$.
Jensen's inequality in the first integral, with mean $2/\pi$, gives
\[
 I_*>\frac{18}{\left(77/10+2/3\right)\left(71/5+2/3\right)}
       =\frac{8100}{55973}>\frac17.
\]
In the second integral, use the probability measure
$\sin2\varphi\dd \varphi$ with mean $\pi/4$ and the decreasing
convex kernel $\left(a+b\right)/\left(a^2b^2\right)$. With
$D_k=297/35$, $D_c=1049/70$ one obtains
\[
 E_*>\frac{12\left(D_k+D_c\right)}{D_k^2D_c^2}
 >\frac{12\left(117/5\right)}{\left(17/2\right)^2\,15^2}>\frac1{58}.
\]
These kernels are convex because they are sums of products of
positive inverse powers of positive affine functions.

Let $\theta$ be the polar phase and
$\vartheta=3\pi/2-\theta\left(x\right)\in\left(0,\pi/2\right)$.
Writing $d=q_2-x$, the backwards phase equation gives
\[
 \vartheta\ge x\left(1-e^{-d/x}\right)\ge\frac d{1+d/x}
 >\frac x8>\frac12.
\]
Indeed $d=x\left(e^{I_*}-1\right)>x/7$, and
$x>\alpha\left(k+\pi\right)>\left(48/77\right)\left(34/5\right)>4$.
The increasing quarter-cycle lengths and the Sturm spacing bound give
\[
 \vartheta\le d<\delta\frac{\kappa-k}{k}
 <\frac{\delta\pi}{\sqrt{k^2-A}}<\frac{319}{259}.
\]
Thus $1<2\vartheta<638/259<\pi$. Concavity of sine shows
\begin{equation}\label{pc:new-sine}
 \sin2\vartheta>\frac{49}{80}.
\end{equation}
For a fully elementary check of the upper endpoint, the complementary
angle $t=\pi-638/259$ lies between $2/3$ and $1$, by the displayed
rational bounds on $\pi$. Hence
$\sin t>t-t^3/6>50/81>49/80$; at the lower endpoint use
$\sin1>5/6>49/80$.

The exact logarithmic derivative at $x$ is
\[
 u_*:=\frac{x\mathcal B'\left(x\right)}{\mathcal B\left(x\right)}
 =\frac{\sin2\vartheta/\left(2x\right)}{1-\sin2\vartheta/\left(2x\right)}.
\]
Since $x<\left(49/76\right)\left(71/10\right)$, \eqref{pc:new-sine} gives
\[
 \frac{\sin2\vartheta}{2x}>\frac{19}{284},\qquad
 u_*>\frac{19}{265}.
\]
Put $\sigma_*:=\mathcal B\left(x\right)/\left(R_1\mathcal B\left(\kappa\right)\right)$.
Then $0<\sigma_*<e^{-E_*}<58/59$, and the normalized first pole is exactly
\[
 U_2:=\frac{3Q\left(\kappa\right)+\alpha X\left(x\right)p\left(x\right)/\kappa}
             {KX\left(\kappa\right)^2}
 =3-\sigma_*\left(3-u_*\right).
\]
If $u_*\ge3$ this is at least three. Otherwise the previous bounds give
\begin{equation}\label{pc:new-U}
 U_2>\frac{19/265+3/58}{1+1/58}
       =\frac{1897}{15635}>\frac3{25}.
\end{equation}

\subsection{Tail control and the coefficient expansion}

We now compare this positive contribution with the possible negative tail.
For each $l\ge2$ put
\[
 U_l=\frac{3Q\left(\kappa_l\right)+\alpha X\left(\alpha\kappa_l\right)
                              p\left(\alpha\kappa_l\right)/\kappa_l}
                {KX\left(\kappa_l\right)^2}.
\]
If $\kappa_l<L$, zero spacing and the zero-ratio theorem put
$\kappa_{l-1}<\alpha\kappa_l<q_l$. Then $X,p$ have the same sign,
so $U_l>0$.
For $\kappa_l\ge L$, set $v=\alpha\kappa_l$,
$Y\left(v\right)=\sqrt vJ_1\left(v\right)$ and
$\mathcal E\left(v\right)=Y'^2+\left(1-A/v^2\right)Y^2$.
The Sonin identity $\mathcal E'=2AY^2/v^3\ge0$, integrated over a
nonempty interval on which $Y$ is not identically zero, gives
$\mathcal E\left(v\right)<\mathcal E\left(\kappa_l\right)=\kappa_lX\left(\kappa_l\right)^2$.
A direct substitution gives the exact numerator
\[
 \kappa_l KX\left(\kappa_l\right)^2U_l
 =3K\kappa_lX\left(\kappa_l\right)^2
 -3\alpha\left\{\mathcal E\left(v\right)+\frac{Y^2}{3v^2}
                                    -\frac{YY'}{3v}\right\}.
\]
Relative to $\mathcal E$, the two correction eigenvalues are bounded
above by the sum
$1/\left[3\left(v^2-A\right)\right]+1/\left[6\sqrt{v^2-A}\right]$.
Using $K/\alpha=R_1>1$ and \eqref{pc:new-D} therefore gives
\[
 U_l>-\frac1{2\sqrt D\,\kappa_l}-\frac1{D\kappa_l^2}.
\]
The zero spacing greater than $\pi$ implies, for $s>1$,
\[
 \sum_{\kappa_l\ge L}\kappa_l^{-s}
 \le L^{-s}+\frac{L^{1-s}}{\left(s-1\right)\pi}.
\]
Consequently, for every $j\ge2$,
\[
 \kappa_2^{2j}\sum_{l\ge2}\frac{U_l}{\kappa_l^{2j}}
 >U_2-\widetilde V,
\]
where
\[
 \widetilde V=
 \frac1{2\sqrt D}\left(\frac{\kappa_2}{L}\right)^4
       \left(\frac1L+\frac1{4\pi}\right)
 +\frac1D\left(\frac{\kappa_2}{L}\right)^4
       \left(\frac1{L^2}+\frac1{5\pi L}\right).
\]
Use $D>3/8$, $\sqrt D>3/5$, $L>25/3$,
$\kappa_2<71/10$, and $\pi>25/8$. These yield the single rational bound
\begin{equation*}
 \widetilde V<\left(\frac16+\frac{184}{3125}\right)
                  \left(\frac{213}{250}\right)^4<\frac3{25}.
\end{equation*}
The last positive difference is
$28105528377/24414062500000$. Together with \eqref{pc:new-U}, this
proves that the weighted sum of the $U_l$ is positive at every order
$j\ge2$, without requiring an individual sign at all far poles.

We obtain the coefficient signs by expressing the Taylor coefficients as sums of paired residues.
The points $0,\pm k$ are removable and $G$ is even meromorphic.
The ODE Taylor expansions, for $\delta=k-r$, begin
\begin{align*}
 p/c={}&\delta+\frac{\delta^2}{2k}
  +\frac{3/k^2-1}{6}\delta^3
  +\left(\frac1{2k^3}-\frac1{12k}\right)\delta^4+O\left(\delta^5\right),\\
 \phi/\left(\beta c\right)={}&\delta+\frac{\delta^2}{2k}
  +\frac{2/k^2-\lambda}{6}\delta^3
  +\left(\frac1{4k^3}-\frac{\lambda}{12k}\right)\delta^4+O\left(\delta^5\right).
\end{align*}
Substitution gives $G\left(k\right)=0$ and $G_r\left(k\right)=-2K/\left(3k\right)$.
At each further zero $\kappa$, one has
$p_r\left(\kappa\right)=X\left(\kappa\right)$,
$p_{rr}\left(\kappa\right)=-X\left(\kappa\right)/\kappa$, and
$\kappa\mathcal L'\left(\kappa\right)=\phi\left(\kappa\right)\phi_r\left(\kappa\right)$.
Thus the paired-pole calculation gives, for $j\ge2$,
\begin{equation}\label{pc:residue}
 a_j:=-\left[r^{2j}\right]G
 =2\sum_{l\ge2}
 \frac{KX\left(\kappa_l\right)^2U_l+2\left(j-2\right)Q\left(\kappa_l\right)}
      {X\left(\kappa_l\right)^2\kappa_l^{2j}}>0.
\end{equation}
To justify the contour passage, use the classical fixed-order asymptotic formula
\begin{equation*}
 J_\nu\left(\zeta\right)=\left(\frac2{\pi\zeta}\right)^{1/2}
 \left\{\cos\left(\zeta-\nu\pi/2-\pi/4\right)
       +O\left(\left|\zeta\right|^{-1}e^{\left|\operatorname{Im}\zeta\right|}\right)\right\},
 \quad\nu=0,1,\quad \left|\arg\zeta\right|\le\pi/2.
\end{equation*}
This fixed-order Hankel expansion and its complex remainder bound are given in \cite[10.17.3, 10.17.5--10.17.6
and \S10.17(iv)]{DLMF}.
On $\left|\zeta\right|=\left(N+3/4\right)\pi$, its leading $J_1$ cosine has modulus bounded
below by a fixed multiple of $e^{\left|\operatorname{Im}\zeta\right|}$.
For $\left|\operatorname{Im}\zeta\right|\ge1$ use $\sinh\left|\operatorname{Im}\zeta\right|$;
otherwise the real part differs from the radius by $O\left(N^{-1}\right)$.
Parity handles the left half-plane. Therefore
$\left|J_1\left(\zeta\right)\right|\ge C N^{-1/2}e^{\left|\operatorname{Im}\zeta\right|}$,
$\mathcal B\left(\zeta\right),\mathcal B\left(\alpha\zeta\right)=O\left(e^{2\left|\operatorname{Im}\zeta\right|}\right)$,
and $G\left(\zeta\right)=O\left(N^2\right)$ there.
The contour integral of $G\left(\zeta\right)/\zeta^{2j+1}$ consequently tends
to zero for $j\ge2$.

The nonzero zeros of $J_1$ are real and simple
\cite[Section 10.21(i)]{DLMF}. On the real zeros the same asymptotic gives
$X\left(\kappa_l\right)^2\asymp\kappa_l^{-1}$,
$\mathcal L\left(\kappa_l\right)=O\left(1\right)$ and
$\phi\left(\kappa_l\right)\phi_r\left(\kappa_l\right)=O\left(\kappa_l^{-1}\right)$.
The paired series is absolutely convergent by zero spacing. At
$\kappa$ the principal coefficients of $K-G$ are
\[
 A_\kappa=\frac{\kappa\mathcal L\left(\kappa\right)}{X\left(\kappa\right)^2},\qquad
 D_\kappa=\frac{2\mathcal L\left(\kappa\right)+\kappa\mathcal L'\left(\kappa\right)}
                    {X\left(\kappa\right)^2}.
\]
The paired residue of $G\left(\zeta\right)/\zeta^{2j+1}$ is
$2\left\{\left(2j-1\right)\mathcal L-\kappa\mathcal L'\right\}/
\left(\kappa^{2j+1}X\left(\kappa\right)^2\right)$, which is exactly the summand in
\eqref{pc:residue}. This proves the coefficient identity.

Finally \eqref{pc:new-a1}, \eqref{pc:residue} and analyticity before
$\kappa_2$ give
\begin{equation}\label{pc:G-coeff}
 G\left(r\right)=2\lambda-K-\sum_{j\ge1}a_jr^{2j},\qquad a_j>0,
 \qquad \left|r\right|<\kappa_2.
\end{equation}
Since $G\left(k\right)=0$, this also proves $G>0$ on $\left(0,k\right)$ and every
positive-order derivative of $G$ with respect to $r^2$ is negative.

\par\medskip\noindent
\textbf{Data availability.}
All quantities used in the proofs are specified in this article.
\par\medskip\noindent
\textbf{Conflict of interest.}
The authors declare that they have no conflicts of interest.


\begin{thebibliography}{99}
\bibitem{AH}
D.~R. Adams and L.~I. Hedberg,
\emph{Function Spaces and Potential Theory},
Grundlehren der mathematischen Wissenschaften 314, Springer, 1996.
\bibitem{AB}
M.~S. Ashbaugh and R.~D. Benguria,
On Rayleigh's conjecture for the clamped plate and its generalization
to three dimensions, Duke Math. J. \textbf{78} (1995), 1--17.
\bibitem{ABL}
M.~S. Ashbaugh, R.~D. Benguria and R.~S. Laugesen,
Inequalities for the first eigenvalues of the clamped plate and buckling
problems, in \emph{General Inequalities 7}, International Series of
Numerical Mathematics 123, Birkh\"auser, 1997, 95--110.
\bibitem{AshbaughBucur}
M.~S. Ashbaugh and D. Bucur,
On the isoperimetric inequality for the buckling of a clamped plate,
Z. Angew. Math. Phys. \textbf{54} (2003), 756--770.
\bibitem{AL}
M.~S. Ashbaugh and R.~S. Laugesen,
Fundamental tones and buckling loads of clamped plates,
Ann. Scuola Norm. Sup. Pisa Cl. Sci. (4) \textbf{23} (1996), 383--402.
\bibitem{BP}
J.~H. Bramble and L.~E. Payne,
Pointwise bounds in the first biharmonic boundary value problem,
J. Math. Phys. \textbf{42} (1963), 278--286.
\bibitem{buckeq:BDV}
L. Brasco, G. De Philippis and B. Velichkov,
Faber--Krahn inequalities in sharp quantitative form,
Duke Math. J. \textbf{164} (2015), 1777--1831.
\bibitem{buckeq:DNT}
C. De Coster, S. Nicaise and C. Troestler,
Nodal properties of eigenfunctions of a generalized buckling problem
on balls, Positivity \textbf{19} (2015), 843--875.
\bibitem{buckeq:HW}
L.-I. Hedberg and T.~H. Wolff,
Thin sets in nonlinear potential theory,
Ann. Inst. Fourier (Grenoble) \textbf{33} (1983), no. 4, 161--187.
\bibitem{LN}
J. Lamboley and M. Nahon,
Boundary regularity of a fourth order Alt--Caffarelli problem and
applications to the minimization of the critical buckling load,
preprint, arXiv:2512.18626, 2025.
\bibitem{Nadirashvili}
N.~S. Nadirashvili,
Rayleigh's conjecture on the principal frequency of the clamped plate,
Arch. Rational Mech. Anal. \textbf{129} (1995), 1--10.
\bibitem{DLMF}
F.~W.~J. Olver, A.~B. Olde Daalhuis, D.~W. Lozier et al., eds.,
\emph{NIST Digital Library of Mathematical Functions}, Chapters 5 and 10,
\url{https://dlmf.nist.gov/}.
\bibitem{PS}
G. P\'olya and G. Szeg\H{o},
\emph{Isoperimetric Inequalities in Mathematical Physics},
Annals of Mathematics Studies 27, Princeton University Press, 1951.
\bibitem{Stollenwerk}
K. Stollenwerk,
On the optimal domain for minimizing the buckling load of a clamped plate,
Z. Angew. Math. Phys. \textbf{74} (2023), article 12.
\bibitem{SW}
K. Stollenwerk and A. Wagner,
Optimality conditions for the buckling of a clamped plate,
J. Math. Anal. Appl. \textbf{432} (2015), 254--273.
\bibitem{Szego50}
G. Szeg\H{o}, On membranes and plates,
Proc. Nat. Acad. Sci. U.S.A. \textbf{36} (1950), 210--216.
\bibitem{Szego58}
G. Szeg\H{o}, Note to my paper ``On membranes and plates'',
Proc. Nat. Acad. Sci. U.S.A. \textbf{44} (1958), 314--316.
\bibitem{Talenti}
G. Talenti, On the first eigenvalue of the clamped plate,
Ann. Mat. Pura Appl. (4) \textbf{129} (1981), 265--280.
\bibitem{TZ}
Q. Tang and H. Zhang,
\emph{Proof of the Planar P\'olya--Szeg\H{o} Buckling Conjecture},
working preprint, 18 October 2026.
\href{https://github.com/QuanyuTang/polya-szego-planar-unconditional-palomar/blob/3b3417486b092363b411ed25ee2edb5e6f87f0ed/manuscript/PS_dim2_v5.tex}{Public manuscript, version 5, commit \texttt{3b3417486b09}}.
\bibitem{Teschl}
G. Teschl, \emph{Ordinary Differential Equations and Dynamical Systems},
Graduate Studies in Mathematics 140, American Mathematical Society, 2012.
\bibitem{Ziemer}
W.~P. Ziemer, \emph{Weakly Differentiable Functions},
Graduate Texts in Mathematics 120, Springer, 1989.
\end{thebibliography}
\end{document}